%% file: main.tex
\documentclass{article}
\usepackage{graphicx} 
\usepackage[english]{babel}
\usepackage{amsfonts,amsmath,amsthm,amssymb,mathtools,nicefrac}
\usepackage{bbold}
\usepackage{fullpage,longtable}
\usepackage{authblk}
\usepackage{enumitem}
\usepackage{todonotes}
\usepackage{hyperref}
\usepackage{bbm,stmaryrd}
\newcommand{\By}[2]{\overset{\mbox{\tiny{#1}}}{#2}}
\newcommand{\ByRef}[2]{   \By{\eqref{#1}}{#2} }
\newcommand{\eqBy}[1]{    \By{#1}{=} }

\newcommand{\eqByRef}[1]{ \ByRef{#1}{=} }

\newcommand{\geByRef}[1]{ \ByRef{#1}{\ge} }
\newcommand{\JUSTIFY}[1]{\fbox{\tiny{#1}}\quad}

\newtheorem{lem}{Lemma}[section]		
\newtheorem{fact}[lem]{Fact}		
		
\newtheorem{rem}[lem]{Remark}
\newtheorem{theo}[lem]{Theorem}
			
\newtheorem{cor}[lem]{Corollary}		
		
\newtheorem{prop}[lem]{Proposition}		
\newtheorem{quest}[lem]{Question}		

\newtheorem*{claimUni}{Claim~\ref*{cor:USTlocalglobal}.A}
\newtheorem*{claim:lem:TreeFunctionsContinuous1}{Claim~\ref*{lem:TreeFunctionsContinuous}.A}
\newtheorem*{claimFirstA}{Claim~\ref*{thm:factorFIgraphs}.A}		
\newtheorem*{claimFirstB}{Claim~\ref*{thm:factorFIgraphs}.B}
\newtheorem*{claimSecondA}{Claim~\ref*{thm:factorFI}.A}		
\newtheorem*{claimSecondB}{Claim~\ref*{thm:factorFI}.B}
\newtheorem*{claim:lem:degmeasmatch:A}{Claim~\ref*{lem:degmeasmatch}.A}
\newtheorem*{claim:lem:degmeasmatch:B}{Claim~\ref*{lem:degmeasmatch}.B}
\newtheorem*{claim:lem:FracIsoIFFProjFracIso:A}{Claim~\ref*{lem:FracIsoIFFProjFracIso}.A}
\newtheorem{defi}[lem]{Definition}

\newcommand{\Dirac}{\mathsf{Dirac}}
\newcommand{\Poisson}{\mathsf{Poi}}
\newcommand{\G}{\mathbb{G}}
\newcommand{\N}{\mathbb{N}}
\newcommand{\Z}{\mathbb{Z}}
\newcommand{\R}{\mathbb{R}}

\renewcommand{\P}{\mathbb{P}} 
\newcommand{\Probability}{\mathbf{P}} 
\newcommand{\Expectation}{\mathbf{E}}
\newcommand{\One}{\mathbbm{1}} 
\newcommand{\M}{\mathcal{M}}
\newcommand{\B}{\mathcal{B}}
\newcommand{\C}{\mathcal{C}}
\newcommand{\D}{\mathcal{D}}
\newcommand{\T}{\mathcal{T}}
\newcommand{\F}{\mathcal{F}}

\newcommand{\cX}{\mathcal{X}}
\newcommand{\fX}{\mathfrak{X}}
\newcommand{\fB}{\mathfrak{B}}
\newcommand{\fU}{\mathfrak{U}}
\newcommand{\diff}{\mathsf{d}}
\newcommand{\OneElementSpace}{\star}
\newcommand{\sm}{\setminus}
\newcommand{\eps}{\varepsilon}

\DeclareMathOperator{\esssup}{ess\, sup}
\DeclareMathOperator{\essinf}{ess\, inf}
\DeclareMathOperator{\degmax}{deg^{max}}
\DeclareMathOperator{\degmin}{deg^{min}}

\DeclareMathOperator{\dist}{dist}

\title{Graphon branching processes and fractional isomorphism}
\author[1]{Jan Hladký}
\author[2]{Eng Keat Hng}
\author[1]{Anna Margarethe Limbach}
\date{}
\affil[1]{Institute of Computer Science of the Czech Academy of Sciences\thanks{Research supported by Czech Science Foundation Project 21-21762X. With institutional support RVO:67985807.}}
\affil[2]{Extremal Combinatorics and Probability Group (ECOPRO), Institute for Basic Science, Daejeon, South Korea. \thanks{Supported by IBS-R029-C4. Most of the work was done while affiliated with the Institute of Computer Science of the Czech Academy of Sciences (institutional support RVO:67985807) and supported by Czech Science Foundation Project 21-21762X.}}

\begin{document}
\maketitle
\begin{abstract}
    In their study of the giant component in inhomogeneous random graphs, Bollobás, Janson, and Riordan introduced a class of branching processes parametrized by a possibly unbounded graphon. We prove that the tree structures underlying two such branching processes have the same distributions if and only if the corresponding graphons are fractionally isomorphic, a notion introduced by Grebík and Rocha.

    A different class of branching processes was introduced by Hladký, Nachmias, and Tran in relation to uniform spanning trees in finite graphs approximating a given connected graphon. We prove that that the tree structures of two such branching processes have the same distributions if and only if the corresponding graphons are fractionally isomorphic up to scalar multiple. Combined with a recent result of Archer and Shalev, this implies that if uniform spanning trees of two dense graphs have a similar local structure, they have a similar scaling limit.

    As a side result we give a characterization of fractional isomorphism for graphs as well as graphons in terms of their connected components.
\end{abstract}
\tableofcontents

\input{body}

\bibliographystyle{abbrv}
\addcontentsline{toc}{section}{References}
\bibliography{references}

\appendix

\section{Summary of notation}\label{app:notation}
\small

\input{summaryofnotation}

\section*{Contact details}
\noindent\begin{tabular}{ll} 
\emph{Postal address:} & Institute of Computer Science of the Czech Academy of Sciences \\ 
 & Pod Vodárenskou věží 2 \\ 
 & 182~00, Prague \\ 
 & Czechia\\
& \\
& Extremal Combinatorics and Probability Group \\
& Institute for Basic Science \\
& 55 Expo-ro, Yuseong-gu \\
& Daejeon 34126 \\
& South Korea \\
& \\
\emph{Email:}&\texttt{hladky@cs.cas.cz}, \texttt{hng@ibs.re.kr}, \texttt{limbach@cs.cas.cz}
\end{tabular}

\end{document}

%% file: body.tex
\section{Introduction}
\subsection{Bollobás--Janson--Riordan random rooted trees}\label{ssec:IntroBJR}
Branching processes related to graphs appear most prominently in the study of the giant component in the Erdős--Rényi random graph $\G(n,\frac{d}{n})$. While the original approach, starting with the seminal work of Erdős and Rényi~\cite{MR0125031}, used enumerative techniques, in the 1990s it was realized that the local structure of $\G(n,\frac{d}{n})$ can be approximated by the random rooted tree $\fX_d$ underlying the Galton--Watson branching process with offspring distribution $\Poisson(d)$. We recall the definition.
\begin{defi}[Galton--Watson random rooted tree with offspring distribution $\Poisson(d)$]
Let $d \ge 0$. The \emph{Galton--Watson random rooted tree} with offspring distribution $\Poisson(d)$ is the random (possibly infinite) rooted tree constructed as follows:
\begin{itemize}
    \item Start with a single root vertex.
    \item Independently for each vertex $v$, generate its number of offspring as a Poisson random variable with parameter $d$.
    \item Attach the corresponding number of child vertices to $v$.
\end{itemize}
Let $\fX_d$ denote the probability measure on the space of isomorphism classes of rooted trees induced by this procedure.
\end{defi}
The size of the giant component can be expressed in terms of $\fX_d$, too. As was kindly pointed out to us by Tomasz \L{}uczak, this idea first appeared in~\cite{MR1068492}. To make the statement precise, we need to introduce additional notation. First, for a random rooted (possibly infinite) tree $\mathfrak{Y}$ whose vertices have finite degrees and $k\in\N_0$, we write $(\mathfrak{Y})_{\restriction k}$\label{notation:restriction} for the $k$-ball of $\mathfrak{Y}$ around the root. Second, we recall the notion of local convergence (also known as the Benjamini--Schramm convergence), which is also relevant for our second main result in Section~\ref{ssec:IntroUST}.
\begin{defi}
Let $\mathfrak{G}$ be a probability distribution on isomorphism classes of rooted (possibly infinite) graphs with finite degrees. We say that a sequence $(F_n)_n$ of graphs \emph{converges locally} to $\mathfrak{G}$ if for every radius $r \in \N$, the distribution of the rooted $r$-neighborhood of a uniformly chosen random vertex in $F_n$ converges to the distribution of the $r$-neighborhood around the root in $\mathfrak{G}$. Here, the set of all possible rooted $r$-neighborhoods is equipped with the discrete topology.
\end{defi}
This convergence notion is in fact metrizable, see~\cite[Section 19.2]{MR3012035} for details. This in turn means that we can talk about the convergence of a sequence of random graphs $\mathbb{F}_1,\mathbb{F}_2,\ldots$ to $\mathfrak{G}$ in probability with respect to the local topology. The local structure of $\G(n,\frac{d}{n})$ can be then described as follows.

\begin{fact}\label{fact:localstructureGnp}
Let $d\ge 0$ be given. The sequence of random graphs $\left(\G(n,\frac{d}{n})\right)_{n=1}^\infty$ converges in probability with respect to the local topology to $\fX_d$. In particular, if $k\in \N$ and $T$ is a rooted graph whose vertices are all at distance at most $k$ from the root, then the probability that the $k$-ball around vertex~1 in $\G(n,\frac{d}{n})$ is isomorphic to $T$ converges to $\Probability[(\fX_d)_{\restriction k}\cong T]$, as $n\to \infty$. 

Furthermore, the order of the largest connected component in $\G(n,\frac{d}{n})$ is $(s+o_\mathrm{p}(1))n$, where $s$ is the survival probability of $\fX_d$. (The term $o_\mathrm{p}(1)$ converges to~0 in probability, as $n\to\infty$.) 
\end{fact}

A significant extension of the Erdős--Rényi random graph was introduced by Bollobás, Janson, and Riordan in~\cite{bollobas2007PhaseTransition}, with further important contributions to the subject in~\cite{MR2599196,MR2659281,MR2816939}. For our purposes, it suffices to introduce a slightly less general version of their model. Let $(X,\B)$\label{notation:XBmu} be a standard Borel space endowed with a Borel probability measure~$\mu$.
A \emph{kernel} is a symmetric measurable function $W:X^2\rightarrow [0,\infty)$. We sometimes call $X$ the \emph{ground space} of $W$. The \emph{degree} of $x\in X$ in $W$ is given by 
\begin{equation}\label{notation:degree}
\deg_W(x)=\int_{y\in X}W(x,y)\diff \mu(y)    \;.
\end{equation}
The \emph{maximum degree} and \emph{minimum degree} of $W$ are defined as the essential supremum and the essential infimum of the degree, \label{notation:minmaxdegre}$\degmin(W)=\essinf_x\deg_W(x)$, $\degmax(W)=\esssup_x\deg_W(x)$. When $\|W\|_1<\infty$, $\degmax(W)<\infty$, $\|W\|_\infty<\infty$, or $\|W\|_\infty\le 1$, we call $W$ an \emph{$L^1$-kernel},  \emph{bounded-degree kernel}, \emph{$L^\infty$-kernel}, or \emph{graphon}, respectively. 

The \emph{sparse inhomogeneous random graph} $\G(n,\frac{W}{n})$ is defined on the vertex set $[n]$ as follows. First, sample elements $x_1,\ldots,x_n\in X$ independently according to $\mu$. For each $\{i,j\}\in \binom{[n]}{2}$ independently, include $ij$ as an edge of $\G(n,\frac{W}{n})$ with probability $\min\{1,\frac{W(x_i,x_j)}{n}\}$. In particular, if $W\equiv d$ for some $d\in[0,\infty)$, we get the Erdős--Rényi random graph $\G(n,\frac{d}{n})$. To introduce a counterpart of Fact~\ref{fact:localstructureGnp} for $\G(n,\frac{W}{n})$, we define the random rooted tree $\fX_W$\label{notation:fXW} as follows. The type of the root is selected according to the distribution $\mu$, and each particle of type $x\in X$ has children whose number and types follow a Poisson point process on $X$ with intensity $W(x,y) \diff\mu(y)$. Observe that when $W$ is an $L^1$-kernel, the $k$-ball $(\fX_W)_{\restriction k}$ is finite almost surely for every $k\in\N$. The counterpart of Fact~\ref{fact:localstructureGnp} for $\G(n,\frac{W}{n})$, proved in~\cite{bollobas2007PhaseTransition}, then reads as follows.

\begin{fact}\label{fact:localstructureInhomo}
Let $W:X^2\rightarrow [0,\infty)$ be an $L^1$-kernel. The sequence of random graphs $\left(\G(n,\frac{W}{n})\right)_{n=1}^\infty$ converges to $\fX_W$ in probability with respect to the local topology. In particular, if $k\in \N$ and $T$ is a rooted graph whose vertices are all at distance at most $k$ from the root, then the probability that the $k$-ball around vertex~1 in $\G(n,\frac{W}{n})$ is isomorphic to $T$ converges to $\Probability[(\fX_W)_{\restriction k}\cong T]$, as $n\to \infty$. 

Furthermore, the order of the largest connected component in $\G(n,\frac{W}{n})$ is $(\gamma(W)+o_\mathrm{p}(1))n$, where $\gamma(W)$ is the survival probability of $\fX_W$.
\end{fact}

The survival probability $s$ of $\fX_d$ is known to be~0 for $d\le 1$, and to be the unique solution $s\in(0,1)$ of the equation $1-s=\exp(-ds)$ for $d>1$ (see e.g.~\cite[\S 10.4]{MR1885388}). In~\cite{bollobas2007PhaseTransition}, an inhomogeneous counterpart is obtained for $\fX_W$. However, the survival probability $\gamma(W)$ is concealed in a solution of a function-valued generalization of the above real-valued equation. Even for rather simple $L^1$-kernels $W$, this functional equation is not tractable. There are other quantities based on Bollobás--Janson--Riordan random rooted trees which are even more mysterious. 
Let us give an example concerning the random minimum spanning tree from~\cite{InhomoMST}. This example is stated using the framework of dense graph limits, as developed in~\cite{GraphLimitsLSz,GraphLimitsBSLSV}. Central to this framework are sequences of finite graphs of growing orders converging to a graphon in the so-called cut distance; we omit details and refer to Chapter~7 of~\cite{MR3012035}.
Suppose $(H_n)_n$ is a sequence of connected graphs converging (in the sense of dense graph convergence) to a connected\footnote{\label{foot:connected}A kernel $W$ is \emph{connected} if $\int_{A\times (X\sm A)}W>0$ for every $A\in\B$ with $\mu(A)\in (0,1)$. See Section~\ref{ssec:GraphonBasics} for details.} kernel $W$. On each edge of each $H_n$, put an independent weight chosen uniformly from $[0,1]$, and consider the minimum spanning tree $T_n$ of $H_n$. Then the total weights of $T_n$ converge in probability to a constant $\kappa(W)=\int_{t=0}^\infty \sum_{k=1}^\infty \frac{\Probability[|\fX_{t\cdot W}|=k]}{k}$, where $|\fX_{t\cdot W}|$ is the order of the random rooted tree $\fX_{t\cdot W}$. Expressing $\kappa(W)$ will be extremely challenging if not impossible. Note that the simplest case $W\equiv d\in [0,\infty)$ can be computed explicitly, see~\cite[\S 6.3]{InhomoMST}, and generalizes a famous result of Frieze~\cite{FRIEZE198547} about the random minimum spanning tree on complete graphs.

Our main result concerning Bollobás--Janson--Riordan random rooted trees characterizes bounded-degree kernels $U$ and $W$ for which $\fX_U$ and $\fX_W$ have the same distribution. We emphasize that by saying that the random rooted trees underlying two branching processes have the same distribution, we mean that they produce the same distribution on \emph{isomorphism classes} of rooted (and possibly infinite) trees.\footnote{\label{foot:equalityOfDistributions}Recall that two distributions $\mathfrak{T}$ and $\mathfrak{D}$ on rooted infinite trees are equal if for every $k\in \N$, $\mathfrak{T}_{\restriction k}$ and $\mathfrak{D}_{\restriction k}$ have the same distribution.} In particular, we disregard the labeling of the particles by their types.

Fractional isomorphism for graphs is a concept introduced in~1986 in~\cite{MR0843938} as a relaxation of graph isomorphism. Among several equivalent definitions, we give a rather intuitive one stated in terms of the color refinement algorithm, which we describe below. This definition naturally generalizes to the definition of fractional isomorphism of kernels that we primarily work with. In Section~\ref{sec:factorFI} we will see a different but equivalent definition of fractional isomorphism for graphs. The color refinement algorithm iteratively colors the vertices of a finite graph. It starts by giving each vertex the same color, and in each step it runs through the colors present and recolors the vertices of the graph in the following way. If all vertices with the current color see the same collection of colors (with multiplicities) on their respective neighbors, they are not recolored. If not, we recolor the vertices with the given color so that vertices get the same color if and only if they see the same colors (with multiplicities) on their respective neighbors. Furthermore, we use only previously unused colors for recoloring. The algorithm stops as soon as no vertices need to be recolored by the above rules. Two finite graphs of the same order are \emph{fractionally isomorphic} if the color refinement algorithm applied to their disjoint union yields a vertex coloring where each color class has the same number of vertices in each of the two graphs. In particular, any two $d$-regular $n$-vertex graphs are fractionally isomorphic, as the algorithm terminates immediately with a uniform coloring.

Remarkably, the above definition of fractional isomorphism is equivalent by results from~\cite{MR3829971,MR2668548,MR1297385} firstly to the existence of a doubly stochastic matrix $P$ so that for the adjacency matrices $A_H$ and $A_G$ of $H$ and $G$ respectively we have $PA_H=A_G P$, and secondly to having equal counts of each tree in $H$ and in $G$.

The above notions of fractional isomorphism were translated to graphons by Grebík and Rocha in~\cite{Grebik2022}, where they also showed that all these notions are equivalent in the graphon setting. In Section~\ref{ssec:FracIso}, we recall their definition of fractional isomorphism for graphons which corresponds to the notion for graphs given in terms of iterated degree sequences. In fact, we recall their definition in a slightly more general setting where we work with a generalized version of graphons that does not require boundedness and permits asymmetry; we call these objects \emph{akernels}.\footnote{The letter `a' in `\emph{a}kernel' indicates that asymmetry is permitted.} The notions of degree, maximum degree, and minimum degree carry over to this setting, but we need to careful about the order of the coordinates due to asymmetry. For example, in the defining formula~\eqref{notation:degree} the first coordinate is fixed and the second one is integrated over. We call a measurable function $W:X^2\rightarrow [0,\infty)$ a \emph{bounded-degree akernel} if $\degmax(W)<\infty$ and an \emph{$L^\infty$-akernel} if $\|W\|_\infty<\infty$. The definition of the random rooted tree $\fX_W$ is sensible for bounded-degree akernels as well. That is, the assumption of symmetry is not needed. In particular, it is the slices $W(x,\cdot)$ rather than $W(\cdot,x)$ that appear both in the definition of bounded degree akernels as well as in the definition of offspring generation in $\fX_W$. The generalization to the asymmetric setting may not seem important in the context of the original setting of Bollobás, Janson, and Riordan. However, this asymmetric version will play a key role in the proof of Theorem~\ref{thm:mainUST}.
\begin{theo}\label{thm:mainBJR}
    Suppose that $U$ and $W$ are bounded-degree akernels. Then $\fX_U$ and $\fX_W$ have the same distribution if and only if $U$ and $W$ are fractionally isomorphic.
\end{theo}
As we show in Proposition~\ref{prop:CorrespondenceTwoBranchings}, the $(\Leftarrow)$ direction of Theorem~\ref{thm:mainBJR} is easy. The $(\Rightarrow)$ direction, which is the main challenge of the theorem, is proved at the end of Section~\ref{ssec:FindingSepTree}.

\subsubsection{Applications}
Combining Theorem~\ref{thm:mainBJR} and Fact~\ref{fact:localstructureInhomo} immediately gives that if two random graph models $\G(n,\frac{U}{n})$ and $\G(n,\frac{W}{n})$ (where $U$ and $W$ are bounded-degree kernels) converge in probability in the local topology to the same limit, then $U$ and $W$ are fractionally isomorphic.

Another application involves a nice percolation result~\cite{MR2599196}. For a graph $F$ and for $p\in[0,1]$, we write $\mathbb{H}(F,p)$ for the random spanning subgraph of $F$ in which each edge is kept with probability $p$, independently of other choices. The main result of~\cite{MR2599196} says that if $a\ge 0$ and $F_1,F_2,\ldots$ is a sequence of graphs of growing orders converging in the cut distance topology\footnote{We need the notion of convergence in the cut distance only briefly. See~\cite[Section 8.2]{MR3012035} for details.} to a graphon $U$, then the random graphs $\mathbb{H}(F_n,\frac{a}{v(F_n)})$ converge in probability in the local topology to $\fX_{aU}$, as $n\to\infty$. Suppose in addition that $b\ge 0$ and $G_1,G_2,\ldots$ is a sequence of graphs of growing orders converging in the cut distance topology to a graphon $W$. Then another application of the same result combined with Theorem~\ref{thm:mainBJR} tell us that the sequences $\left(\mathbb{H}(F_n,\frac{a}{v(F_n)})\right)_n$ and $\left(\mathbb{H}(G_n,\frac{b}{v(G_n)})\right)_n$ converge in probability to the same local limit if and only if $aU$ is fractionally isomorphic to $bW$.

Let us give a different example, this time concerning the numerical parameter $\kappa(\cdot)$ defined above. Theorem~\ref{thm:mainBJR} allows us to say that the survival probabilities $\gamma(W_1)$ and $\gamma(W_2)$ are equal whenever $W_1$ and $W_2$ are fractionally isomorphic, even though we will typically be unable to determine their values. Since fractional isomorphism of $W_1$ and $W_2$ easily implies fractional isomorphism of $t\cdot W_1$ and $t\cdot W_2$ for every $t\ge 0$, we also get $\kappa(W_1)=\kappa(W_2)$. This is in particular useful for $d$-regular kernels. Recall that a kernel $W$ is $\emph{$d$-regular}$ (for $d\in[0,\infty)$) if for $\mu$-almost every $x\in X$ we have $\deg_W(x)=d$. The class of all $d$-regular kernels may seem quite complicated, but the above implies that any parameter derived from $\fX_W$ is constant on $d$-regular kernels, and in particular is equal to its value for $\fX_d$, where it may be tractable computationally, as was the case with $\gamma(\fX_d)$ and $\kappa(\fX_d)$. Note that in this example we only used the easy $(\Leftarrow)$ direction of Theorem~\ref{thm:mainBJR}.

\subsection{Local structure of uniform spanning trees}\label{ssec:IntroUST}
Our second main result is related to random rooted trees arising from the study of the uniform spanning tree. Recall that the \emph{uniform spanning tree} of a finite connected graph $G$ is the uniform measure on spanning trees of $G$. While much of the study of uniform spanning trees is concerned with large but sparse graphs with lattice structure, here we are concerned with dense graphs.
The first result that relates to the uniform spanning tree of dense graphs is by Kolchin~\cite{Kolchin1977} and Grimmett~\cite{Gri:RandomTree}. The result says that the uniform spanning trees of $K_n$ converge, as $n\to\infty$, locally to the random rooted tree underlying the Galton--Watson branching process with offspring distribution $\Poisson(1)$, conditioned on survival. It is well-known that such a random rooted tree can be constructed by taking a one-way infinite path rooted at its endvertex and attaching $\fX_1$ to each vertex of the aforementioned one-way infinite path.\footnote{Let us include details. Section~5 of \cite{MR2908619} contains a construction of a modified Galton--Watson branching process by introducing a size-biased version of the offspring distribution. Since we started with a $\Poisson(1)$ Galton--Watson branching process, the modified Galton--Watson branching process has offspring distribution in which each particle produces $k$ offspring with probability $k\Probability[\Poisson(1)=k]$. It is easy to see that the size-biased $\Poisson(1)$ distribution is equal to $1+\Poisson(1)$, and hence corresponds to the above construction with a one-way infinite path. Theorem 22.2 in~\cite{MR2908619} shows that the modified Galton--Watson random rooted tree is indeed the $\Poisson(1)$ Galton--Watson random rooted tree conditioned on survival.} This corresponds to the random rooted tree $\fU_W$ defined below, when we take $W\equiv 1$.

The above result was generalized from complete graphs to sequences converging to a graphon by Hladký, Nachmias, and Tran~\cite{MR3876899}. Let us introduce their definition in the more general setting of $L^1$-kernels. We say that an $L^1$-kernel $W$ is \emph{nondegenerate} if $\deg_W(x)>0$ for almost every $x\in X$.
\begin{defi}\label{def:HNTbranching}
Given a nondegenerate $L^1$-kernel $W$, we define a multitype Galton--Watson random rooted tree $\fU_W$ with type space $\{(\mathsf{anc},x),(\mathsf{oth},x)\;:\:x\in X\}$ as follows. Here, ``$\mathsf{anc}$'' stands for ``ancestral'' and ``$\mathsf{oth}$'' stands for ``other''.
	\begin{enumerate}[label=(\roman*)]
		\item The root has type $(\mathsf{anc},x)$, where the distribution of $x$ is $\mu$.
		\item\label{en:HNTdesc} If a particle has type $(\mathsf{oth},x)$, then its children are $\{(\mathsf{oth},x_1),\ldots,(\mathsf{oth},x_k)\}$ where $\{x_1,\ldots,x_k\}$ is given by a Poisson point process on $X$ with intensity $\frac{W(x,y)}{\deg_W(y)}$ at $y\in X$. In particular, the number of children of this particle has distribution $\Poisson\left(\int_y\frac{W(x,y)}{\deg_W(y)}\right)$.
		\item\label{HNT3} If a particle has type $(\mathsf{anc},x)$, then its children are $\{(\mathsf{anc},x_0),(\mathsf{oth},x_1),\ldots,(\mathsf{oth},x_k))\}$ where $\{x_1,\ldots,x_k\}$ is given by a Poisson point process on $X$ with the same intensity as above and $x_0$ is an independent element which is distributed according to the probability measure on $X$ that has density at each $y\in X$ equal to $\frac{W(x,y)}{\deg_W(x)}$. In particular, the number of children of this particle has distribution $1+\Poisson\left(\int_y\frac{W(x,y)}{\deg_W(y)}\right)$.
	\end{enumerate}
 (Recall also that while types are used to construct $\fU_W$, the labels of the types are eventually removed; see text around Footnote~\ref{foot:equalityOfDistributions}.)
\end{defi}
\begin{figure}\centering
	\includegraphics[scale=0.6]{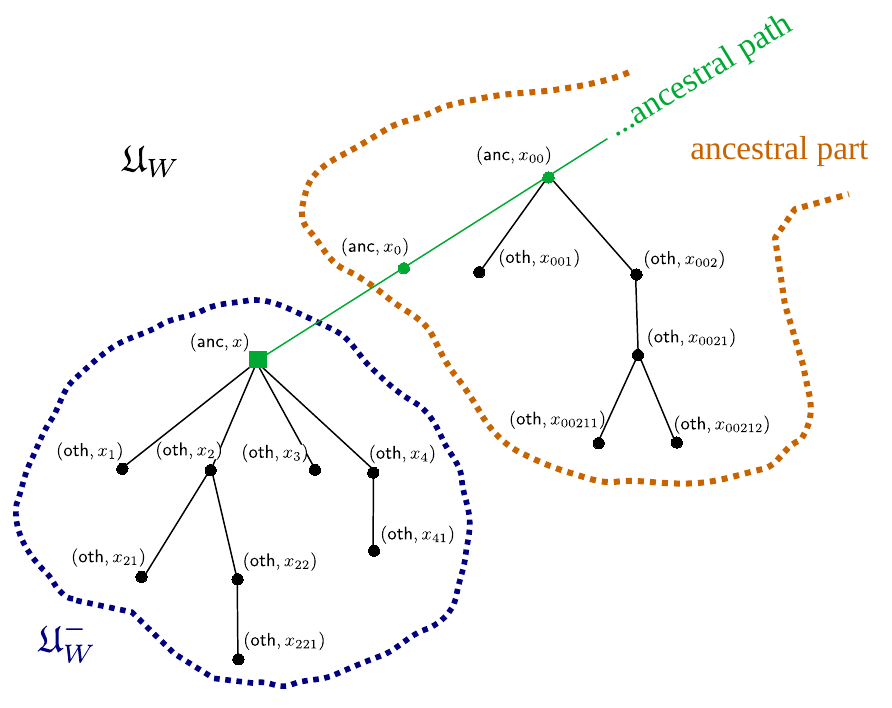}
	\caption{The random rooted tree $\fU_W$ from Definition~\ref{def:HNTbranching}. Let \label{notation:HNTminus}$\fU^-_W$ be obtained from $\fU_W$ by removing the unique $\mathsf{anc}$-vertex neighboring the root and the whole subtree appended to it. The complement of $\fU_W^-$ is called the \emph{ancestral part} of $\fU_W$. The one-way infinite path rooted at the root of the tree consisting of $\mathsf{anc}$-vertices is called the \emph{ancestral path}.}
	\label{fig:HNTbranching}
\end{figure}

In Figure~\ref{fig:HNTbranching} we give an example of the random rooted tree $\fU_W$. In that figure we also introduce the random rooted tree $\fU_W^-$ and the ancestral part of $\fU_W$. The main result of~\cite{MR3876899} then reads as follows.

\begin{theo}\label{thm:HNT}
Let $W:X^2\rightarrow [0,\infty)$ be an arbitrary nondegenerate graphon. Suppose that $(G_n)_n$ is a sequence of connected graphs of growing orders that converges to $W$ in the cut distance. Then the uniform spanning trees of $G_n$ converge in probability with respect to the local topology to $\fU_W$, as $n\to \infty$.
In particular, if $r\in \N$ and $T$ is an arbitrary rooted graph whose vertices are all at distance at most $r$ from the root, then the probability that the $r$-ball around a randomly selected vertex in the uniform spanning tree of $G_n$ is isomorphic to $T$ converges to $\P[(\fU_W)_{\restriction r}\cong T]$, as $n\to \infty$.
\end{theo}

Like in Theorem~\ref{thm:mainBJR}, it can be checked that $\fU_U$ and $\fU_W$ have the same distribution if $U$ and $W$ are fractionally isomorphic. There is however another operation which is easily seen to preserve the distribution of $\fU_W$, namely multiplication by a scalar. Our main result says that there are no further ways of creating a connected kernel $Z$ with the same distribution as $\fU_Z$. This leads to the following definition.
\begin{defi}\label{def:projfraciso}
We say that two kernels $U$ and $W$ are \emph{projectively fractionally isomorphic} if there exists a constant $t>0$ (called the \emph{projective constant}) so that $U$ and $tW$ are fractionally isomorphic. In the non-degenerate case $U,W\neq0$, this is equivalent to $\frac{U}{\|U\|_1}$ and $\frac{W}{\|W\|_1}$ being fractionally isomorphic. 
\end{defi}
We state the theorem only for connected kernels. A full characterization including disconnected kernels is given as Theorem~\ref{thm:mainUSTGeneral} in Section~\ref{sec:HNT}.
\begin{theo}\label{thm:mainUST}
    Suppose that $U$ and $W$ are connected $L^\infty$-kernels with positive minimum degrees. Then $\fU_U$ and $\fU_W$ have the same distribution if and only if $U$ and $W$ are projectively fractionally isomorphic.
\end{theo}

\subsubsection{An application: local limit versus scaling limit}
Theorem~\ref{thm:mainUST} has an interesting corollary which connects the local limit and the scaling limit of uniform spanning trees on dense graphs. We briefly recall the latter concept. Aldous~\cite{MR1085326} famously constructed a certain random metric space which he called `the continuum random tree' $\mathcal{T}$. He showed that $\mathcal{T}$ is the Gromov--Hausdorff--Prokhorov limit of the uniform spanning tree on $K_n$ when each edge is given length $n^{-1/2}$, as $n\to\infty$. This type of convergence is called the `scaling limit'. Since then, $\mathcal{T}$ has been shown to be the scaling limit of many other sequences of uniform spanning trees. For dense graphs, a recent result of Archer and Shalev~\cite{GlobalUSTGraphon} asserts that in the setting of Theorem~\ref{thm:HNT} (with an additional assumption that $W$ is connected), the scaling limit of the uniform spanning trees is $c_W\cdot \mathcal{T}$, where the rescaling constant $c_W$ is defined by (see~\cite[Equation~(1)]{GlobalUSTGraphon})
\begin{equation}\label{eq:constantEllieMatan}
c_W:=\sqrt{\frac{\int_x \deg_W(x)^2 \diff\mu(x)}{\|W\|_1^2}}\;.
\end{equation}
Let us state the main result from~\cite{GlobalUSTGraphon} using the epsilon-delta quantification (rather than using the language of convergence in which the result is stated in~\cite{GlobalUSTGraphon}). We state this theorem without explaining the Gromov--Hausdorff--Prokhorov (GHP) distance between metric measure spaces and refer the reader to~\cite{GlobalUSTGraphon}.
\begin{theo}\label{thm:AS}
    Given a connected graphon $U$, for every $\eps>0$ there exist $n(U,\eps)\in\N$ and $\delta(U,\eps)>0$ such that if $G$ is a connected graph on at least $n(U,\eps)$ vertices and at cut distance less than $\delta(U,\eps)$ from $U$, then the Lévy--Prokhorov distance with respect to the GHP distance between the random metric space of the uniform spanning tree of $G_n$ (with edge lengths $v(G_n)^{-1/2}$) and the rescaled continuum random tree $c_U\cdot \mathcal{T}$ is at most $\eps$.
\end{theo}

Crucially for us, if $U$ and $W$ are projectively fractionally isomorphic, then $c_U=c_W$; we prove this in Lemma~\ref{lem:EllieMatan}.

The aforementioned connection between the local and the scaling limit applies to a class of graphs which do not have sparse cuts. More precisely, for $\rho>0$ we say that a graph $G$ is \emph{$\rho$-robust} if for every $U\subset V(G)$ we have $e_G(U,V(G)\sm U)\ge \rho |U||V(G)\sm U|$.
\begin{cor}\label{cor:USTlocalglobal}
Let $\rho>0$ be an arbitrary constant. Suppose that $(G_n)_n$ is a sequence of $\rho$-robust graphs of growing orders $v(G_n)$. Suppose that the distribution of the uniform spanning trees of $G_n$ converges locally, as $n\to \infty$. Then the uniform spanning tree on $G_n$ whose each edge is given length $v(G_n)^{-1/2}$ converges in the sense of the scaling limit, as $n\to\infty$.
\end{cor}
\begin{proof}
Let $\mathcal{U}$ be the collection of graphons that arise as accumulation points of the sequence $(G_n)_n$ with respect to the cut distance topology. Theorems~\ref{thm:HNT} and~\ref{thm:mainUST} tell us that $\mathcal{U}$ is contained within one equivalence class of projective fractional isomorphism. Lemma~\ref{lem:EllieMatan} tells us that the constant defined in~\eqref{eq:constantEllieMatan} is universal for $\mathcal{U}$, say $c_{\mathcal{U}}$. The assumption of robustness implies that all the graphons in $\mathcal{U}$ are connected.

To prove the statement, it suffice to show that for every $\eps>0$, there exists $N_\eps\in\N$, such that for every $n>N_\eps$, the Lévy--Prokhorov distance between the random metric space of the uniform spanning tree of $G_n$ (with edge lengths $v(G_n)^{-1/2}$) and the rescaled continuum random tree $c_{\mathcal{U}}\cdot \mathcal{T}$ is at most $\eps$. To this end, fix $\eps>0$.

We claim that the constants $n(U,\eps)$ and $\delta(U,\eps)$ in Theorem~\ref{thm:AS} can be chosen uniformly over $U\in\mathcal{U}$.
\begin{claimUni}
There exist numbers $n(\mathcal{U},\eps)\in\N$ and $\delta(\mathcal{U},\eps)>0$ such that for all $U\in\mathcal{U}$, Theorem~\ref{thm:AS} holds with $n(U,\eps)=n(\mathcal{U},\eps)$ and $\delta(U,\eps)=\delta(\mathcal{U},\eps)$.
\end{claimUni}
\begin{proof}
    Suppose for a contradiction that the claim does not hold. Then when we consider the smallest working $n(U,\eps)$ in Theorem~\ref{thm:AS} for each $U\in\mathcal{U}$ and separately the largest working $\delta(U,\eps)$ for each $U\in\mathcal{U}$, we have $\sup_{U\in \mathcal{U}}n(U,\eps)=\infty$ or $\inf_{U\in \mathcal{U}}n(U,\eps)=0$. We shall focus on the former case and note that the latter case can be treated analogously.

    Since $\sup_{U\in \mathcal{U}}n(U,\eps)=\infty$, there exists a sequence $U_1,U_2,\ldots\in\mathcal{U}$ such that $\lim_{i\to \infty} n(U_i,\eps)=\infty$. By the Lov\'asz--Szegedy compactness theorem, this sequence has a subsequence convergent in the cut distance. We may, without loss of generality, assume that the sequence $U_1,U_2,\ldots$ itself is cut distance convergent, and call its limit graphon $U^*$. Since each $U_i$ is an accumulation point of the sequence $(G_n)_n$, we have a sequence of graphs $(G_{k_i})_i$ such that for each $i$ the cut distance between $U_i$ and $G_{k_i}$ is less than $1/i$. Now $U^*$ is an accumulation point of the sequence $(G_{k_i})_i$, so we have $U^*\in\mathcal{U}$. Theorem~\ref{thm:AS} applied to $U^*$ gives us numbers $n^*:=n(U^*,\eps)$ and $\delta^*:=\delta(U^*,\eps)$. In particular, Theorem~\ref{thm:AS} holds for our fixed $\eps$ and any graphon $U$ at cut distance less than $\delta^*/2$ from $U^*$ with numbers $n(U,\eps)=n^*$ and $\delta(U,\eps)=\delta^*/2$. In particular, for large enough $i$ we have $n(U_i,\eps)\le n^*$, which contradicts the divergence of the numbers $n(U_i,\eps)$.
\end{proof}

For each $n\in \N$, let $d_n$ be the infimum of the cut distance of $G_n$ to $U$, taken over $U\in\mathcal{U}$. We have $\lim_{n\to\infty} d_n=0$. Indeed, by contradiction, if this were not the case, then we could take a sequence with $\lim_{n_k\to\infty} d_{n_k}>0$. But any cut distance accumulation point of the graphs $(G_{n_k})_k$ (existence of which is guaranteed by the Lov\'asz--Szegedy compactness theorem) is contradictory.

In particular, for all sufficiently large $n$ we have $d_n<\delta(\mathcal{U},\eps)$ and $v(G_n)>n(\mathcal{U},\eps)$. Thus, there exists $U\in\mathcal{U}$ whose cut distance from $G_n$ is less than $\delta(\mathcal{U},\eps)$. Theorem~\ref{thm:AS} applied to the pair $U$ and $G_n$ gives the assertion about the Lévy--Prokhorov distance that we need.
\end{proof}

\subsection{A factorization result for fractional isomorphism}
In Section~\ref{sec:factorFI} we state and prove Theorem~\ref{thm:factorFI}, which is a factorization result for fractional isomorphism of disconnected graphons. Roughly, it says that two graphons $U$ and $W$ are fractionally isomorphic if and only if, for each connected graphon $\Gamma$, the total measure of connected components which are fractionally isomorphic to $\Gamma$ is the same in both $U$ and $W$. This result is an addition to the theory worked out in~\cite{Grebik2022} and does not directly concern graphon branching processes. The reason we include it here is that it allows us to extend Theorem~\ref{thm:mainUST} to disconnected graphons as we do in Theorem~\ref{thm:mainUSTGeneral}. We also include Theorem~\ref{thm:factorFIgraphs}, which is a graph counterpart of Theorem~\ref{thm:factorFI}. While Theorem~\ref{thm:factorFIgraphs} is not needed for our results about branching processes, we consider it an important contribution to the theory of graph fractional isomorphism.

\subsection{Organization of the paper}
Section~\ref{sec:notation} contains preliminaries for the proofs of our main results. In Sections~\ref{ssec:GraphonBasics}--\ref{ssec:StoneWeierstrass} we provide tools that will be useful for our proof of Theorem~\ref{thm:mainBJR}: general background on graphons, fractional isomorphism of graphons (mostly from~\cite{Grebik2022}, with minor adjustments), and a suitable version of the Stone--Weierstraß Theorem. New material appears in Section~\ref{sec:ProofThmMainBJR}. Specifically, Sections~\ref{ssec:BranchingDegree}--\ref{ssec:FindingSepTree} introduce further useful random rooted trees and establish identities which are useful for the proof of Theorem~\ref{thm:mainBJR}. The $(\Leftarrow)$ direction of Theorem~\ref{thm:mainBJR} is proved in Proposition~\ref{prop:CorrespondenceTwoBranchings}, while the $(\Rightarrow)$ direction of Theorem~\ref{thm:mainBJR} is proved at the end of Section~\ref{ssec:FindingSepTree}.


It turns out that the proofs of our results about random rooted trees $\fU_W$, Theorem~\ref{thm:mainUST} and its generalization for disconnected kernels given in Theorem~\ref{thm:mainUSTGeneral}, require more preparation. In Section~\ref{ssec:MarkovRenormalization} we introduce a transformation $W^\dagger$ of a graphon $W$ which is suitable to study the subtree $\fU_W^-$ of the random rooted tree $\fU_W$. In particular, the main result of Section~\ref{sssec:Extinction} says that $\fU_W^-$ goes extinct almost surely. We also need to recall the theory of discrete time Markov chains on uncountable measure spaces, which we do in Section~\ref{ssec:RWbackground}. In Section~\ref{sec:factorFI} we then state the factorization theorem for fractional isomorphism needed for extending Theorem~\ref{thm:mainUST} to Theorem~\ref{thm:mainUSTGeneral}. With these preparations, Section~\ref{sec:HNT} contains the proof of Theorem~\ref{thm:mainUSTGeneral}. 

\subsection{Acknowledgments}
We thank Ellie Archer and Jan Grebík for helpful discussions, and an anonymous referee for their valuable comments.

\section{Notation, preliminaries, and tools for proof of Theorem~\ref{thm:mainBJR}}\label{sec:notation}
We write $\N=\{1,2,\ldots\}$, $\N_0=\{0,1,2,\ldots\}$, $\N_\infty=\{1,2,\ldots,\infty\}$, and $\N_{0,\infty}=\{0,1,2,\ldots,\infty\}$. For $n\in\N_0$ we write $[n] = \{1,2,\ldots,n\}$ and $[n]_0 = \{0,1,2,\ldots,n\}$. We write $\Probability$ and $\Expectation$ for probability and expectation, respectively. 

All measures in this paper are assumed to be nonnegative. Given a measure space $(X,\B,\mu)$ and $p\in[1,\infty]$, we write $L^p(X,\B,\mu)$ for the space of all $\B$-measurable functions of finite $L^p$-norm with respect to $\mu$. When $\B$ or $\mu$ are clear from the context, we write just $L^p(X,\mu)$ or $L^p(X)$.

After introducing connected components of kernels in Section~\ref{ssec:GraphonBasics}, in Section~\ref{ssec:FracIso} we introduce the notion of fractional isomorphism which is used in the statement of our main results. Objects such as $\P^c$ and $\nu_W$ introduced in this section play a key role in the main proofs. The proof of Theorem~\ref{thm:mainBJR} relies on the version of the Stone--Weierstraß Theorem that we introduce in Section~\ref{ssec:StoneWeierstrass}. 

A summary of the notation introduced in the paper is also given in Appendix~\ref{app:notation}.

\subsection{Graph limits basics}\label{ssec:GraphonBasics}
We recall the notions of connected kernels (sometimes also called irreducible kernels) and connected components of kernels. These notions were thoroughly studied in~\cite{ConnectednessGraphons} and have since become standard. 
\begin{defi}\label{def:connected}
Let $I$ be a finite or countable set which does not contain~$0$. For a kernel $W$ on $(X,\B,\mu)$ we say that a decomposition $X=\Lambda_0\sqcup \bigsqcup_{i\in I}\Lambda_i$ into finitely or countably many $\B$-measurable sets is a \emph{decomposition into connected components of $W$ with isolated elements $\Lambda_0$} if the following hold.
\begin{enumerate}[label=(\roman*)]
    \item For $\mu$-almost every $x\in\Lambda_0$ we have $\deg_W(x)=0$.
    \item\label{en:connectedMain} For every $i\in I$ we have $\mu(\Lambda_i)>0$. Further, for every $B\subset \Lambda_i$ with $\mu(B)\in (0,\mu(\Lambda_i))$ we have $\int_{B\times(\Lambda_i\sm B)}W>0$.
    \item\label{en:connectedSep} For every $i\in I$ we have $\int_{\Lambda_i\times (X\setminus \Lambda_i)}W=0$.
\end{enumerate}
The sets $\Lambda_i$ ($i\in I$) are called the \emph{connected components} of $W$.
\end{defi}

It is shown in~\cite{ConnectednessGraphons} that a decomposition into connected components exists and is unique modulo nullsets. Further, if $W$ is nondegenerate, then $\Lambda_0$ is obviously a nullset. We say that $W$ is \emph{connected} if $\Lambda_0$ is a nullset and $|I|=1$, a definition consistent with the one we included in Footnote~\ref{foot:connected}.

\subsection{Fractional isomorphism}\label{ssec:FracIso}
In this section, we give a definition of fractional isomorphism for bounded-degree akernels. The theory developed by Grebík and Rocha~\cite{Grebik2022} gives five equivalent definitions of fractional isomorphism for graphons, and it straightforwardly generalizes in its entirety from graphons to kernels of bounded $L^\infty$-norm. We use one of those five definitions and generalize it to the slightly more general setting of bounded-degree akernels.\footnote{\label{foot:Example} 
An example of a bounded-degree kernel $W:(0,1)\rightarrow[0,\infty)$ which does not have finite $L^\infty$-norm is as follows. For $x,y$ such that $x,y\in(2^{-k},2^{-(k-1)})$ for some $k\in\N$, define $W(x,y)=2^k$. Otherwise, set $W(x,y)=0$.}
However, we do not claim that generalizations of their other definitions are equivalent or even well-defined. Later, in Section~\ref{ssec:MorePrelFI} we recall a different characterization of fractional isomorphism for the purposes of the proof of Theorem~\ref{thm:mainUST}.

To prepare for our definition of fractional isomorphism for bounded-degree akernels, we will do the following. First, we introduce a compact metric space $\P^c$ whose elements represent encapsulations of iterated degree information. Second, for each bounded-degree akernel $W$ we define a map $i_W \colon X \to \P^c$, where $i_W(x)$ captures the iterated degree information of $x$ in $W$. Finally, we assign a measure $\nu_W$ on $\P^c$ to each $W$; this measure plays a key role in our definition of fractional isomorphism.

\subsubsection{The space $\P^c$}
We recall a general definition of pushforward measure. Given a measurable map $f:A\to B$ between measurable spaces and a measure $\alpha$ on $A$, the \emph{pushforward} of $\alpha$ via $f$, denoted by $(f)_*\alpha$, is a measure on $B$ satisfying $(f)_*\alpha(S):=\alpha(f^{-1}(S))$ for every measurable $S\subset B$.

Let $c\in \R_+$. For a compact metric space $K$ we write $\M_{\leq c}(K)$ for the set of all Borel measures on $K$ with total mass at most $c$ and we define $\M_{= c}(K)$ analogously. We equip both $\M_{\leq c}(K)$ and $\M_{=c}(K)$ with the topology of weak convergence. We recall the basics and refer to Sections~1,~2 (basics) and Section~5 (Prokhorov's Theorem) of~\cite{MR1700749} for details. The weak topology of measures is metrizable and thus also characterized by convergent sequences; a sequence $\zeta_1,\zeta_2,\ldots\in \M_{\leq c}(K)$ \emph{converges weakly} to a measure $\zeta$ if for every bounded continuous function $f:K\rightarrow\R$ we have
\begin{equation}\label{eq:defWeakConv}
    \lim_n \int f\diff\zeta_n=\int f\diff\zeta\;.
\end{equation}
We recall that $\M_{\leq c}(K)$ is compact; this result is sometimes referred to in literature as Prokhorov's Theorem, the Prokhorov--Banach--Alaoglu Theorem, or simply as weak compactness of the space of bounded measures.

To construct $\P^c = \P^c_{\infty}$ we shall recursively define collections $\{L_k^c\}_{k\in\N_0}$ and $\{\P_k^c\}_{k\in\N_0}$ of spaces as well as canonical projections $p_{\ell,k}^c$ from $ \P_k^c$ to $ \P_\ell^c$ for $k>\ell$. These definitions are intertwined: first a trivial definition $L_0^{c}$ is given, then one can construct, in this order, $\P_0^{c}$ using~\eqref{eq:Pfinite}, $p_{0,0}^c$ using~\eqref{notation:projection}, $L_1^c$ using~\eqref{notation:Llittle}, $\P_1^{c}$ using~\eqref{eq:Pfinite}, $p_{0,1}^c,p_{1,1}^c$ using~\eqref{notation:projection}, $L_2^c$ using~\eqref{notation:Llittle}, and so on. Immediately after defining each space, we provide an argument establishing its compactness. Each such argument relies on the compactness of the spaces defined in earlier steps. Let $L_0^{c}=\{\OneElementSpace\}$ be the one-point space. This space is obviously compact.
For $k\in\N_0$, let 
\begin{equation}
\label{eq:Pfinite}
\P_k^{c}:=\Bigg\{\alpha\in\prod_{i=0}^{k}L_i^{c}\ \Bigg\vert\  \underbrace{\forall \ell \in [k-1]:\ \alpha(\ell)=(p^c_{\ell-1,\ell})_* \alpha(\ell+1)}_{\mathsf{(C)}_k}\Bigg\},
\end{equation}
where the condition $\forall \ell \in [k-1]$ is vacuous for $k=0,1$.
The space $\P_k^{c}$ is compact since $\mathsf{(C)}_k$ specifies a closed subset of a product of compact spaces.
Let 
\begin{equation}
\label{notation:projection}
    p_{h,k}^c\colon \P_k^c\to \P_h^c
\end{equation}
be the canonical projection for $h\in[k]_0$, and let 
\begin{equation}
\label{notation:Llittle}
L_{k+1}^{c}:=\M_{\leq c}(\P_k^{c})\;.    
\end{equation}
The space $L_{k+1}^{c}$ is compact by Prokhorov's Theorem.
We note that the canonical projections $p_{h,k}^c$ are well defined, that is, if $\alpha\in \prod_{i=0}^{k}L_i^{c}$ satisfies~$\mathsf{(C)}_k$, then $p_{h,k}^c(\alpha)\in \prod_{i=0}^{h}L_i^{c}$ satisfies~$\mathsf{(C)}_h$. 

The canonical projections, for $n,k\in\N_0$ with $n\leq k$, naturally give rise to the pushforwards $(p_{\ell,k}^c)_*\colon L_{k+1}^c\to L_{\ell+1}^c$.
Finally, let 
\begin{equation}\label{eq:Pinfty}
\P^{c}=\P^c_\infty=\left\{\left.\alpha\in\prod_{i=0}^{\infty}L_i^{c}\right\vert\ \forall k\in\N:\ \alpha(k)=(p^c_{k-1,k})_* \alpha(k+1)\right\}
\end{equation}
and for $k\in \N_0$ let 
\begin{equation}
    p_{k,\infty}^c\colon \P_\infty^c\to \P_k^c\;
\end{equation}
be the corresponding canonical projections.
The space $\P^{c}$ is compact for the same reason that spaces $\P_k^{c}$ are.  We write $\B(\P^c)$\label{notation:BorelP} and $\B(\P_n^c)$ for the Borel sets on $\P^c$ and $\P_n^c$, respectively.

For each $k\in \N_\infty$ and $\alpha\in \P_k^c$, we have $\alpha(1)(\P^{c}_0)=\alpha(2)(\P^{c}_1)=\ldots=\alpha(k)(\P^{c}_{k-1})\leq c$. Throughout the paper, we denote this value by $D(\alpha)$.\label{notation:degreeMeasure} 

\subsubsection{An example}
Elements of $\P^c$ will be used to record iterated degrees in bounded-degree akernels. For example, consider $\alpha\in\P^{0.8}$ such that
\begin{align}
\begin{split}\label{eq:Example}
\alpha(0)&=\OneElementSpace\;,\\
\alpha(1)&=0.8\cdot \Dirac(\OneElementSpace)\;,\\
\alpha(2)&=0.5\cdot \Dirac\big((\OneElementSpace,0.4\cdot \Dirac(\OneElementSpace))\big)+0.3\cdot \Dirac\big((\OneElementSpace,0.7\cdot \Dirac(\OneElementSpace))\big)\;, \\
\ldots\;.
\end{split}
\end{align}
If such an $\alpha$ encodes information about the iterated degrees in an akernel $W$ at a vertex $x$, then this means that
\begin{align}
\label{eq:degEx1}
&\deg_W(x)=0.8,\\
\label{eq:degEx2}
&\int_{y:\deg_W(y)=0.4} W(x,y)=0.5\;,\;
\int_{y:\deg_W(y)=0.7} W(x,y)=0.3\;\mbox{,}\;
\int_{y:\deg_W(y)\not\in\{0.4,0.7\}} W(x,y)=0.
\end{align}
In particular, the condition $\alpha(k)=(p^c_{k-1,k})_* \alpha(k+1)$ in~\eqref{eq:Pinfty} expresses consistency of the iterated degrees. To give a specific example of violating this consistency condition, consider~\eqref{eq:degEx2} but modify~\eqref{eq:degEx1} to $\deg_W(x)=0.9$.

\subsubsection{Iterated degree measures and the definition of fractional isomorphism}
We now recall a certain operation of lifting elements of $\P^c$ to measures in $\M_{\le c}(\P^c)$ from~\cite{Grebik2022}. Take the exemplar $\alpha\in\P^{1}$ from~\eqref{eq:Example} to gain some intuition. This $\alpha$ is in a way expressing an object whose total mass is~0.8, and this total mass is split in the ratio $0.5:0.3$ between objects whose total mass is $0.4$ and $0.7$, respectively. We could equivalently encode such information into an object $\mu_\alpha\in \M_{\le 1}(\P^{1})$ such that $\mu_\alpha(\P^{1})=0.8$, and $\mu_\alpha(Y_{0.4})=0.5$ and $\mu_\alpha(Y_{0.7})=0.3$, where for $r\ge 0$, $Y_r:=\{\beta\in\P^1:D(\beta)=r\}$.\label{notation:lifting} The general definition is as follows. For $\alpha\in\P^c$ let $\mu_\alpha\in \M_{\le c}(\P^c)$ be the unique measure that satisfies 
\begin{equation}\label{eq:lifting}
(p_{k,\infty})_*\mu_\alpha=\alpha(k+1)    
\end{equation}
for every $k\in\N$. The soundness of this definition follows from Kolmogorov's extension theorem together with the consistency condition $\alpha(k)=(p^c_{k-1,k})_* \alpha(k+1)$ in~\eqref{eq:Pinfty}.

Let $W$ be a bounded-degree akernel with $\degmax(W)\le c$. By modifying $W$ on a $(\mu\times\mu)$-null set, we may assume without loss of generality that $W$ satisfies $\deg_W(x) \le c$ for all $x\in X$. We remark that the collection of $x\in X$ for which the object $i_W(x)$ (introduced below) changes is a $\mu$-null set, so in particular the key object $\nu_W$ defined below will remain unchanged. For $k=0,1,2\ldots$ we recursively define $i_{W,k}:X\to \P^c_k$ in the following way.
Let $i_{W,0}\colon X\to \P^c_0=\{\OneElementSpace\}$ be the constant map. For $k\in \N_0$, let $i_{W,k+1}$ be constructed from $i_{W,k}$ so that
\begin{align}
    i_{W,k+1}(x)(j)&=i_{W,k}(x)(j)&\text{ for }0\leq j\leq k \text{ and}\\
    \label{eq:idm_lift}i_{W,k+1}(x)(k+1)(A)&=\int_{i^{-1}_{W,k}(A)}W(x,y)\diff\mu(y) &\text{ for every measurable $A\subset \P^c_{k}$}.
\end{align}
That is, the function $i_{W,1}$ carries information about degrees of vertices in the akernel $W$, the function $i_{W,2}$ carries information about a refined version of degrees such as `the degree of a vertex into the set of vertices whose degree lies in the interval $[0.4,0.5]$', and so on.
Let \label{notation:i}$i_W\colon X\to \P^c$ be defined by $i_{W}(x)(k)=i_{W,n}(x)(k)$. It follows from the definitions given above that the map $i_W$ is well defined, i.\,e., we have $i_W(x)\in\P^c$ for every $x\in X$; the key details can be found in the proof of~\cite[Proposition~6.8]{Grebik2022}. Finally, let $\nu_W=(i_W)_*\mu$ be the pushforward of $\mu$ via $i_W$. Obviously, we have $\nu_W\in\M_{= 1}(\P^c)$. We call~$\nu_W$\label{notation:iterateddegreemeasure} the \emph{iterated degree measure} of the bounded-degree akernel $W$.

\begin{rem}[The bound $c$ in the definition of $\P^c$]\label{rem:PcBounded}
Our main results assume some type of boundedness of the kernel or akernel in question. The mildest such assumption is $\degmax(W)\le c$. This assumption guarantees that the support of~$\nu_W$ is contained in $\P^c$. To work with kernels without bounded degrees, we could easily repeat the construction and build a space
$\P^{<\infty}$ using sets $L_{k+1}^{<\infty}:=\M_{<\infty}(\P_k^{<\infty})$. However, by extending from $\P^c$ to $\P^{<\infty}$ we lose compactness, a prominent feature needed for an application of the Stone--Weierstraß Theorem, which in turn plays an important role in our proof of Theorem~\ref{thm:mainBJR}. A discussion of possible generalizations can be found in Section~\ref{ssec:boudnednessassumption}. 

This is the only role the bound $c$ plays, and its value is fairly immaterial. This is reflected by the fact that we usually assume the value of $c$ implicitly fixed and remove it from the main quantification of our lemmas. 
\end{rem}

Out of several equivalent definitions of fractional isomorphism from~\cite{Grebik2022}, we recall the one based on iterated degree measures.  
\begin{defi}\label{def:FracIso}
    We say that two bounded-degree akernels $U$ and $W$ are \emph{fractionally isomorphic} if $\nu_U=\nu_W$.
\end{defi}

The last result in this section rectifies an omission from Section~\ref{ssec:IntroUST}. 
\begin{lem}\label{lem:EllieMatan}
Suppose that $U$ and $W$ are two projectively fractionally isomorphic $L^\infty$-graphons. Then the constants $c_U$ and $c_W$ defined in~\eqref{eq:constantEllieMatan} are equal.
\end{lem}
\begin{proof}
Suppose that $U$ is fractionally isomorphic to $tW$ for some $t>0$, i.\,e., $\nu_U=\nu_{tW}$. Looking at~\eqref{eq:constantEllieMatan}, we can replace $\deg_W(x)$ by $D(i_W(x))$. Also, $\|W\|_1=\int_x \deg_W(x)\diff \mu(x)=\int D(\alpha)\diff \nu_W(\alpha)$.
Hence, we obtain
\begin{align*}
    c_U^2&=\frac{\int_x \deg_U(x)^2 \diff\mu(x)}{\|U\|_1^2}=\frac{\int D(\alpha)^2\diff \nu_U(\alpha)}{\left(\int D(\alpha)\diff \nu_U(\alpha)\right)^2}
    =\frac{\int D(\alpha)^2\diff \nu_{tW}(\alpha)}{\left(\int D(\alpha)\diff \nu_{tW}(\alpha)\right)^2}=\frac{\int_x \deg_{tW}(x)^2 \diff\mu(x)}{\|tW\|_1^2}\\
    &=\frac{\int_x t^2\deg_{W}(x)^2 \diff\mu(x)}{t^2\|W\|_1^2}=c_W^2
\end{align*}
as required.
\end{proof}

\subsection{The Stone--Weierstraß Theorem}\label{ssec:StoneWeierstrass}
We recall the setting of the Stone--Weierstraß Theorem.
Given a topological space $Z$, we write $C(Z,\R)$ for the set of all continuous functions from $Z$ to $\R$. Recall that a family $\mathcal{E}\subset C(Z,\R)$ is \emph{multiplicative} if for every $f,g\in \mathcal{E}$ the function $h$ defined as $h(z)=f(z)g(z)$ satisfies $h\in\mathcal{E}$.
A family $\mathcal{E}\subset C(Z,\R)$ is \emph{multiplicative up to constants}, if for every $f,g\in \mathcal{E}$ there exists an $a\in \R\sm\{0\}$ such that the function $h$ defined as $h(z)=af(z)g(z)$ lies in $\mathcal{E}$.
Given distinct elements $x,y\in Z$, the family $\mathcal{E}$ \emph{separates} $x$ and $y$ if there exists $f\in\mathcal{E}$ such that $f(x)\neq f(y)$; the family $\mathcal{E}$ \emph{separates points} if it separates every pair of distinct elements $x,y\in Z$.

Like in \cite{Grebik2022} we use a corollary of the Real Stone--Weierstraß Theorem, which can be found in \cite[Theorem~7.32]{rudin1976principles}. 

\begin{cor}\label{cor:StoneWeierstrass}
    Let $K$ be a compact metric space and $c>0$ be a real number. Suppose that $\mathcal{E}\subset C(K,\mathbb{R})$ is a family of functions that is multiplicative up to constants, contains the constant-1 function, and separates points. Then for every distinct $\alpha,\beta\in \M_{\leq c}(K)$ there is $f\in \mathcal{E}$ such that
$$\int_K f \diff\alpha\neq \int_K f \diff\beta\;.$$
\end{cor}

\section{Proof of Theorem~\ref{thm:mainBJR}}
\label{sec:ProofThmMainBJR}

For a kernel $W$, the offspring in the random rooted tree $\fX_W$ are generated via a Poisson point process on $W$. Unfolding the definition of fractional isomorphism, Theorem~\ref{thm:mainBJR} says that it is not the entire complexity of $W$ that determines $\fX_W$ but rather less complex information encoded in the iterated degrees.

In Section~\ref{ssec:BranchingDegree} we show that the fractional isomorphism type of a bounded-degree akernel $W$ determines the distribution of $\fX_W$ by constructing a process $\fB_W$ that depends only on the iterated degrees of vertices and not on the vertices themselves. In particular, this establishes the $(\Leftarrow)$ direction of Theorem~\ref{thm:mainBJR}.

In Section~\ref{ssec:ProbGenTree} we express the probabilities that $\fB_W$ (or rather, a related process $\fB(\alpha)$) generates a given tree. This is in particular used to establish, in Lemma~\ref{lem:FMult}, that certain families $\overline{\F_n}$ of functions on $\P_n^c$ which are related to these tree-generation probabilities are multiplicative up to constants. We also show in Lemma~\ref{lem:TreeFunctionsContinuous} that the functions in $\overline{\F_n}$ are continuous. In Section~\ref{ssec:FindingSepTree} we state Proposition~\ref{prop:keyseparation}, which asserts that $\overline{\F_n}$ also separates points. This puts us in a position where we can use the Stone--Weierstraß Theorem to show that for  each pair of non-fractionally-isomorphic akernels $U$ and $W$ there is a rooted tree $T$ and $k\in\N_0$ such that $\Probability\left[(\fX_U)_{\restriction k}\cong T\right]\neq \Probability\left[(\fX_W)_{\restriction k}\cong T\right]$, hence establishing the $(\Rightarrow)$ direction of Theorem~\ref{thm:mainBJR}. We include the deferred proof of Proposition~\ref{prop:keyseparation} in Section~\ref{ssec:proofSeparation}.

\subsection{The $(\Leftarrow)$ direction of Theorem~\ref{thm:mainBJR}}\label{ssec:BranchingDegree}

For a bounded-degree akernel $W$, we introduce a random rooted tree $\fB_W$ with type space $\P^c$. We show in Proposition~\ref{prop:CorrespondenceTwoBranchings} that this random rooted tree is closely related to the random rooted tree $\fX_W$.

Recall that for each $\beta\in\P^c$ we have a unique measure $\mu_\beta\in \M_{\le c}(\P^c)$ satisfying~\eqref{eq:lifting}. Given any $\alpha\in\P^c$, the random rooted tree \label{notation:fBallpha}$\fB(\alpha)$ is given as follows. The type space of $\fB(\alpha)$ is the set $\P^c$. The root has type $\alpha$, and a particle of type $\beta\in\P^c$ has children distributed as a Poisson point process on $\P^c$ with intensity $\mu_\beta$. Given a bounded-degree akernel $W$, the random rooted tree \label{notation:fBW}$\fB_W$ with type space $\P^c$ is given as follows. The root has type according to the distribution $\nu_W$ and a particle of type $\beta\in\P^c$ has children distributed as a Poisson point process on $\P^c$ with intensity $\mu_\beta$.

The following fact encapsulates the close connection between the random rooted trees $\fB_W$ and $\fX_W$ via the function $i_W$ defined in Section~\ref{ssec:FracIso}. In particular, Proposition~\ref{prop:CorrespondenceTwoBranchings}\ref{corespond3} establishes the $(\Leftarrow)$ direction of Theorem~\ref{thm:mainBJR}.
 For $x\in X$, we write $\fX_W(x)$ for the random rooted tree $\fX_W$ conditioned to start at $x$.

\begin{prop}\label{prop:CorrespondenceTwoBranchings}
\begin{enumerate}[label=(\roman*)]
\item\label{corespond1} Suppose that $W$ is a bounded-degree akernel. Then for each $x\in X$ the random rooted trees $\fX_W(x)$ and $\fB(i_W(x))$ have the same distribution.
\item\label{corespond2} Suppose that $W$ is a bounded-degree akernel.
Then the random rooted trees $\fX_W$ and $\fB_W$ have the same distribution.
\item\label{corespond3} Suppose that $U$ and $W$ are two fractionally isomorphic bounded-degree akernels. Then the random rooted trees $\fX_U$ and $\fX_W$ have the same distribution.
\end{enumerate}
\end{prop}
\begin{proof}
The implication \ref{corespond1}$\Rightarrow$\ref{corespond2} is clear because $\nu_W$ is the pushforward of $\mu$ via $i_W$. The implication \ref{corespond2}$\Rightarrow$\ref{corespond3} is also clear because two fractionally isomorphic bounded-degree akernels $U$ and $W$ have the same iterated degree measure $\nu_U = \nu_W$ by definition. It remains to prove~\ref{corespond1}.

While we generally ignore type labels and compare random rooted trees via their distributions on type-unlabeled rooted trees, here we will retain the type labels. More specifically, we prove the desired statement by showing that we may couple the distributions of $\fX_W(x)$ and $\fB(i_W(x))$ on type-labeled rooted trees by using the natural transformation $x\mapsto i_W(x)$ on particle types.

We show by induction on $k$ that for every $k\in\N_0$ and $x\in X$, the $k$-balls $(\fX_W(x))_{\restriction k}$ and $(\fB(i_W(x)))_{\restriction k}$ have the same distribution on type-labeled rooted trees. This is trivial for $k=0$. For step $k\in\N$ observe that $(\fX_W(x))_{\restriction k}$ consists of the root of type $x$, children of types $x_1,\ldots,x_h\in X$ following a Poisson point process on $X$ with intensity $\pi$ given by $\diff\pi = W(x,\cdot) \diff\mu$, and trees $(\fX_W(x_1))_{\restriction k-1}$, \ldots, $(\fX_W(x_h))_{\restriction k-1}$ attached to them. On the other hand, $(\fB(i_W(x)))_{\restriction k}$ consists of the root of type $i_W(x)$, children of types $\alpha_1,\ldots,\alpha_\ell\in \P^c$ following a Poisson point process on $\P^c$ with intensity $\mu_\alpha$, and trees $(\fB(\alpha_1))_{\restriction k-1}$, \ldots, $(\fB(\alpha_\ell))_{\restriction k-1}$ attached to them. Now by the Mapping Theorem for point processes, $i_W(x_1),\ldots,i_W(x_h)$ follows a Poisson point process on $\P^c$ with intensity $(i_W)_*\pi$. By the definitions of pushforward measure, $i_W$ and $\mu_{i_W(x)}$ (see~\eqref{eq:idm_lift} and~\eqref{eq:lifting}), we have $(i_W)_*\pi = \mu_{i_W(x)}$. Hence, the inductive step is completed by the inductive hypothesis.
\end{proof}

\subsection{The probability of generating a tree}\label{ssec:ProbGenTree}
In this section we introduce notation that will allow us to express the probability that the first $k$ levels of $\fB(\alpha)$ yield a given tree. Constructions to this end are very similar to \cite[Section 7.2]{Grebik2022}, in which measures $\alpha\in\P^1$ were used to express the rooted homomorphism density of a rooted tree. While there is some resemblance between our branching processes and homomorphism densities of trees, the two concepts are not directly related.

For a finite rooted tree $T$ with root $r$, its \emph{height} is given by $h(T):=\max(\dist(r,v)\mid v\in V(T))$.
For $n\in \N$ let \label{notation:Tn}$\T_n$ be the set of isomorphism classes of finite rooted trees of height at most $n$. We denote the $1$-vertex rooted tree by $\square$, i.\,e., $\T_0=\{\square\}$. 
    
We introduce two basic operations for building a rooted tree from smaller ones. The first one plants a new root, while the second one merges several trees into one by identifying their roots.
\begin{itemize}
    \item Suppose that $T$ is a rooted tree with root $r$. Then let $T^\uparrow$\label{notation:uparrow} denote the tree consisting of the unrooted version of $T$ with a new root attached to the previous root. In particular, we have $v(T^\uparrow)=v(T)+1$ and $h(T^\uparrow)=h(T)+1$.
    \item Suppose that $T_1,\ldots,T_\ell$ are finite rooted trees. The rooted tree \label{notation:oplus}$T_1\oplus\ldots \oplus T_\ell$ is constructed by taking their disjoint union while identifying all the roots. In particular, we have $v(T_1\oplus\ldots\oplus T_\ell)=v(T_1)+\cdots+v(T_\ell)-\ell+1$ and $h(T_1\oplus\ldots\oplus T_\ell)=\max\{h(T_1),\ldots,h(T_\ell)\}$. 
\end{itemize} 
For each finite rooted tree $T$ and each $k\in\N_0$, we define two functions $f_{T,k},\overline{f_{T,k}}\colon \P_k^c\to \R$. 

Lemma~\ref{lem:generatingBranchings} connects our algebraic Definition~\ref{def_treefunc} of $f_{T,n}(p_{n,\infty}(\alpha))$ to the branching process $\fB(\alpha)$. The way we define $\overline{f_{T,n}}(\beta)$ (in Definition~\ref{def_treefunc}\ref{f_bar}) based on $f_{T,n}(\beta)$ is quite bland, yet convenient for the application of the Stone--Weierstraß Theorem.

Suppose that $F$ is a finite rooted tree with root $r$. Let $\mathcal{L}_F=\{\ell_1,\ldots\ell_t\}$ be the multiset of multiplicities of isomorphism types of connected components of $F-r$. For example, $\mathcal{L}_F=\{1,3,3\}$ means that $F-r$ has 7 connected components in total, comprising one triple of mutually isomorphic connected components, another triple of mutually isomorphic connected components, and a further unique connected component. Note that $\sum\mathcal{L}_F=\deg_F(r)$. Define a constant $e_F$ using elements in the multiset $\mathcal{L}_F$,
\begin{equation}\label{eq:constE}
e_{F}:=\prod_{\ell\in\mathcal{L}_F}\frac{1}{\ell!}\;.
\end{equation}

\begin{defi}\label{def_treefunc}
    \begin{enumerate}[label=(\roman*)]
    \item\label{en:froot} For each $k\in\N_0$ and each $\beta \in \P_k^c$, we have  $$f_{\square,k}(\beta)=\begin{cases}
        1& \text{ if $k=0$,}\\
        \exp(-D(\beta))& \text{ if $k>0$.}
    \end{cases}$$
    
    \item\label{new_root} For a finite rooted tree $T$ of height at most $k\in\N_0$ 
    and $\beta\in\P_{k+1}^c$ define
        $$f_{T^\uparrow,k+1}(\beta)=\exp(-D(\beta))\cdot\int_{\sigma\in \P_k^{c}} f_{T,k}(\sigma) \diff\big(\beta(k+1)\big)(\sigma)\;.$$
    \item\label{multiplication} For finite rooted trees $T_1,\ldots, T_\ell$ of positive heights not exceeding $k\in\N$ and whose roots have degree exactly~1, and for $\beta\in\P_k^c$, $$f_{T_1\oplus\ldots \oplus T_\ell,k}(\beta)=e_{T_1\oplus\ldots\oplus T_\ell}\cdot\exp(D(\beta))^{\ell-1}\prod_{i=1}^\ell f_{T_i,k}(\beta)\;.$$
Note that this definition is consistent for $\ell=1$, that is, it says $f_{T_1,k}(\beta)=1\cdot\exp(D(\beta))^{0} f_{T_1,k}(\beta)$.
    \item\label{f_bar} For each finite rooted tree $T$ and each $k\in\N_0$, we define $\overline{f_{T,k}}(\beta)=\exp(D(\beta))f_{T,k}(\beta)$. In particular, we have $\overline{f_{\square,1}}(\beta)=1$.
    \end{enumerate}
\end{defi}

Every rooted tree can be decomposed into a collection of subtrees rooted at the children of the original root, and we can perform this decomposition recursively. Hence, the function $f_{T,k}$ is given for every $k\in \N_0$ and for every finite tree $T$ with $h(T)\leq k$. Let \label{notation:cF}$\F_n:=\{f_{T,n}\mid T\in \T_n\}$ and \label{notation:overlineFdash}$\overline{\F_n}:=\{\overline{f_{T,n}}\mid T\in \T_n\}$.

The next lemma is needed to apply Corollary~\ref{cor:StoneWeierstrass}.
\begin{lem}\label{lem:FMult}
    For each $n\in\N$, the set $\overline{\F_n}$ is multiplicative up to constants. 
\end{lem}

\begin{proof}
    Let $T$ and $F$ be finite rooted trees of height at most $n$. We shall show that $\overline{f_{T, n}}\cdot \overline{f_{F,n}}=\frac{e_Te_F}{e_{T\oplus F}}\cdot \overline{f_{T\oplus F, n}}$. As $h(T\oplus F)=\max(h(T),h(F))$, the function $\overline{f_{T\oplus F, n}}$ is contained in $\overline{\F_n}$. This will show multiplicativity up to constants.
    
    Let $T_1,\ldots,T_\ell$ be the connected components of $T$ after the removal of the root. Let $F_1,\ldots,F_t$ be the connected components of $F$ after the removal of the root. Note that the connected components of $T\oplus F$ after the removal of the root are $T_1,\ldots,T_\ell,F_1,\ldots,F_t$. For each $\alpha\in \P_n^c,$ we have
    \begin{align*}
    \overline{f_{T, n}}(\alpha)\cdot \overline{f_{F,n}}(\alpha)&=\exp(D(\alpha))^2f_{T,n}(\alpha)\cdot f_{F,n}(\alpha)\\
    &=\exp(D(\alpha))^2\cdot e_{T}\cdot\exp(D(\alpha))^{\ell-1}\prod_{i=1}^\ell f_{T_i^\uparrow,n}(\alpha)\cdot
    e_{F}\cdot\exp(D(\alpha))^{t-1}\prod_{i=1}^t f_{F_i^\uparrow,n}(\alpha)\\
    &=\frac{e_Te_F}{e_{T\oplus F}}\cdot \exp(D(\alpha))\cdot f_{T\oplus F, n}(\alpha)=\frac{e_Te_F}{e_{T\oplus F}}\cdot \overline{f_{T\oplus F, n}}(\alpha)\;. \qedhere
    \end{align*}
\end{proof}

The next lemma connects the above analytic definition of functions $f_{T,n}(\cdot)$ with branching processes $\fB(\cdot)$.
\begin{lem}\label{lem:generatingBranchings}
	Given $k\in \N_0$, a rooted tree $T$ of height at most $k$, and $\alpha\in\P^c$, the quantity $f_{T,k}(p^c_{k,\infty}(\alpha))$ is the probability that the first $k$ levels of $\fB(\alpha)$ are isomorphic to $T$.
\end{lem}
\begin{proof}
The claim is trivial for $k=0$. Indeed, the zeroth level of $\fB(\alpha)$ is always isomorphic to $T=\square$, the only tree of height at most $0$, and we have $f_{\square,0}(p^c_{0,\infty}(\alpha))=1$.

For $k>0$ we proceed by induction on $k$. Let $r$ be the root of $T$. Write $s=\deg_T(r)$. Let $\mathcal{L}_T=\{\ell_1,\ldots\ell_t\}$ be the multiset of multiplicities of isomorphism types of connected components of $T-r$. Let $T_1,\ldots,T_t$ be the isomorphism types of these connected components of $T-r$ corresponding to these multiplicities.

First, let us argue when $D(\alpha)=0$. Indeed then on the one hand, $\mu_\alpha$ is the trivial measure and $\fB(\alpha)$ is almost surely just the root. Let us now turn to evaluating $f_{T,k}(p^c_{k,\infty}(\alpha))$. If $T=\square$ then Definition~\ref{def_treefunc}\ref{en:froot} gives that $f_{T,k}(p^c_{k,\infty}(\alpha))=\exp(-D(\alpha))=1$. If $T\neq\square$, then $r$ has at least one nontrivial subtree attached to it. Using the product formula of Definition~\ref{def_treefunc}\ref{multiplication}, in order to show that $f_{T,k}(p^c_{k,\infty}(\alpha))=0$, it suffices to see nullification on this one subtree. This can be seen from the integral formula in Definition~\ref{def_treefunc}\ref{new_root} as the measure $\alpha(k)$ is trivial.

Thus, in the following we can assume that $D(\alpha)>0$. Generate the first generation of $\fB(\alpha)$. If the first $k$ levels of $\fB(\alpha)$ are isomorphic to $T$ then there are $s$ many offspring generated in the first generation. We shall assume that this occurred. Let $[s]$ be an enumeration of the offspring of the root considered in a uniformly random order. For a given partition $\mathcal{U}\in\binom{[s]}{\mathcal{L}_F}$, let $E_\mathcal{U}$ be the event for each $i\in[t]$ the $\ell_i$ many elements of $[s]$ corresponding to the $i$th type of $\mathcal{U}$ generate a tree isomorphic to $T_i$. The events $E_\mathcal{U}$ are disjoint, their union is the event that the first $k$ levels of $\fB(\alpha)$ are isomorphic to $T$, and (because the ordering is random) they all occur with the same probability. In particular, we can focus on the partition $\mathcal{U}^*=(\{1,\ldots,\ell_1\},\{\ell_1+1,\ldots,\ell_1+\ell_2\},\ldots,\{\ell_1+\ell_2+\ldots+\ell_{t-1}+1,\ldots,s\})$. This gives that
\begin{align}\label{eq:vlak1}
\Probability\left[\fB(\alpha)_{\restriction k}\cong T\right]
=\binom{s}{\mathcal{L}_F}\cdot\frac{D(\alpha)^{s}\exp(-D(\alpha))}{s!}\cdot\prod_{i=1}^t p_i^{\ell_i}\;,
\end{align}
where $p_i$ is the probability that $T_i$ is isomorphic to the first $k-1$ level of the branching process $\fB(\beta)$, where $\beta\in\P_{k-1}^c$ is chosen at random according to the probability distribution $\nicefrac{\alpha(k)}{D(\alpha)}$. That is, using the induction hypothesis, we have 
\begin{equation}\label{eq:vlak2}
p_i= \int_{\beta\in \P_{k-1}^{c}} f_{T_i,k-1}(\beta)\;\diff\left(\nicefrac{\alpha(k)}{D(\alpha)}\right)(\beta)\;.
\end{equation}
Putting~\eqref{eq:vlak1} and~\eqref{eq:vlak2} together we have
\begin{equation}
\label{eq:Konec1}
\Probability\left[\fB(\alpha)_{\restriction k}\cong T\right]=\frac{D(\alpha)^{s}\exp(-D(\alpha))}{\prod_{i=1}^t\ell_i!}\cdot\prod_{i=1}^t p_i^{\ell_i}=
\frac{\exp(-D(\alpha)) }{\prod_{i=1}^t\ell_i!}\cdot\prod_{i=1}^t \left(\int_{\sigma\in\P_{k-1}^c} f_{T_i,k-1}(\sigma)\;\diff(\alpha(k))(\sigma)\right)^{\ell_i}\;.
\end{equation}
We now shift to expanding $f_{T,k}(\alpha)$. We first use Definition~\ref{def_treefunc}\ref{multiplication} and then Definition~\ref{def_treefunc}\ref{new_root}.
\begin{align}
f_{T,k}(p^c_{k,\infty}(\alpha))&=\frac{1}{\prod_{i=1}^t\ell_i!}\cdot\exp(D(\alpha))^{s-1}\prod_{i=1}^t \left(f_{T_i^\uparrow,k}(p^c_{k,\infty}(\alpha))\right)^{\ell_i}\\
\nonumber
&=\frac{1}{\prod_{i=1}^t\ell_i!}\cdot\exp(D(\alpha))^{s-1}\cdot 
\exp(-D(\alpha))^s\cdot\prod_{i=1}^t\left(\int_{\sigma\in \P_k^{c}} f_{T_i,k-1}(\sigma) \diff\big(\alpha(k)\big)(\sigma)\right)^{\ell_i}\\
\label{eq:Konec2}
&=\frac{\exp(-D(\alpha))}{\prod_{i=1}^t\ell_i!}\cdot\prod_{i=1}^t\left(\int_{\sigma\in \P_k^{c}} f_{T_i,k-1}(\sigma) \diff\big(\alpha(k)\big)(\sigma)\right)^{\ell_i}\;.
\end{align}
This finishes the proof.
\end{proof}

The Stone--Weierstraß Theorem works in the setting of continuous functions. The next lemma shows that this is indeed the case for the functions $f_{T,n}$.
\begin{lem}\label{lem:TreeFunctionsContinuous}
    Given $k\in\N_0$ and a rooted tree $T$ of height at most $k$, the functions $f_{T,k}:\P_k^c\to \R$ and $\overline{f_{T,k}}:\P_k^c\to \R$ are continuous.
\end{lem}
\begin{proof}
Let us first think about the topology on $\P^c_k$. As a matter of fact, this is the only place in the paper where we need to study the topology on $\P^c_k$. The only topological fact about $\P^c_k$ used at other places is that of compactness. Recall that the topologies on $L_i^c$ are metrizable, and hence so is the product topology on $\P^c_k$. Thus the topology is characterized by convergent sequences.
\begin{claim:lem:TreeFunctionsContinuous1}
Let $k\in\N$ be arbitrary. Suppose that $\alpha_1,\alpha_2,\ldots\in \P^c_k$ converge in $\P^c_k$ to $\alpha\in \P^c_k$. Then for every bounded continuous function $f: \P^c_{k-1}\rightarrow\R$ we have 
\begin{equation}\label{eq:consWT}
\lim_{n\to\infty}\int_{\sigma} f(\sigma)\diff \big(\alpha_n(k)\big)(\sigma)=\int_{\sigma} f(\sigma)\diff \big(\alpha(k)\big)(\sigma)\;. 
\end{equation}
\end{claim:lem:TreeFunctionsContinuous1}
\begin{proof}
Indeed, since $\alpha_1,\alpha_2,\ldots$ converge in $\P^c_k$ to $\alpha$, then in particular $\alpha_1(k),\alpha_2(k),\ldots$ converge in $L^c_k$ to $\alpha(k)$. By~\eqref{notation:Llittle}, we have that $L^c_k=\M_{\leq c}(\P_{k-1}^{c})$ is equipped with the weak topology, of which~\eqref{eq:consWT} is the defining property (c.f.~\eqref{eq:defWeakConv}).
\end{proof}

Observe that the total mass function $D:\P_k^c\to [0,c]$ is continuous.
Hence, it is enough to focus on functions $f_{T,k}:\P_k^c\to \R$, and the continuity of the functions $\overline{f_{T,k}}$ will follow. For each $k$, the function $f_{\square,k}$ is continuous by its definition given in Definition~\ref{def_treefunc}\ref{en:froot}. For more complicated trees, we proceed by induction on $k$. Reflecting Definition~\ref{def_treefunc}, we distinguish whether the root of $T$ has degree~1 or more. In the former case we can write $T=F^\uparrow$. We apply Definition~\ref{def_treefunc}\ref{new_root},
$$
f_{T,k}(\beta)
=
\exp(-D(\beta))\cdot\int_{\sigma\in \P_{k-1}^{c}} f_{F,k-1}(\sigma) \diff\big(\beta(k)\big)(\sigma)\;.
$$
The term $\exp(-D(\beta))$ is continuous in $\beta$ as we noted above. The term $\int_{\sigma\in \P_{k-1}^{c}} f_{F,k-1}(\sigma) \diff\big(\beta(k)\big)(\sigma)$ is continuous in $\beta$ by Claim~\ref*{lem:TreeFunctionsContinuous}.A and the inductive assumption about the continuity of the function $f_{F,k-1}$.

If the root of $T$ has degree more than~1, then we use Definition~\ref{def_treefunc}\ref{multiplication}. It tells us that we can write $f_{T,k}$ as a product of a constant, a power of $\exp(-D(\beta))$, and a product of functions $f_{T_i,k}$ where the degree of the root in each $T_i$ is~1. In particular, the functions $f_{T_i,k}$ are continuous as per the previous case. Hence, $f_{T,k}$ is continuous.
\end{proof}

\subsection{Finding a tree separating iterated degree sequences}\label{ssec:FindingSepTree}

The following proposition is the key step to proving the $(\Rightarrow)$ direction of Theorem~\ref{thm:mainBJR}. Its proof is given in Section~\ref{ssec:proofSeparation}.
\begin{prop}\label{prop:keyseparation}
	For each $n\in \N_0$ the sets $\F_n$ and $\overline{\F_n}$ both separate points of $\P_n^c$.
\end{prop}

To finish the proof of Theorem~\ref{thm:mainBJR}, the following concept will be useful. Given $c,d>0$, $n\in \N_\infty$ and a measure $\nu\in \M_{\le d}(\P_n^c)$, its \emph{exponential tilting}\label{notation:exponentialtilting} is a measure $\widehat{\nu}$ on $\P_n^c$ defined for each Borel set $S\subset\P_n^c$ by
\[
\widehat{\nu}(S):=\int_{\alpha\in S} \exp(-D(\alpha)) \diff \nu(\alpha)
\;.
\]

The next lemma asserts that exponential tilting is injective on $\M_{\leq d}(\P_n^c)$.
\begin{lem}\label{lem:exptilting}
Suppose that $c,d>0$, $n\in \N_\infty$ and that $\nu_1,\nu_2\in\M_{\leq d}(\P_n^c)$ are two different measures. Then their exponential tiltings $\widehat{\nu_1}$ and $\widehat{\nu_2}$ are different measures in $\M_{\leq d}(\P_n^c)$.
\end{lem}
\begin{proof}
The standard rules for Radon--Nikodym derivatives imply that exponential tilting is invertible, that is, for $i=1,2$ we have $\nu_i(S):=\int_{\alpha\in S} \exp(D(\alpha)) \diff \widehat{\nu_i}(\alpha)$ for each measurable subset $S\subset \P_n^c$. Since $\nu_1$ and $\nu_2$ are different measures, this implies that $\widehat{\nu_1}$ and $\widehat{\nu_2}$ are different. Furthermore, for $i=1,2$ the fact that $\exp(-D(\alpha))\leq 1$ for each $\alpha\in \P_n^c$ yields
\[\widehat{\nu_i}(\P_n^c) = \int_{\alpha\in \P_n^c} \exp(-D(\alpha))\diff\nu_i(\alpha) \leq \int_{\alpha\in \P_n^c} 1 \diff\nu_i(\alpha) = \nu_i(\P_n^c) \leq d\;,\]
so we have $\widehat{\nu_i}\in\M_{\leq d}(\P_n^c)$ as required.
\end{proof}

\begin{proof}[Proof of the $(\Rightarrow)$ direction of Theorem~\ref{thm:mainBJR}]
Suppose that $U$ and $W$ are not fractionally isomorphic, that is, we have $\nu_U\neq \nu_W$. Thus, there is an $n\in \N$ such that $\nu_{U,n}\neq \nu_{W,n}$, where $\nu_{W,n}:=(p_{n,\infty}^c)_*\nu_W$. Fix that $n$ and let $\widehat{\nu_{U,n}}$ and $\widehat{\nu_{W,n}}$ be exponential tiltings of $\nu_{U,n}$ and $\nu_{W,n}$ respectively. By Lemma~\ref{lem:exptilting} we have $\widehat{\nu_{U,n}}\neq \widehat{\nu_{W,n}}$.

The family $\overline{\F_n}$ consists of continuous functions (Lemma~\ref{lem:TreeFunctionsContinuous}), is multiplicative up to constants (by Lemma~\ref{lem:FMult}), contains the constant-1 function (by Definition~\ref{def_treefunc}\ref{f_bar}), and separates points of $\P_n^c$ (by Proposition~\ref{prop:keyseparation}). Therefore, the Stone--Weierstraß Theorem (Corollary~\ref{cor:StoneWeierstrass}) tells us that there exists a finite rooted tree $T$ of height at most $n$ such that
\begin{equation} \label{eq:intdiff}
\int_{\alpha\in\P_n^c}\overline{f_{T,n}}(\alpha)\diff\widehat{\nu_{U,n}}(\alpha) \neq \int_{\alpha\in\P_n^c}\overline{f_{T,n}}(\alpha)\diff\widehat{\nu_{W,n}}(\alpha)\;.
\end{equation}
By unpacking the definitions of $\overline{f_{T,n}}$ and $\widehat{\nu_{U,n}}$, and then applying Lemma~\ref{lem:generatingBranchings} with $f_{T,n}$, we obtain
\begin{equation*}
\Probability[(\fX_U)_{\restriction n}\cong T] = \int_{\alpha\in\P^c} f_{T,n}(p^c_{n,\infty}(\alpha))\diff\nu_U(\alpha) =\int_{\sigma\in\P^c_n} f_{T,n}(\sigma)\diff\nu_{U,n}(\sigma) = \int_{\sigma\in\P_n^c}\overline{f_{T,n}}(\sigma)\diff\widehat{\nu_{U,n}}(\sigma)\;.
\end{equation*}
Analogously, we have $\Probability[(\fX_W)_{\restriction k}\cong T] = \int_{\alpha\in\P_n^c}\overline{f_{T,n}}(\sigma)\diff\widehat{\nu_{W,n}}(\sigma)$. Hence, by~\eqref{eq:intdiff} we have $\Probability[(\fX_U)_{\restriction k}\cong T] \neq \Probability[(\fX_W)_{\restriction k}\cong T]$, completing the proof.
\end{proof}

\subsection{Proof of Proposition~\ref{prop:keyseparation}}\label{ssec:proofSeparation}
We shall prove by induction on $n\in\N_0$ that $\F_n$ and $\overline{\F_n}$ both separate any pair $\alpha\neq\beta\in\P_n^c$. By consistency, $\alpha\neq\beta$ implies $\alpha(n) \neq \beta(n)$.

The statement is trivial for $n=0$ because $\alpha(0) = \beta(0) = \P_0 = \{\OneElementSpace\}$ has only one element. Now suppose that the claim has been shown for some $n\in\N_0$. Let $\alpha\neq \beta\in\P_{n+1}^c$. Let $\alpha'$ and $\beta'$ be the exponential tiltings of $\alpha(n+1)$ and $\beta(n+1)$ respectively. By Lemma~\ref{lem:exptilting}, $\alpha'$ and $\beta'$ are different measures in $\M_{\leq c}(\P_n^c)$.

By the inductive hypothesis, $\overline{\F_n}$ separates elements of $\P_n^c$. By Lemma~\ref{lem:FMult}, $\overline {\F_n}$ is multiplicative up to constants. Clearly $\overline{\F_n}$ contains the constant-1 function. Hence, the Stone--Weierstraß Theorem (Corollary~\ref{cor:StoneWeierstrass}) implies that there is a $T\in \T_{n}$ such that
\[\int_{\P_n^c}\overline{f_{T,n}}(\gamma)\diff\alpha'(\gamma)\neq\int_{\P_n^c}\overline{f_{T,n}}(\gamma)\diff\beta'(\gamma)\;.\]
We have $\overline{f_{T,n}}(\gamma) = \exp(D(\gamma))f_{T,n}(\gamma)$ by Definition~\ref{def_treefunc}\ref{f_bar}. Together with the definitions of $\alpha'$ and $\beta'$, this yields $$\int_{\P_n^c}f_{T,n}(\gamma)\exp(D(\gamma))\exp(-D(\gamma))\diff\alpha(n+1)\neq\int_{\P_n^c}f_{T,n}(\gamma)\exp(D(\gamma))\exp(-D(\gamma))\diff\beta(n+1),$$ 
    which cancels to 
    \begin{equation}\label{eq:domu}
    \overline{f_{T^\uparrow,n+1}}(\alpha)\overset{Def.\  \ref{def_treefunc}\ref{new_root}, \ref{f_bar}}{=}\int_{\P_n^c}f_{T,n}(\gamma)\diff\alpha(n+1)\neq\int_{\P_n^c}f_{T,n}(\gamma)\diff\beta(n+1)=\overline{f_{T^\uparrow,n+1}}(\beta)\;.
    \end{equation}
    Thus, $\overline{\F_{n+1}}$ separates $\alpha$ and $\beta$.

    To complete the inductive step, it remains to show that $\F_{n+1}$ also separates $\alpha$ and $\beta$. We have from before a rooted tree $T\in\T_{n}$ such that $\overline{f_{T^\uparrow,n+1}}(\alpha)\neq \overline{f_{T^\uparrow,n+1}}(\beta)$. We now distinguish two cases depending on whether $f_{T^\uparrow,n+1}(\alpha) \neq f_{T^\uparrow,n+1}(\beta)$ or $f_{T^\uparrow,n+1}(\alpha) = f_{T^\uparrow,n+1}(\beta)$. Obviously, $\F_{n+1}$ separates $\alpha$ and $\beta$ in the former case. So, let us assume the latter case. For all $\gamma\in\P^c$ we have 
    \[f_{T^\uparrow\oplus T^\uparrow,n+1}(\gamma) \overset{Def.\  \ref{def_treefunc}\ref{multiplication}, \ref{f_bar}}{=} e_{T^\uparrow\oplus T^\uparrow}f_{T^\uparrow,n+1}(\gamma)\overline{f_{T^\uparrow,n+1}}(\gamma)\;.\] 
    Using this for $\gamma=\alpha,\beta$, it follows that $f_{T^\uparrow\oplus T^\uparrow,n+1}(\alpha) \neq f_{T^\uparrow\oplus T^\uparrow,n+1}(\beta)$. This shows that $\F_{n+1}$ separates $\alpha$ and $\beta$, completing the proof.

\section{More tools for the proof of Theorem~\ref{thm:mainUST}}

\subsection{Random walks on general state spaces}\label{ssec:RWbackground}
One of the key techniques in our proof is random walks on graphons. As far as we could find, the only work dealing with this topic is \cite{MR4334549}, which focuses on linking random walks on dense graphs to those on their limiting graphon. In this paper, we do not need this link and study just the graphon setting. To this end, standard theory of discrete-time Markov chains on general probability spaces suffices. We use~\cite{MeynTweedieBook} as our reference. 

We begin with some definitions. Let $(X,\B)$ be a standard Borel space. A function $P\colon X\times\B\to[0,1]$ is a \emph{transition probability kernel} if for each $x\in X$ the measure $P_x$ on $\B$ given by $P_x(A) = P(x,A)$ is a probability measure and for each $A\in\B$ the function $f_A\colon X \to [0,1]$ given by $f_A(x) = P(x,A)$ is a measurable function on $X$. A time-homogeneous \emph{Markov chain} with transition probability kernel $P$ and initial distribution $\chi$ is a stochastic process $\Phi=\{\Phi(n)\}_{n\in\N_0}$ defined on $\left(\prod_{i\in\N_0}X,\bigotimes_{i\in\N_0}\B\right)$ whose finite dimensional distributions are such that for all $n\in\N$ and measurable $A_0,\dots,A_n\subset X$ we have
\[ \Probability\left[(\Phi(0)\in A_0,\dots,\Phi(n)\in A_n\right] =\int_{y_0\in A_0}\cdots\int_{y_{n-1}\in A_{n-1}} P(y_{n-1},A_n)\ \diff P_{y_{n-2}}(y_{n-1})\dots\diff P_{y_0}(y_1)\ \diff\chi(y_0)\;. \]

Given a transition probability kernel $P$, we recursively define the corresponding $n$-step transition probability kernel $P^n\colon X\times\B\to[0,\infty)$ as follows. Set $P^0(x,A) = \One_A(x)$ for all $x\in X$ and $A\in\B$. For $n\in\N$ we recursively define
\begin{equation} \label{eq:stepTPKrecursive}
    P^n(x,A) = \int_{y\in X}P^{n-1}(y,A)\diff P_x(y)
\end{equation}
for all $x\in X$ and $A\in\B$. We write $P^n_x$ for the probability measure on $\B$ given by $P^n_x(A):=P^n(x,A)$.

A key plank of Markov chain convergence is the property of irreducibility. Given a probability measure $\phi$ on $\B$, a Markov chain $\Phi$ is \emph{$\phi$-irreducible} if for every $x\in X$ and every set $A\in\B$ with $\phi(A)>0$ there is a positive probability that $\Phi$ started at $x$ satisfies $\Phi(n)\in A$ for some $n\in\N$. We call $\phi$ an \emph{irreducibility measure} for $\Phi$. The property of $\phi$-irreducibility is somewhat weak and fails to fully represent the usual notion of irreducibility for Markov chains on discrete state spaces, so we give a further definition to fully capture the range of behavior of Markov chains. We say that a probability measure $\psi$ on $\B$ is a \emph{full irreducibility measure} for $\Phi$ if there is an irreducibility probability measure $\phi$ for $\Phi$ such that, writing $P$ for the transition probability kernel of $\Phi$, for all $A\in\B$ we have
\[\psi(A) := \int_X \sum_{j=1}^\infty 2^{-j-1}P^j(y,A) \diff\phi(y)\;.\]
The concept of a full irreducibility probability measure is a cornerstone of the theory of Markov chains on general state spaces and provides the right setting for the study of Markov chain convergence. Indeed, it is shown in~\cite[Proposition 4.2.2]{MeynTweedieBook} that full irreducibility probability measures have some additional properties that fully capture the range of behavior of Markov chains. That said, we need full irreducibility probability measures solely to enable the application of Theorem~\ref{thm:MCconverge} and will not directly utilize or mention any of the aforementioned properties. A more detailed discussion can be found in~\cite[Section 4.2]{MeynTweedieBook}.

A probability measure $\pi$ on $\B$ is \emph{invariant} if for all $A\in\B$ we have $\pi(A) = \int P(x,A)\diff\pi(x)$. A chain $\Phi$ is \emph{positive recurrent} if it has a full irreducibility probability measure and an invariant probability measure.

Let $\Phi$ be a Markov chain with transition probability kernel $P$ and full irreducibility probability measure $\psi$, and let $p\in\N$. A \emph{$p$-cycle} for $\Phi$ is a collection of disjoint sets $D_{p+1}=D_1,\dots,D_p\in\B$ such that
\begin{itemize}
    \item for all $i\in[p]$ and all $x\in D_i$ we have $P(x,D_{i+1})=1$, and
    \item the set $N=X\sm(\bigcup_{i\in[p]}D_i)$ satisfies $\psi(N)=0$.
\end{itemize}
The \emph{period} of $\Phi$ is the largest $p\in\N$ for which $\Phi$ has a $p$-cycle. It follows from Theorems 5.2.2 and 5.4.4 in~\cite{MeynTweedieBook} that there exists $p\in\N$ such that $\Phi$ has a $p$-cycle and any $q\in\N$ for which $\Phi$ has a $q$-cycle must be a divisor of $p$; in particular, the period of $\Phi$ is well-defined.

A \emph{signed measure} on $\B$ is a countably additive set function $\eta\colon \B \to (-\infty,\infty)$; in particular, this extends the usual notion of a measure by permitting negative values. The \emph{total variation norm} of a signed measure $\eta$ on $\B$ is defined by
\begin{equation} \label{eq:TVNorm}
    \|\eta\|_{TV} = \sup_{f:X\rightarrow [-1,1]}\left|\int f(x)\diff\eta(x)\right| = \sup_{A\in\B}\eta(A) - \inf_{A\in\B}\eta(A)\;.
\end{equation}

The main result from the theory of Markov chains we will need is the result about convergence of a Markov chain to the stationary distribution. Recall that in the more familiar setting of Markov chains on a finite or countable state space, in the aperiodic and irreducible case, the distribution of a random walk after $n$ steps converges to a uniquely defined stationary distribution in the total variation distance as $n\to\infty$, independently of the initial distribution. The following is a version of the convergence result for Markov chains on a general state space which handles periodicity by taking an average over each periodic interval.
\begin{theo}[Theorem 13.3.4(ii) in~\cite{MeynTweedieBook}] \label{thm:MCconverge}
    Suppose that $\Phi$ is a positive recurrent Markov chain with period $p\in\N$ with an invariant probability measure $\pi$. Then for every initial probability distribution $\chi$ which is absolutely continuous with respect to $\pi$ we have
    \[ \left\|\int \frac{1}{p}\sum_{r=0}^{p-1} P^{np+r}(x,\cdot)\  \diff\chi(x)-\pi(\cdot)\right\|_{TV} \to 0 \textrm{ as } n \to \infty. \]
\end{theo}
\subsection{Markov renormalization $W^\dag$}\label{ssec:MarkovRenormalization}
In this section, we introduce a way to transform a kernel $W$ into an akernel $W^\dag$, which we call the Markov renormalization. The Markov renormalization will play a key role in the proof of Theorem~\ref{thm:mainUST}. The defining formula comes from the density appearing in Definition~\ref{def:HNTbranching}\ref{en:HNTdesc}.

\begin{defi}\label{def:MarkovRenorm}
Let $W$ be an $L^\infty$-kernel. The \emph{Markov renormalization} of $W$, denoted \label{notation:Markovrenormalization}$W^\dagger$, is an akernel $W^\dagger:X^2\rightarrow[0,\infty)$ defined at $(x,y)\in X^2$ by
    \[
    W^\dagger(x,y):=\frac{W(x,y)}{\deg_W(y)} \quad\text{with the convention $\frac{0}{0}=0$}\;.
    \]
\end{defi}

The reason why we call $W^\dagger$ the Markov renormalization is that for each $y\in X$ we have $\int_{x\in X}W^\dagger(x,y)=1$. Thus we can interpret $W^\dagger(\cdot,y)$ as transition densities from state $y$ of a Markov chain on $X$.

\begin{lem}\label{lem:MarkovBoundedDegree}
Suppose that $W$ is an $L^\infty$-kernel with minimum degree $\delta>0$. Then $W^\dagger$ is an $L^\infty$-akernel with $\|W^\dagger\|_\infty \le \frac{\|W\|_\infty}{\delta}$ and minimum degree at least $\frac{\delta}{\|W\|_\infty}$.
\end{lem}
\begin{proof}
The quantity $\|W\|_\infty$ is an upper bound on the maximum degree of $W$, so for $\mu$-almost every $x\in X$ we have 
\[ \deg_{W^\dagger}(x)=\int_y \frac{W(x,y)}{\deg_W(y)}\diff\mu(y) \ge \int_y \frac{W(x,y)}{\|W\|_\infty} \diff\mu(y) \ge \frac{\delta}{\|W\|_\infty} \;.\]
This yields the desired lower bound on the minimum degree. The upper bound on $\|W^\dagger\|_\infty$ follows from the observation that we have $W^\dagger(x,y) = 0$ or $W^\dagger(x,y) = \frac{W(x,y)}{\deg_W(y)} \le \frac{\|W\|_\infty}{\delta}$ for $(\mu\times\mu)$-almost every $(x,y)\in X\times X$.
\end{proof}

The following fact follows directly by comparing the definitions of $\fU^-_W$ and $\fX_{W^\dagger}$.
\begin{fact}\label{fact:CorrespondendceDownUST}
Suppose that $W$ is an $L^\infty$-kernel with positive minimum degree. Then the branching process $\fU^-_W$ defined in Figure~\ref{fig:HNTbranching} has the same distribution as $\fX_{W^\dagger}$.
\end{fact}

\subsubsection{Extinction}\label{sssec:Extinction}
The main result of this section is Lemma~\ref{lem:extinct}, which asserts that for each kernel $W$ the branching process $\fX_{W^\dagger}$ goes extinct almost surely.

\begin{lem}\label{lem:zeroIsPossible}
Given $\rho>0$ there exists a number $g_\rho>0$ such that the following holds.
Suppose that $W$ is a kernel satisfying $\frac{\degmin(W)}{\degmax(W)}\ge \rho$. Then for each $x\in X$ the probability that a particle type $x$ in the branching process $\fX_{W^\dagger}$ has no children is at least $g_\rho$.
\end{lem}
\begin{proof}
Indeed, the number of children of such a particle has distribution $\Poisson(\int_{z\in X} \frac{W(x,z)}{\deg_W(z)})$. The lemma follows by noting that $\int_{z\in X} \frac{W(x,z)}{\deg_W(z)}\le \frac{1}{\rho}$.
\end{proof}

The next lemma tells us that in each generation the expected number of particles is~1. It actually gives a more precise description of the distribution of the particles.
\begin{lem}\label{lem:eachgenerationone}
Suppose that $W$ is a nondegenerate $L^1$-kernel. Consider the branching process $\fX_{W^{\dagger}}$. For $k\in \N_0$ and $A\in \B$, let $Y_{k,A}$ be the number of particles in generation $k$ and of type $A$. Then we have $\Expectation[Y_{k,A}]=\mu(A)$.
\end{lem}
\begin{proof}
We prove the claim by induction. The base case $k=0$ is clear since the root of $\fX_{W^{\dagger}}$ is chosen according to $\mu$. Let us move to the induction step $k\to k+1$. The number of particles in generation $k+1$ of type in $A$ which are born from a given particle in generation $k$ of type $x\in X$ has distribution $\Poisson\left(\int_{z\in A} \frac{W(x,z)}{\deg_W(z)}\right)$. In particular, the expected number of particles in generation $k+1$ of type in $A$ which are born from a given particle in generation $k$ of type $x$ is $\int_{z\in A} \frac{W(x,z)}{\deg_W(z)}$. By the induction hypothesis, the expected number of particles in generation $k+1$ of type in $A$ is obtained by integrating the above uniformly over the choice of $x$,
\begin{align*}
    \Expectation[Y_{k+1,A}]&=\int_{x\in X}\int_{z\in A} \frac{W(x,z)}{\deg_W(z)}\ \diff\mu(z)=\int_{z\in A}\frac1{\deg_W(z)}\int_{x\in X} W(x,z)\ \diff\mu(z)\\
    &=\int_{z\in A}\frac1{\deg_W(z)}\cdot \deg_W(z)\ \diff\mu(z)=\mu(A)\;,
\end{align*}
as was needed.
\end{proof}

\begin{lem}\label{lem:extinct}
Suppose that $W$ is a bounded-degree kernel with positive minimum degree. Then the branching process $\fX_{W^\dagger}$ goes extinct almost surely.
\end{lem}
\begin{proof}
    For $k,j\in\N_0$, let $q_{k,j}$ be the probability that in $\fX_{W^\dagger}$ there are exactly $j$ particles in the $k$-th generation. By the law of total probability, we have for each $k\in\N$ that
    \begin{equation}\label{eq:total_prob}
        \sum_{j=0}^\infty q_{k,j}=1\;.
    \end{equation}
We have from Lemma~\ref{lem:eachgenerationone} that for each $k\in\N$, $\sum_{j=0}^\infty j q_{k,j}=1$. In particular, for every $J\in \N$,
\begin{equation}\label{eq:tail}
    \sum_{j=1}^J q_{k,j}>1-q_{k,0}-\tfrac{1}{J}\;.
\end{equation}

The sequence $(q_{k,0})_{k=0}^\infty$ is nondecreasing and hence has a limit $L=\lim_{k\to\infty} q_{k,0}$. The lemma amounts to proving that $L=1$. Assume for contradiction that $L<1$. Let
\begin{equation}\label{eq:IloveJ}
J:=\left\lceil\frac{20}{1-L}\right\rceil\;.
\end{equation}
Let $g>0$ be from Lemma~\ref{lem:zeroIsPossible} for our kernel $W$. Define $g^*:=g^J$. Note that for each $j\le J$, $g^*$ is a lower bound for going extinct in generation $k+1$ ($k$ arbitrary) given that in generation $k$ there were $j$ particles. 

Let $K\in\N$ be such that for every $k\ge K$ we have that 
\begin{equation}\label{eq:inthebox}
    q_{k,0}\in\left[L-\frac{g^*(1-L)}{10},L\right]\;.
\end{equation}

Note that we have~0 particles in a generation $k+1$ if there were 0 particles in generation $k$ or the number of particles in generation $k$ was between $1$ and $J$ and none of these particles had any children (this is obviously not if and only if). By the definition of $g^*$ we have
\begin{align}
\begin{split}\label{eq:split}
    q_{k+1,0}&\ge q_{k,0}+g^*(q_{k,1}+q_{k,2}+\ldots+q_{k,J})\\
    \JUSTIFY{by~\eqref{eq:tail}}&\ge q_{k,0}+g^*\left(1-q_{k,0}-\tfrac1J\right)=(1-\tfrac1J)g^*+(1-g^*)q_{k,0}\;.
    \end{split}
\end{align}
Take $k\ge K$ and use~\eqref{eq:inthebox} for $q_{k,0}$ and $q_{k+1,0}$ in~\eqref{eq:split}.
\begin{equation*}
    L\ge (1-\tfrac1J)g^*+(1-g^*)\left(L-\frac{g^*(1-L)}{10}\right)\;.
\end{equation*}
Pedestrian manipulations (in which we divide both sides by $g^*>0$) give
\begin{align*}
    0 \ge 1-\frac1J-\frac{1-L}{10}-L+\frac{g^*(1-L)}{10}\geByRef{eq:IloveJ}
    1-\frac{1-L}{20}-\frac{1-L}{10}-L=0.85-0.85L>0\;,
\end{align*}
a contradiction.
\end{proof}

\subsection{Weak isomorphism of kernels}
Last, we will use the notion of weak isomorphism of $L^\infty$-kernels. While the original definition goes via homomorphism densities, here we will rather use another property. As is usual in the theory of graph limits, given a kernel $Z:X^2\to \R$ and measure preserving map $\pi:X\to X$, we write $Z^\pi$ for the kernel defined by $Z^\pi(x,y)=Z(\pi(x),\pi(y))$.
\begin{lem}[Corollary~10.35(a) in~\cite{MR3012035}]\label{lem:weakIso}
Suppose that $Z_1$ and $Z_2$ are two weakly isomorphic $L^\infty$-kernels. Then there exist measure preserving maps $\pi_1,\pi_2:X\to X$ such that $Z_1^{\pi_1}(x,y)=Z_2^{\pi_2}(x,y)$ for almost every $(x,y)\in X^2$.
\end{lem}

\subsection{More on fractional isomorphism and sub-sigma-algebras}\label{ssec:MorePrelFI}

Let $(X,\B)$ be a standard Borel space endowed with a Borel probability measure $\mu$. Let $W$ be an $L^2$-kernel on $X$. The corresponding Hilbert--Schmidt integral operator $T_W\colon L^2(X,\mu)\to L^2(X,\mu)$ is defined by\begin{equation}\label{eq:IntegralOperator}
T_W(f)(x) = \int_X W(x,y)f(y)\diff\mu(y)\;.    
\end{equation}

 We say that $\C\subset\B$ is a \emph{$\mu$-relatively complete sub-sigma-algebra} of $\B$ if it is a sub-sigma-algebra such that every $Z\in\B$ satisfies $Z\in\C$ whenever there is $Z_0\in\C$ such that $\mu(Z\triangle Z_0)=0$. Let $\Theta_\mu$ be the collection of all $\mu$-relatively complete sub-sigma-algebras of $\B$. For $\cX\subset\B$ we write $\langle\cX\rangle$ for the unique relatively complete sub-sigma-algebra generated by $\cX$; this notion is well-defined as a consequence of~\cite[Claim~5.4]{Grebik2022}.

For an $L^\infty$-akernel $W$, we recall from~\cite[Section~5]{Grebik2022} the construction of a sequence \label{not:CW}$(\C^W_n)_{n\in\N_0}$ of non-decreasing sigma-algebras which are all subalgebras of $\B$. We begin with some necessary definitions. Let $\C\in\Theta_\mu$. For $\D\in\Theta_\mu$ we say that $(\C,\D)$ is \emph{$W$-invariant} if $T_W(L^2(X,\C,\mu))\subset L^2(X,\D,\mu)$; in particular, we say that $\C$ is \emph{$W$-invariant} if $(\C,\C)$ is $W$-invariant. Let $Z_{\C}$ be the collection of $\D\in\Theta_\mu$ such that $(\C,\D)$ is $W$-invariant. Write $m(\C)$ for the collection of sets $S\in\B$ such for all $\D\in Z_{\C}$ we have $S\in\D$. Now the collection $(\C^W_n)_{n\in\N_0}$ of sigma-algebras is given follows. Let $\C^W_0 = \langle\{\varnothing,X\}\rangle$ and inductively define $\C^W_n = m(\C^W_{n-1})$. Finally, define $\C(W) = \left\langle\bigcup_{n\in\N_0}\C^W_n\right\rangle$ as the unique relatively complete sigma-algebra generated by the union of these sigma-algebras.

The following lemma represents a straightforward generalization of~\cite[Corollary 6.7]{Grebik2022} and describes the relationship between the function $i_W$ and the sub-sigma-algebra $\C(W)$. Observe that the lemma is stated here for $L^\infty$-akernels with bounded maximum degree, while~\cite[Corollary 6.7]{Grebik2022} is stated for akernels taking values in $[0,1]$. We remark that that the key property utilized in the proof of~\cite[Corollary 6.7]{Grebik2022} is that akernels taking values in $[0,1]$ have bounded $L^\infty$-norm and bounded maximum degree, so in particular the proof easily generalizes to all $L^\infty$-akernels with bounded maximum degree. As such, we shall omit the proof of the following lemma.

\begin{lem}\label{lem:generate_CW}
    Suppose $W$ is an $L^\infty$-akernel with maximum degree $c$. Then $i_W$ is measurable and we have \[\langle\{i_W^{-1}(A)\mid A\in\B(\P^c)\}\rangle=\C(W),\] i.\,e., the minimum relatively complete sub-sigma-algebra of $\B$ that makes the map $i_W$ measurable is $\C(W)$.
\end{lem}

Let us also mention here that we will need to borrow further results from~\cite[Section 6]{Grebik2022} in Section~\ref{ssec:ProofofLemmaFracIsoProjFracIso}. As with the lemma above, these results will be stated for $L^\infty$-akernels with bounded maximum degree, but they will correspond to results in~\cite[Section 6]{Grebik2022} stated for akernels taking values in $[0,1]$. In a similar vein, their proofs will easily generalize and so we will omit them too.

\subsection{Conditional expectation}

Let $R$ be an integrable random variable on a standard Borel space $(X,\B)$ equipped with a Borel probability measure $\mu$. Recall that $\Theta_\mu$ is the collection of all $\mu$-relatively complete sub-sigma-algebras of $\B$. Let $\C\in\Theta_\mu$. The conditional expectation $\Expectation[R\vert\C]$ of $R$ given $\C$ is a $\C$-measurable integrable random variable on $(X,\B,\mu)$ such that for all $A\in\C$ we have
\begin{equation} \label{eq:condexpec-defining}
    \int_A \Expectation[R\vert\C] \diff\mu = \int_A R \diff\mu\;.
\end{equation}

The following is a useful standard fact about conditional expectations.

\begin{fact} \label{fact:condexpec-knownfactor}
    Let $R$ and $S$ be integrable random variables on a standard Borel space $(X,\B)$ equipped with a Borel probability measure $\mu$. Let $\C\in\Theta_\mu$. Suppose that $RS$ is integrable and $R$ is $\C$-measurable. Then $\Expectation[RS\vert\C] = R\Expectation[S\vert\C]$ holds $\mu$-almost everywhere.
\end{fact}

\section{A factorization result for fractional isomorphism}\label{sec:factorFI}
This section deals with Theorem~\ref{thm:factorFI}, which is a factorization result for fractional isomorphism of disconnected kernels. It says that two kernels $U$ and $W$ are fractionally isomorphic if and only if for each connected kernel $\Gamma$ the total measure of connected components which are fractionally isomorphic to $\Gamma$ is the same in both $U$ and $W$. This result allows us to extend Theorem~\ref{thm:mainUST} to disconnected kernels in Theorem~\ref{thm:mainUSTGeneral}. We also give a factorization result for finite graphs.

To state the kernel result and also for subsequent parts of the paper, it is useful to introduce the concept of subgraphons or subkernels. Given a kernel $U$ on $X$ and a measurable subset $Y\subset X$ with $\mu(Y)>0$, we define \label{notation:restrictedgraphon}$U_{\restriction Y}$ as the restriction of $U$ to $Y\times Y$. We view $U_{\restriction Y}$ again as a kernel, that is, we equip $Y$ with the probability measure $\frac{\mu(Y\cap \cdot)}{\mu(Y)}$. Next, we introduce a rescaling \label{notation:restrictedrescaledgraphon}$U\llbracket Y\rrbracket:Y\times Y\rightarrow \R$ defined by $U\llbracket Y\rrbracket(x,y)=\mu(Y)U_{\restriction Y}(x,y)$. This rescaling is chosen so that if $Y$ is a connected component of $U$ (or a union of connected components) then for every $y\in Y$ we have $\deg_{U\llbracket Y\rrbracket}(y)=\deg_U(y)$. More generally, for each $y\in Y$ we have $i_{U\llbracket Y\rrbracket}(y)=i_U(y)$ and for every $A\in\B(\P^c)$ we have
\begin{equation}\label{eq:transferOfDegreeMeasureIntoSubgraphs}
\mu(Y)\nu_{U\llbracket Y\rrbracket}(A)=\mu\left(i_U^{-1}(A)\cap Y\right)\;.
\end{equation}

The main result now reads as follows.
\begin{theo}\label{thm:factorFI}
Suppose that $U$ and $W$ are two $L^\infty$-kernels. Let $X=\Lambda_0\sqcup\bigsqcup_{i=1}^{N_U}\Lambda_i$ and $X=\Omega_0\sqcup\bigsqcup_{i=1}^{{N_W}}\Omega_i$ be decompositions into the connected components of $U$ and $W$ respectively, with $\Lambda_0$ and $\Omega_0$ being the respective sets of isolated elements. Here $N_U,N_W\in \N\cup\{\infty\}$ are not necessarily equal. Then $U$ and $W$ are fractionally isomorphic if and only if for each connected kernel $\Gamma$ we have
\begin{equation}\label{eq:factorFI}
\sum_{i\ge 1:\text{$U\llbracket\Lambda_i\rrbracket$ is frac.\ iso.\ to $\Gamma$}}\mu(\Lambda_i)
=
\sum_{i\ge 1:\text{$W\llbracket\Omega_i\rrbracket$ is frac.\ iso.\ to $\Gamma$}}\mu(\Omega_i)
\;.    
\end{equation}
\end{theo}

Note that if~\eqref{eq:factorFI} holds, then $\mu(\Lambda_0)=\mu(\Omega_0)$. Indeed, there are at most countably many mutually fractionally nonisomorphic akernels $\Gamma$ for which the value in~\eqref{eq:factorFI} is positive. The sum of these values is, when reasoning from the left-hand side, $1-\mu(\Lambda_0)$, and when reasoning from the right-hand side, $1-\mu(\Omega_0)$.

Let us now turn to finite graphs. We emphasize that the discussion and factorization result below are not needed for our results about branching processes. Nonetheless, we view them as an important contribution to the theory of fractional isomorphism of graphs. The obvious counterpart to Theorem~\ref{thm:factorFI} fails for fractional isomorphism. To see this, let $G$ be the disjoint union of a 6-cycle and a 4-cycle, and let $H$ be a 10-cycle. Then $G$ and $H$ are fractionally isomorphic, yet each of the three components of $G\cup H$ lies in a separate class of fractional isomorphism. There is, however, an interesting factorization result for finite graphs. Furthermore, the proof of this result can be viewed as a finite and easier to digest counterpart to the proof of Theorem~\ref{thm:factorFI}. To state the result, we introduce the term `practional isomorphism', a portmanteau of `proportional' and `fractional'. To the best of our knowledge, this notion is new. Let us give details. Suppose that $G$ is a graph and $\mathcal{P}=\{P_j\}_{j=1}^k$ is a partition of $V(G)$ into nonempty sets. Suppose that $\mathbf{D}=(D_{j,\ell})_{j,\ell=1}^k$ is a matrix of integers, and $\mathbf{p}=(p_j)_{j=1}^k$ is a vector of positive reals. We say that $\mathcal{P}$ is an \emph{equitable partition for template $(\mathbf{D},\mathbf{p})$} if for every $j,\ell\in [k]$ we have
\begin{align}
\label{eq:eq_part1}
    \frac{|P_j|}{v(G)}&=p_j\quad\mbox{, and }\\
\label{eq:eq_part2}
    \deg_G(v,P_\ell)&=D_{j,\ell}\quad\mbox{for every $v\in P_j$.}
\end{align}
Furthermore, we say that $\mathcal{P}$ is the \emph{coarsest equitable partition of $G$} if every other equitable partition refines it. It is easy to check that every graph has a unique (up to the order of the cells) coarsest equitable partition.

We say that two graphs $G$ and $H$ are \emph{practionally isomorphic} if there exist $k\in\N$, and $\mathbf{D}$ and $\mathbf{p}$ as above so that both $G$ and $H$ have an equitable partition for template $(\mathbf{D},\mathbf{p})$. Again, it is well-known and easy to check that this is equivalent to the coarsest equitable partitions of $G$ and of $H$ having the same template.

So, while traditionally, two fractionally isomorphic graphs need to have the same number of vertices, now, for example, any two $3$-regular graphs are practionally isomorphic, even if their orders differ.\footnote{Take $k=1$, $D_{1,1}=3$, $p_1=1$, and trivial equitable partitions $\{V(G)\}$, and $\{V(H)\}$.} 

We remark that two finite graphs $G$ and $H$ being practionally isomorphic is not equivalent to graphon representations $W_G$ and $W_H$ of those two graphs being fractionally isomorphic. For example, graphon representations of two $3$-regular graphs of different orders are not fractionally isomorphic.

By our choice of normalization in~\eqref{eq:eq_part1}, the factorization result for graphs has a particularly elegant form.
\begin{theo}\label{thm:factorFIgraphs}
Suppose that $G$ and $H$ are two graphs. Let $V(G)=A_1\sqcup A_2\sqcup\ldots\sqcup A_n$ and $V(H)=B_1\sqcup B_2\sqcup \ldots\sqcup B_m$ be the connected components of $G$ and of $H$, respectively. Then $G$ and $H$ are practionally isomorphic if and only if for each connected graph $\Gamma$ we have
\begin{equation}\label{eq:factorFIgraphs}
\frac{1}{v(G)}\cdot\sum_{i:\text{$G[A_i]$ is prac.\ iso.\ to $\Gamma$}}|A_i|
=
\frac{1}{v(H)}\cdot\sum_{i:\text{$H[B_i]$ is prac.\ iso.\ to $\Gamma$}}|B_i|
\;.    
\end{equation}
\end{theo}

\subsection{Proof of Theorem~\ref{thm:factorFIgraphs}}\label{ssec:ProofGrDifficult}
We start with the routine $(\Leftarrow)$ direction. The idea is to stitch together corresponding cells of equitable partitions of components of $G[A_i]$ (and later $H[B_i]$) that are practionally isomorphic to the same connected graph $\Gamma$. Let $\mathcal{R}$ be a family  of arbitrary representatives of practional isomorphism classes of those connected graphs which appear among components of $G$ (and equivalently, by \eqref{eq:factorFIgraphs}, among components of $H$). That is, for each $G[A_i]$ there is a unique $\Gamma\in\mathcal{R}$ such that $G[A_i]$ and $\Gamma$ are practionally isomorphic. For each $\Gamma\in\mathcal{R}$, let $(\mathbf{D}^{\Gamma},\mathbf{p}^{\Gamma})$ be the template of the coarsest equitable partition of $\Gamma$. We write $\mathbf{D}^{\Gamma}=(D^{\Gamma}_{j,\ell})_{j,\ell=1}^{k_{\Gamma}}$ and $\mathbf{p}^{\Gamma}=(p^{\Gamma}_j)_{j=1}^{k_{\Gamma}}$. Let us now define a template $(\mathbf{D},\mathbf{p})$ of dimension $k=\sum_{\Gamma\in\mathcal R}k_\Gamma$. Rather than indexing the entries of the template by elements of $\{1,\ldots,k\}$, it is convenient to index them by $S=\{(\Gamma,j):\Gamma\in \mathcal{R},j\in[k_\Gamma]\}$. That is, we define $\mathbf{D}=(D_{s,q})_{s,q\in S}$ and $\mathbf{p}=(p_s)_{s\in S}$ by
\begin{align}
\nonumber    D_{(\Gamma_1,j_1),(\Gamma_2,j_2)}&=D^{\Gamma_1}_{j_1,j_2}&\text{if $\Gamma_1=\Gamma_2$},\\
\nonumber    D_{(\Gamma_1,j_1),(\Gamma_2,j_2)}&=0&\text{if $\Gamma_1\neq\Gamma_2$},\\
\label{eq:GAZ}    p_{(\Gamma,j)}&=\tfrac{\sum_{i:\text{$G[A_i]$ is frac.\ iso.\ to $\Gamma$}}|A_i|}{v(G)}\cdot p^{\Gamma}_j\;.
\end{align}
We now define a partition $\mathcal{P}=(P_s)_{s\in S}$ of $V(G)$ as follows,
\[
P_{(\Gamma,j)}=\bigcup_{i:\text{$G[A_i]$ is frac.\ iso.\ to $\Gamma$}}P^{i}_j\;,
\]
where $\mathcal{P}^i=(P^{i}_j)_{j=1}^{k_\Gamma}$ is a partition of $A_i$ corresponding to template $(\mathbf{D}^{\Gamma},\mathbf{p}^{\Gamma})$, which exists as $\Gamma$ and $G[A_i]$ are practionally isomorphic. It is routine to check that $\mathcal{P}$ is an equitable partition for template $(\mathbf{D},\mathbf{p})$. We can now repeat the whole construction for the graph $H$. The template we create this way is the same; this crucial fact uses~\eqref{eq:factorFIgraphs} in~\eqref{eq:GAZ}. We conclude that $G$ and $H$ are practionally isomorphic.

The bulk of the work is the $(\Rightarrow)$ direction. Assume that $G$ and $H$ are practionally isomorphic.
Let $(\mathbf{D},\mathbf{p})$ be their common template, and let $\mathcal{P}=\{P_j\}_{j=1}^k$ and $\mathcal{Q}=\{Q_j\}_{j=1}^k$ be the corresponding equitable partitions for $G$ and for $H$, respectively.

Consider an auxiliary graph $\Delta$ on vertex set $[k]$. Make a pair $j\ell$ an edge of $\Delta$ if $D_{j,\ell}>0$ (observe that this is equivalent to $D_{\ell,j}>0$). For a set $S\subset [k]$ write $V^{(G)}_{S}=\bigcup_{j\in S} P_j$ and $V^{(H)}_{S}=\bigcup_{j\in S} Q_j$.

\begin{claimFirstA}
For each connected component $A_i$ of $G$ there exists a connected component $C$ of $\Delta$ such $A_i\subset V^{(G)}_C$. Furthermore, for each $j,\ell\in C$, we have $\frac{|A_i\cap P_j|}{|P_j|}=\frac{|A_i\cap P_\ell|}{|P_\ell|}$.
\end{claimFirstA}
Let $C$ be a connected component of $\Delta$ for which $V^{(G)}_C$ contains at least~1 vertex of $A_i$. To prove the first part of the statement, we need to prove that for any other connected component $C'$ of $\Delta$, we have $e_G(A_i\cap V^{(G)}_C,A_i\cap V^{(G)}_{C'})=0$. Indeed, we have 
$$e_G(A_i\cap V^{(G)}_C,A_i\cap V^{(G)}_{C'})\le e_G(V^{(G)}_C,V^{(G)}_{C'})=\sum_{j\in C}\sum_{\ell\in C'}e_G(P_j,P_\ell)=0\;,$$
where the last equality uses~\eqref{eq:eq_part2} and the fact that all the corresponding numbers $D_{j,\ell}$ are~0.

For the furthermore part, observe that it is enough to prove the equality  $\frac{|A_i\cap P_j|}{|P_j|}=\frac{|A_i\cap P_\ell|}{|P_\ell|}$ only for those $j,\ell\in C$ that form an edge in $\Delta$. Using the fact that there are no edges between $A_i\cap P_j$ and $P_\ell\setminus A_i$ we have
\begin{equation}\label{eq:nigh1}    
e_G(A_i\cap P_j,A_i\cap P_\ell)=e_G(A_i\cap P_j, P_\ell)\eqByRef{eq:eq_part2}D_{j,\ell}|A_i\cap P_j|\;,
\end{equation}
and similarly 
\begin{equation}\label{eq:night2}
e_G(A_i\cap P_j,A_i\cap P_\ell)=D_{\ell,j}|A_i\cap P_\ell|\;.
\end{equation}
We also have
\begin{equation}\label{eq:night3}
D_{j,\ell}|P_j|=e_G(P_j,P_\ell)=D_{\ell,j}|P_\ell|\;.
\end{equation}
Putting equations~\eqref{eq:nigh1}, \eqref{eq:night2} and~\eqref{eq:night3} together (using that we may divide by $e_G(P_j,P_\ell)\neq 0$), we get the desired statement.
\begin{claimFirstB}
Suppose that $A_i$ is a connected component of $G$ and $C$ is a connected component in $\Delta$ such $A_i\subset V^{(G)}_C$. Then $G[A_i]$ is practionally isomorphic to $G[V^{(G)}_C]$.
\end{claimFirstB}
\begin{proof}
We show that $\{A_i\cap P_j\}_{j\in C}$ and $\{P_j\}_{j\in C}$ are equitable partitions for $G[A_i]$ and for $G[V^{(G)}_C]$ with the same template.
For~\eqref{eq:eq_part1}, we proceed as follows. Write $\xi$ for the positive real number equal to $\frac{|A_i\cap P_t|}{|P_t|}$ for all $t\in C$; this is well-defined by Claim~\ref*{thm:factorFIgraphs}.A. Now for each $j\in C$ we have
\begin{align*}
    \frac{|A_i\cap P_j|}{|A_i|}=\frac{|A_i\cap P_j|}{\sum_{t\in C}|A_i\cap P_t|}=\frac{\xi|P_j|}{\sum_{t\in C}\xi|P_t|}=\frac{|P_j|}{|V^{(G)}_C|} \;.
\end{align*}
This shows that the corresponding ratios in~\eqref{eq:eq_part1} are equal for the partition $\{A_i\cap P_j\}_{j\in C}$ in the graph $G[A_i]$ and for the partition $\{P_j\}_{j\in C}$ in the graph $G[V^{(G)}_C]$. Let us turn to verifying~\eqref{eq:eq_part2}. Let $j,\ell\in C$ be given. Let us start with the graph $G[A_i]$. Let $v\in A_i\cap P_j$. We have $\deg_{G[A_i]}(v,A_i\cap P_{\ell})=\deg_{G}(v,P_{\ell})=D_{j,\ell}$, where the first equality uses that $A_i$ is a connected component. The calculation for the graph $G[V^{(G)}_C]$ is even more straightforward. Let $v\in P_j$. We have $\deg_{G[V^{(G)}_C]}(v,P_{\ell})=\deg_{G}(v,P_{\ell})=D_{j,\ell}$, as $V^{(G)}_C$ is a union of connected components of $G$ that includes $v$.
\end{proof}
We can now quickly finish the proof. Indeed, Claims~\ref*{thm:factorFIgraphs}.A and~\ref*{thm:factorFIgraphs}.B tell us that the left-hand side of~\eqref{eq:factorFIgraphs} is equal to $\frac{1}{v(G)}\cdot\sum_{C}|V^{(G)}_C|$, where the summation runs over connected components $C$ of $\Delta$ such that $G[V^{(G)}_C]$ is practionally isomorphic to $\Gamma$. By repeating the argument for $H$, we see that the right-hand side of~\eqref{eq:factorFIgraphs} is equal to $\frac{1}{v(H)}\cdot\sum_{C}|V^{(H)}_C|$, where the summation runs over connected components $C$ of $\Delta$ such that $H[V^{(H)}_C]$ is practionally isomorphic to $\Gamma$. The last observation we need to equate both sums, which follows easily from the common template, is that for any nonempty set $S\subset [k]$, the graphs $G[V^{(G)}_S]$ and $H[V^{(H)}_S]$ are practionally isomorphic and $\frac{|V^{(G)}_S|}{v(G)}=\frac{|V^{(H)}_S|}{v(H)}$.

\subsection{Proof of Theorem~\ref{thm:factorFI}}\label{ssec:ProofGraphon}
In kernels, the role of the equitable partition from the proof of Theorem~\ref{thm:factorFIgraphs} is replaced by the sigma-algebra $\C(W)$ which we reintroduced in Section~\ref{ssec:MorePrelFI}. The proof of Theorem~\ref{thm:factorFI} proceeds by adopting this formalism, but the main ideas remain. In particular, there is a clear analogy between Claim~\ref*{thm:factorFIgraphs}.A and Claim~\ref*{thm:factorFI}.A and between Claim~\ref*{thm:factorFIgraphs}.B and Claim~\ref*{thm:factorFI}.B.

We start with the routine $(\Leftarrow)$ direction. The measure $\nu_U$ can be expressed as a convex combination of measures $\{\nu_{U\llbracket\Lambda_i\rrbracket}\}_{i=1}^{N_U}$, 
\begin{equation}\label{eq:rewriteme}
\nu_U=\sum_{i=1}^{N_U}\mu(\Lambda_i)\cdot \nu_{U\llbracket\Lambda_i\rrbracket}+\mu(\Lambda_0)\cdot\Dirac(\oslash)\;,    
\end{equation}
where $\oslash\in\P^0$ is the only element of $\P^0$ (corresponding to the 0-degree). In analogy with the proof of Theorem~\ref{thm:factorFIgraphs}, let $\mathcal{R}$ be an arbitrary family of representatives of fractional isomorphism classes of connected kernels appearing among components of $U$. If $U\llbracket\Lambda_i\rrbracket$ is fractionally isomorphic to $\Gamma\in\mathcal{R}$ then $\nu_{U\llbracket\Lambda_i\rrbracket}=\nu_{\Gamma}$. Hence, we can rewrite~\eqref{eq:rewriteme},
\begin{equation*}
\nu_U=\sum_{\Gamma\in\mathcal{R}}\left(\sum_{i\ge 1:\text{$U\llbracket\Lambda_i\rrbracket$ is frac.\ iso.\ to $\Gamma$}}\mu(\Lambda_i)\right)\cdot \nu_{\Gamma}+\mu(\Lambda_0)\cdot\Dirac(\oslash)\;.    
\end{equation*}
We can now repeat the whole process for the kernel $W$ and end up with 
\begin{equation*}
\nu_W=\sum_{\Gamma\in\mathcal{R}}\left(\sum_{i\ge 1:\text{$W\llbracket\Omega_i\rrbracket$ is frac.\ iso.\ to $\Gamma$}}\mu(\Omega_i)\right)\cdot \nu_{\Gamma}+\mu(\Omega_0)\cdot\Dirac(\oslash)\;.    
\end{equation*}
We use~\eqref{eq:factorFI} to conclude that $\nu_U=\nu_W$. Hence, $U$ and $W$ are fractionally isomorphic.

Let us now turn to the $(\Rightarrow)$ direction. Write $\Delta:=\Expectation[U \vert \C(U)\times\C(U) ]$. By Claim~5.7 in~\cite{Grebik2022} we have
\begin{equation}\label{eq:ExpectationCCBC}
\Delta=\Expectation[U \vert \B(X)\times\C(U) ]=\Expectation[U \vert \C(U)\times \B(X) ]\;.
\end{equation}
\begin{claimSecondA}
For each connected component $\Lambda_i$ of $U$ there exists a connected component $C$ of $\Delta$ such $\Lambda_i\subset C$ (up to a $\mu$-nullset). Furthermore, let $h_i:=\Expectation[\One_{\Lambda_i}\vert \C(U)]$. Then we have that $h_i=\frac{\mu(\Lambda_i)}{\mu(C)}\cdot \One_C$. 
\end{claimSecondA}
\begin{proof}
Let $C$ be a connected component of $\Delta$ such that $\mu(\Lambda_i\cap C)>0$ (such a connected component exists). In order to prove the first part of the statement, we need to prove that for any other connected component $C'$ of $\Delta$, we have $\int_{(\Lambda_i\cap C)\times (\Lambda_i\cap C')}U(x,y)\diff\mu^2(x,y)=0$. Indeed, we have 
\[0 \le \int_{(\Lambda_i\cap C)\times (\Lambda_i\cap C')}U(x,y)\diff\mu^2(x,y) \le \int_{C\times C'}U(x,y)\diff\mu^2(x,y) = \int_{C\times C'}\Delta(x,y)\diff\mu^2(x,y) = 0\;.\]

Let us turn to the furthermore part. We establish the claim in three parts: 
\begin{enumerate}[label=(p\arabic{*})]
\item\label{bonton:1} $h_i$ is constant-0 almost everywhere outside $C$,
\item\label{bonton:2} $h_i$ is constant almost everywhere on $C$, and
\item\label{bonton:3} the value of $h_i$ on $C$ is $\frac{\mu(\Lambda_i)}{\mu(C)}$.
\end{enumerate} 

\emph{Proof of~\ref{bonton:1}.} This is obvious as $\Lambda_i\subset C$.

\emph{Proof of~\ref{bonton:2}.} Suppose for contradiction that $h_i$ is not constant on $C$. That means that there exists $r\ge 0$ such that the sets $R^+=\{x\in C\mid h_i(x)>r\}$ and $R^-=\{x\in C\mid h_i(x)\le r\}$ both have positive measure. As $h_i$ is the conditional expectation with respect to $\C(U)$, it is $\C(U)$-measurable. In particular, $R^+$ and $R^-$ are $\C(U)$-measurable. Note also that since $R^+\sqcup R^-$ is a nontrivial partition of a connected component of $\Delta$, we have
\begin{equation}\label{eq:EverybodyKnows}
    \int_{R^+\times R^-}\Delta(x,y)\diff\mu^2(x,y)>0\;.
\end{equation}

Using that $\Lambda_i$ is a connected component of $U$, we have
\begin{align}
\begin{split}\label{eq:welovepizzaanddoublecounting}
\int_{R^+\times R^-}\One_{\Lambda_i}(x)U(x,y)\diff\mu^2(x,y)&=\int_{R^+\times R^-}\One_{\Lambda_i}(x)\One_{\Lambda_i}(y)U(x,y)\diff\mu^2(x,y)\\
&=\int_{R^+\times R^-}\One_{\Lambda_i}(y)U(x,y)\diff\mu^2(x,y)\;.    
\end{split}
\end{align}
Let us express the first integral of~\eqref{eq:welovepizzaanddoublecounting}.
\begin{align*}
\int_{R^+\times R^-}\One_{\Lambda_i}(x)U(x,y)\diff\mu^2(x,y)&=\int_{R^+}\One_{\Lambda_i}(x)\int_{R^-}U(x,y)\diff\mu(y)\diff\mu(x)\\
\JUSTIFY{\eqref{eq:ExpectationCCBC}, and $R^-\in \C(U)$}
&=\int_{R^+}\One_{\Lambda_i}(x)\int_{R^-}\Delta(x,y)\diff\mu(y)\diff\mu(x)\;.
\end{align*}
Since the function $x\mapsto \int_{R^-}\Delta(x,y)\diff\mu(y)$ is $\C(U)$-measurable by~\eqref{eq:ExpectationCCBC}, by Fact~\ref{fact:condexpec-knownfactor} we have 
$$\int_{R^+}\One_{\Lambda_i}(x)\int_{R^-}\Delta(x,y)\diff\mu(y)\diff\mu(x)=\int_{R^+}h_i(x)\int_{R^-}\Delta(x,y)\diff\mu(y)\diff\mu(x)\;.$$ That is, we have derived $\int_{R^+\times R^-}\One_{\Lambda_i}(x)U(x,y)\diff\mu^2(x,y)=\int_{R^+\times R^-}h_i(x)\Delta(x,y)\diff\mu^2(x,y)$. Similarly, starting with the last integral of~\eqref{eq:welovepizzaanddoublecounting}, we get $\int_{R^+\times R^-}\One_{\Lambda_i}(y)U(x,y)\diff\mu^2(x,y)=\int_{R^+\times R^-}h_i(y)\Delta(x,y)\diff\mu^2(x,y)$. Substituting these to~\eqref{eq:welovepizzaanddoublecounting}, we get
\begin{equation*}
\int_{R^+\times R^-}h_i(x)\Delta(x,y)\diff\mu^2(x,y)=\int_{R^+\times R^-}h_i(y)\Delta(x,y)\diff\mu^2(x,y)\;.
\end{equation*}
We have $h_i(x)> r$, $h_i(y)\le r$, and further these terms $h_i(x)$ and $h_i(y)$ multiply a nonzero main term by~\eqref{eq:EverybodyKnows}. This is a contradiction.

\emph{Proof of~\ref{bonton:3}.} Let $H$ be the value of $h_i$ on $C$. The definition of $h_i$ via the conditional expectation gives $\int_X h_i\diff \mu=\int_X \One_{\Lambda_i}\diff\mu$. Together with $\int_X h_i\diff \mu=H\mu(C)$ and $\int_X \One_{\Lambda_i}\diff\mu=\mu(\Lambda_i)$, the claim follows.
\end{proof}

\begin{claimSecondB}
For a connected component $\Lambda_i$ of $U$ and the connected component $C$ of $\Delta$ such $\Lambda_i\subset C$ (up to a $\mu$-nullset), we have that $U\llbracket\Lambda_i\rrbracket$ is fractionally isomorphic to $\Delta\llbracket C\rrbracket$.
\end{claimSecondB}
\begin{proof}
We verify fractional isomorphism using Definition~\ref{def:FracIso}. Let $A\subset \P^c$ be measurable. By~\eqref{eq:transferOfDegreeMeasureIntoSubgraphs} we have
\begin{equation*}
\mu(\Lambda_i)\nu_{U\llbracket\Lambda_i\rrbracket}(A)=\mu\left(i_U^{-1}(A)\cap \Lambda_i\right)=\int_{i_U^{-1}(A)} \One_{\Lambda_i}(x)\diff\mu(x)\;.
\end{equation*}
The set $i_U^{-1}(A)$ is $\C(U)$-measurable by Lemma~\ref{lem:generate_CW}. Hence we have
\begin{align*}
    \mu(\Lambda_i)\nu_{U\llbracket\Lambda_i\rrbracket}(A)&=\int_{i_U^{-1}(A)} \One_{\Lambda_i}(x)\diff\mu(x)=\int_{i_U^{-1}(A)} h_i(x)\diff\mu(x)\\
    \JUSTIFY{by Claim~\ref*{thm:factorFI}.A}&=\frac{\mu(\Lambda_i)}{\mu(C)} \cdot \mu\left(i_U^{-1}(A)\cap C\right) 
\eqByRef{eq:transferOfDegreeMeasureIntoSubgraphs}
    \mu(\Lambda_i)\nu_{\Delta\llbracket C\rrbracket}(A)\;,
\end{align*}
as desired.
\end{proof}

We can now quickly finish the proof. Indeed, Claims~\ref*{thm:factorFI}.A and~\ref*{thm:factorFI}.B tell us that for the left-hand side of~\eqref{eq:factorFI} we have 
\begin{equation}\label{eq:LHS}
LHS\eqref{eq:factorFI}=\sum_{C:(*)}\mu(C)\;,
\end{equation}
where the summation runs over connected components $C$ of $\Delta$ such that $\Delta\llbracket C\rrbracket$ is fractionally isomorphic to $\Gamma$. By repeating the argument for $W$ and for $\Delta':=\Expectation[W \vert \C(W)\times\C(W) ]$, we see that
\begin{equation}\label{eq:RHS}
RHS\eqref{eq:factorFI}=\sum_{C:(**)}\mu(C)\;,
\end{equation}
where the summation runs over connected components $C$ of $\Delta'$ such that $\Delta'\llbracket C\rrbracket$ is fractionally isomorphic to $\Gamma$. We now use one of the main characterizations of fractional isomorphism. Namely, Theorem~1.2(3) in~\cite{Grebik2022} tells us that kernels $U/{\C (U)}$ and $W/{\C (W)}$ (introduced therein) are isomorphic.
By routinely expanding the definitions of $U/{\C (U)}$ and $W/{\C (W)}$, this implies that $\Delta$ and $\Delta'$ are weakly isomorphic. By Lemma~\ref{lem:weakIso}, there exist measure preserving maps $\pi$ and $\pi'$ such that $\Delta^{\pi}=(\Delta')^{\pi'}$ almost everywhere. We replace the host kernel $\Delta$ by $\Delta^{\pi}$ in~\eqref{eq:LHS}. Obviously, the sum, which we call $\sum_{C:(*mpm)}\mu(C)$ has not changed. Similarly, we write $\sum_{C:(**mpm)}\mu(C)$ has not for~\eqref{eq:RHS} in which $\Delta'$ is replaced by $(\Delta')^{\pi'}$. We have
\[
LHS\eqref{eq:factorFI}\eqByRef{eq:LHS}\sum_{C:(*)}\mu(C)
=
\sum_{C:(*mpm)}\mu(C)
\eqBy{L\ref{lem:weakIso}}
\sum_{C:(**mpm)}\mu(C)
=\sum_{C:(**)}\mu(C)\eqByRef{eq:RHS}RHS\eqref{eq:factorFI}\;,
\]
as was needed.

\section{Theorem~\ref{thm:mainUST}: generalization and proof}\label{sec:HNT}

In this section we provide a proof of Theorem~\ref{thm:mainUST}. In fact, we state and prove Theorem~\ref{thm:mainUSTGeneral}, which represents a generalization of Theorem~\ref{thm:mainUST} to possibly disconnected kernels.

In Section~\ref{ssec:USTGeneral} we introduce the necessary terminology and provide the statement of Theorem~\ref{thm:mainUSTGeneral}. We also discuss the more challenging $(\Rightarrow)$ direction and explain how it suffices to consider distributional information on the finite non-ancestral part $\fU^-_W$ of the branching process $\fU_W$, which is in turn closely linked to the Markov renormalization $W^\dag$. This reduces the key $(\Rightarrow)$ direction to the problem of proving Lemma~\ref{lem:FracIsoIFFProjFracIso}.

In Section~\ref{ssec:TheEasyDirectionUST} we provide a proof of the $(\Leftarrow)$ direction. In Section~\ref{ssec:ProofofLemmaFracIsoProjFracIso} we set up the necessary framework to provide a proof of Lemma~\ref{lem:FracIsoIFFProjFracIso}, which is given at the end of the subsection.
 
\subsection{Generalization of Theorem~\ref{thm:mainUST} to disconnected kernels} \label{ssec:USTGeneral}
Theorem~\ref{thm:mainUST} does not generalize verbatim to disconnected kernels. To see this, partition $X=X_1\sqcup X_2$ with $\mu(X_1)=0.2$ and $\mu(X_2)=0.8$. Let $W$ be constant-13 on $X_1\times X_1$, constant-7 on $X_2\times X_2$, and constant-0 otherwise. No matter whether $x\in X$ from Definition~\ref{def:HNTbranching} is sampled from $X_1$ or from $X_2$, the entire branching process $\fU_W$ will stay confined to that connected component, in which it will behave as in $U\equiv 1$. In line with this, we introduce below the notion of piecewise projective fractional isomorphism, which exactly captures the fact that the projective constants in Definition~\ref{def:projfraciso} may be different on different connected components.

Let $U$ and $W$ be nondegenerate bounded-degree kernels with ground spaces $X$ and $Y$ respectively. Let $X=\Lambda_1\sqcup \Lambda_2\sqcup\ldots$ and $Y=\Omega_1\sqcup \Omega_2\sqcup\ldots$ be decompositions into the connected components of $U$ and $W$ respectively. We say that $U$ and $W$ are \emph{piecewise projectively fractionally isomorphic} if for each connected bounded-degree kernel $\Gamma$ we have
\begin{equation}\label{eq:defPPFI}
\sum_{i:\text{$U_{\restriction \Lambda_i}$ is proj.\ frac.\ iso.\ to $\Gamma$}}\mu(\Lambda_i)
=
\sum_{i:\text{$W_{\restriction \Omega_i}$ is proj.\ frac.\ iso.\ to $\Gamma$}}\mu(\Omega_i)
\;.    
\end{equation}
Obviously, piecewise projective fractional isomorphism is an equivalence relation on nondegenerate bounded-degree kernels. We can now formulate the generalization of Theorem~\ref{thm:mainUST} to disconnected kernels.
\begin{theo}\label{thm:mainUSTGeneral}
Suppose that $U$ and $W$ are $L^\infty$-kernels with positive minimum degrees. Then $\fU_U$ and $\fU_W$ have the same distribution if and only if $U$ and $W$ are piecewise projectively fractionally isomorphic.
\end{theo} 
\begin{proof}
    
The $(\Leftarrow)$ direction, while somewhat tedious, follows by expanding the definitions. We postpone this verification to Section~\ref{ssec:TheEasyDirectionUST}.
Here, we focus on the $(\Rightarrow)$ direction.

Suppose that $\fU_U$ and $\fU_W$ have the same distribution. We first claim that $\fU^-_U$ and $\fU^-_W$ have the same distributions. To illustrate a difficulty and a way to circumvent it, let us consider a toy example from the realm of real numbers rather than branching processes. Suppose that $a,b,a^+,a^-,b^+,b^-\in \R$ are such that $a = a^+ + a^-$ and $b = b^+ + b^-$. Simply being told that $a=b$ alone is not adequate information to conclude that $a^-=b^-$. However, we do get this conclusion if we are also told that, say, $a^+,b^+\in\Z$ and $a^-,b^-\in[0,1)$. We can apply a similar idea of splitting into big and small to the branching processes $\fU_U$ and $\fU_W$ whose equality in distribution is given to us. Here, the `big parts' are the ancestral parts and the `small parts' are $\fU^-_U$ and $\fU^-_W$. The ancestral parts are infinite since they contain an infinite $\mathsf{anc}$-path. To complement this, Lemma~\ref{lem:extinct} tells us that $\fU^-_U$ and $\fU^-_W$ are finite almost surely (here, we also use Fact~\ref{fact:CorrespondendceDownUST} to transfer from $\fX_{U^\dag}$ and $\fX_{W^\dag}$ to $\fU^-_U$ and $\fU^-_W$). It indeed follows that $\fU^-_U$ and $\fU^-_W$ have the same distribution.

Now Theorem~\ref{thm:mainBJR} applies and tells us that $U^\dagger$ and $W^\dagger$ are fractionally isomorphic. The last piece in the puzzle is the following key lemma.

\begin{lem}\label{lem:FracIsoIFFProjFracIso}
Suppose that $U$ and $W$ are $L^\infty$-kernels with positive minimum degrees such that $U^\dag$ and $W^\dag$ are fractionally isomorphic. Then $U$ and $W$ are piecewise projectively fractionally isomorphic.
\end{lem}

We give a proof of Lemma~\ref{lem:FracIsoIFFProjFracIso} in Section~\ref{ssec:ProofofLemmaFracIsoProjFracIso}.
\end{proof}

\subsection{Proof of the $(\Leftarrow)$ direction of Theorem~\ref{thm:mainUSTGeneral}}\label{ssec:TheEasyDirectionUST}

Given a kernel $W$ with maximum degrees are bounded from above by $b<\infty$ and minimum degrees are bounded from below by $a>0$, we shall give a construction of $i_{W^\dagger}$ from $i_W$. Then, in Lemmas~\ref{lem:distributionofancestoralpaths} and~\ref{lem:frompathstoUSTs}, we conclude that $\fU_U$ and $\fU_W$ have indeed the same distribution.

First, we construct a space, denoted $\P^{a,b}$, which is a subset of $\P^b$. Define $\mathbb{D}^{(0)}_{a,b}:=\{\alpha\in\P^b:D(\alpha)< a\}$, and inductively, $\mathbb{D}^{(\ell)}_{a,b}:=\left\{\alpha\in\P^b\mid\mu_\alpha(\mathbb{D}^{(\ell-1)}_{a,b})>0\right\}$. Let \label{notation:Pab}$\P^{a,b}:=\P^b\setminus \bigcup_{\ell=0}^\infty \mathbb{D}^{(\ell)}_{a,b}$. So, $\P^{a,b}$ are elements of $\P^b$ whose degree is at least $a$, and which only see elements of degrees at least $a$, and which only see elements which only see elements of degrees at least $a$, \dots. Obviously, this property holds for elements of the support of $\nu_Z$ of any akernel $Z$ with $\degmax(Z)\le b$ and $\degmin(Z)\ge a$.

We  now introduce a notion of degree tilting. For each $\alpha\in\P^{a,b}$ we define \emph{degree tilting}\label{notation:degreetilting} $\alpha^{\ddagger}\in\P^{a/b,b/a}$ by $\alpha^{\ddagger}(0)=\OneElementSpace$ and inductively for $n\in\N$ and for each Borel set $A\subset \P_{n-1}^{a/b,b/a}$ by
\begin{equation}\label{eq:defdoubledagger}
\alpha^{\ddagger}(n)(A)= \int_{\{\gamma\in \P^{a,b} \mid \gamma^\ddagger(n-1)\in A\}}\frac {1}{D(\beta)}\diff\mu_{\alpha}(\beta)\;.    
\end{equation}
Note that this definition is not circular as we only need the $(n-1)$-st component of all $\gamma^\ddagger$ to define the $n$-th component of $\alpha^\ddagger$. It is a straightforward exercise, that indeed $\alpha^{\ddagger}\in\P^{a/b,b/a}$. 

We now return to the kernel $W$. For each $x\in X$, define \label{notation:lambda}$\lambda_{x,W}\in\M_{\le b}(X)$ by $\lambda_{x,W}(Y)=\int_Y W(x,z)\diff\mu(z)$.

\begin{lem} \label{lem:pushforwarddegmeas}
     For each $x\in X$, $\mu_{i_W(x)}$ is the pushforward of $\lambda_{x,W}$ via $i_W$.
\end{lem}
\begin{proof}
    As $\B(\P^b)$ is generated by $\bigcup_{n\in\N}\{p_{n,\infty}^{-1}(A)\mid A\in \B(\P_n^b)\}$, it suffices to show that $\mu_{i_W(x)}(p_{n,\infty}^{-1}(A))=(i_W)_*\lambda_{W,x}(p_{n,\infty}^{-1}(A))$ for each $n\in \N$ and $A\in \B(\P_n^b)$. We have
    \begin{align*}
        \mu_{i_W(x)}(p_{n,\infty}^{-1}(A))&=(p_{n,\infty})_*\mu_{i_W(x)}(A)=i_{W}(x)(n+1)(A)\\
        &=\int_{\{y\in X\mid p_{n,\infty}(i_W(y))\in A\}} W(x,z)\diff\mu(z)=\int_{i_W^{-1}(p_{n,\infty}^{-1}(A))} W(x,z)\diff\mu(z)\\
        &=\lambda_{W,x}(i_W^{-1}(p_{n,\infty}^{-1}(A)))=((i_W)_*\lambda_{W,x})(p_{n,\infty}^{-1}(A))\;,
    \end{align*} 
as was needed.
\end{proof}

\begin{lem} \label{lem:transformdegmeas}
    For each $x\in X$ we have $i_{W^\dagger}(x)=(i_W(x))^\ddagger$.
\end{lem}
\begin{proof}
    Recall that we have $a\le \degmin(W)$ and $b\ge \degmax(W)$. 
    We proceed by induction on $n\in \N$. For $n=0$ we have $i_{W^\dagger}(0)=\OneElementSpace=(i_W)^\ddagger(0)$.
    As an induction hypothesis assume that $i_{W^\dagger}(y)(k)=(i_W(y))^\ddagger(k)$ for each $y\in X$ and $k\in[n]_0$.
    
    Suppose we are given $A\in\B(\P_n^{b/a})$. Define $B:=\{\alpha\in \P^{a,b}\mid \alpha^\ddagger(n)\in A\}$. For each $x\in X$,
    \begin{align*}
        (i_{W^\dagger}(x)(n+1))(A)=&\int_{\{y\in X\mid i_{W^\dagger}(y)(n)\in A\}}W^\dagger(x,z)\diff\mu(z)\\
        \JUSTIFY{induction hypothesis, Definition~\ref{def:MarkovRenorm}}= &\int_{\{y\in X\mid i_{W}(y)^\ddagger(n)\in A\}}\frac{W(x,z)}{D(i_{W}(z))}\diff\mu(z)\\
        \JUSTIFY{definition of $B$}= &\int_{\{y\in X\mid i_{W}(y)\in B\}}\frac{W(x,z)}{D(i_{W}(z))}\diff\mu(z)\\
        \JUSTIFY{definition of $\lambda_{x,W}$}= &\int_{\{y\in X\mid i_{W}(y)\in B\}}\frac{1}{D(i_{W}(z))}\diff\lambda_{x,W}(z) \\
        \JUSTIFY{change of variables, Lemma~\ref{lem:pushforwarddegmeas}}= &\int_{B}\frac{1}{D(\beta)}\diff((i_W)_*\lambda_{x,W})(\beta) \\
        =&\int_{B}\frac{1}{D(\beta)}\diff\mu_{i_W(x)}(\beta)\\
        \JUSTIFY{by \eqref{eq:defdoubledagger}}=&\left(\left((i_{W}(x)\right)^\ddagger(n+1)\right)(A)\;,
    \end{align*} 
    as was needed.
\end{proof}

Given an akernel $W$ with minimum degree $a>0$ and maximum degree $b$, we introduce a stochastic process \label{notation:frakA}$\mathfrak{A}_W=(\mathfrak{A}^{(0)}_W,\mathfrak{A}^{(1)}_W,\ldots)\in (\P^{b/a})^{\N_0}$ as follows. First, generate $z_0$ from $X$ according to $\mu$. Given $z_{k-1}$, generate $z_k$ according to the intensity measure $\frac{W(z_{k-1},\cdot)}{\deg_W(z_{k-1})}$. Now, for each $k\in\N_0$, let $\mathfrak{A}^{(k)}_W:=\big(i_W(z_k)\big)^\ddagger$. That is, $(\mathfrak{A}^{(0)}_W,\mathfrak{A}^{(1)}_W,\ldots)$ has the distribution of the labels of the ancestral path in $\fU_W$ (as they move away from the root), when viewed through the map $\big(i_W(\cdot)\big)^\ddagger$.

The $(\Leftarrow)$ direction of Theorem~\ref{thm:mainUSTGeneral} follows from the following two lemmas.
\begin{lem}\label{lem:distributionofancestoralpaths}
Suppose that $U$ and $W$ are bounded-degree kernels with positive minimum degrees. If $U$ and $W$ are piecewise projectively fractionally isomorphic then $\mathfrak{A}_U$ and $\mathfrak{A}_W$ have the same distribution.
\end{lem}
\begin{lem}\label{lem:frompathstoUSTs}
Suppose that $U$ is a bounded-degree kernel with positive minimum degree. Then the distribution $\fU_U$ is determined by the distribution $\mathfrak{A}_U$.

In particular, suppose that $U$ and $W$ are bounded-degree kernels with positive minimum degrees with the property that $\mathfrak{A}_U$ and $\mathfrak{A}_W$ have the same distribution. Then $\fU_U$ and $\fU_W$ have the same distribution.
\end{lem}

\begin{proof}[Proof of Lemma~\ref{lem:distributionofancestoralpaths}]
Let $(\Gamma_j)_{j\in J}$ be a collection of connected mutually non-projectively fractionally isomorphic kernels for which~\eqref{eq:defPPFI} is positive, and let $f_j$ be the value of~\eqref{eq:defPPFI} for $\Gamma_j$. We have $\sum_{j\in J}f_j=1$. We couple the random choice of the root $z_0$ in $U$ and of the root $z'_0$ in $W$ as follows. First, we pick a random $j\in J$ from the distribution $(f_j)_{j\in J}$. Then, in $U$ we pick a random connected component $U_{\restriction \Lambda_i}$ that is projectively fractionally isomorphic to $\Gamma_j$. That is, $i$ is taken from the distribution 
\[\left(\frac{\mu(\Lambda_i)}{f_j}\right)_{i:\text{$U_{\restriction \Lambda_i}$ is proj.\ frac.\ iso.\ to $\Gamma_j$}}\;.\]
Similarly, in $W$ we pick a random connected component $W_{\restriction \Omega_{i'}}$ that is projectively fractionally isomorphic to $\Gamma_j$. Write $U^*:=U_{\restriction \Lambda_i}/\|U_{\restriction \Lambda_i}\|_1$, $W^*:=W_{\restriction \Omega_{i'}}/\|W_{\restriction \Omega_{i'}}\|_1$.
Last, we take $z_0\in \Lambda_i$ and $z'_0\in\Omega_{i'}$ uniformly at random but coupled in way so that 
$$i_{U^*}(z_0)
=i_{W^*}(z'_0)\;.
$$
This is possible because of projective fractional isomorphism. Note that $z_0$ and $z'_0$ are uniformly chosen in $X$, and hence can be used for generating the root particles for $\mathfrak{A}_U$ and $\mathfrak{A}_W$, respectively. Now, suppose that for $k=0,1,2,\ldots$ we have $z_k\in \Lambda_i$ and $z'_k\in\Omega_{i'}$ such that
\begin{equation}\label{eq:indas}
i_{U^*}(z_k)
=i_{W^*}(z'_k)\;.
\end{equation}
For extending $\mathfrak{A}_U$ and $\mathfrak{A}_W$, we should generate $z_{k+1}$ and $z'_{k+1}$ according to $\frac{\lambda_{z_k,U}(\cdot)}{\deg_U(z_k)}$ and $\frac{\lambda_{z'_k,W}(\cdot)}{\deg_W(z'_k)}$. By~\eqref{eq:indas} and Lemma~\ref{lem:pushforwarddegmeas} we can couple the choice of $z_{k+1}$ and $z'_{k+1}$ in a way that we have again~\eqref{eq:indas}. 

To conclude, we have generated sequences $z_0,z_1,z_2,\ldots$ and $z'_0,z'_1,z'_2,\ldots$ so that we have~\eqref{eq:indas} for each $k$. In particular, that means that we have $(i_{U^*}(z_k))^\ddagger
=(i_{W^*}(z'_k))^\ddagger$. By the virtue of Lemma~\ref{lem:transformdegmeas}, we have $i_{(U^*)^\dagger}(z_k)
=i_{(W^*)^\dagger}(z'_k)$. Definition~\ref{def:MarkovRenorm} tells us that $\dagger$ is invariant under multiplying a kernel by any nonzero scalar. Also, note that for an akernel $Z$ and an element $x$ of its ground space, $i_Z(x)$ depends only on the component of $Z$ containing $x$. These facts combined give
\[i_{U^\dag}(z_k)=i_{U\llbracket\Lambda_i\rrbracket^\dag}(z_k) =i_{W\llbracket\Omega_{i'}\rrbracket^\dag}(z'_k)=i_{W^\dagger}(z'_k)\;.\]
Using Lemma~\ref{lem:transformdegmeas} again, we get $(i_U(z_k))^\ddagger=(i_W(z'_k))^\ddagger$, as was needed.
\end{proof}

\begin{proof}[Proof of Lemma~\ref{lem:frompathstoUSTs}]
The ancestral path in $\fU_U$ is generated in the same way as we generate the sequence $z_0,z_1,\ldots$ in the definition of $\mathfrak{A}_U$. That is, if we view each type $(\mathsf{anc},\cdot)$ through the lens of the operator $\big(i_U(\cdot)\big)^{\ddagger}$, we get the distribution of $\mathfrak{A}_U$. By Lemma~\ref{lem:transformdegmeas}, this is equivalent to viewing the types $(\mathsf{anc},\cdot)$ on the ancestral path through the lens of $i_{U^\dagger}(\cdot)$. Definition~\ref{def:HNTbranching} tells us that we can construct $\fU_U$ from by attaching an independent branching process $\fX_{U^\dagger}(\cdot)$ to each vertex of type $(\mathsf{anc},\cdot)$ on its ancestral path. By Proposition~\ref{prop:CorrespondenceTwoBranchings}\ref{corespond1}, this is the same as attaching independent branching processes $\fB(i_{U^\dagger}(\cdot))=\fB((i_{U}(\cdot))^\ddagger)$ (c.f. Lemma~\ref{lem:transformdegmeas}). That is, the distribution $\mathfrak{A}_U$ fully determines the distribution $\fU_U$.
\end{proof}

\subsection{Proof of Lemma~\ref{lem:FracIsoIFFProjFracIso}}\label{ssec:ProofofLemmaFracIsoProjFracIso}

We begin with a definition. Given an $L^\infty$-kernel $W$ on $(X,\B,\mu)$ with positive minimum degree and a decomposition $X = \bigsqcup_{i\in I}\Lambda_i$ into connected components, let $W^{\heartsuit}$ be given by
\begin{equation}\label{not:heart}
    W^{\heartsuit}(x,y) = 
    \begin{cases}
         \frac{W(x,y)\mu(\Lambda_i)}{\int_{\Lambda_i\times \Lambda_i}W(u,v)\diff\mu(u)\diff\mu(v)} &\textrm{ if $x,y\in \Lambda_i$ for some $i\in I$,}\\
         0 &\textrm{ otherwise.}
    \end{cases}
\end{equation}
It is easy to verify that $W^{\heartsuit}$ is an $L^\infty$-kernel with positive minimum degree.

Our proof of Lemma~\ref{lem:FracIsoIFFProjFracIso} has two main steps. First, we show that there is a transformation from $i_{W^\dag}$ to $i_{W^{\heartsuit}}$ that depends only on the fractional isomorphism class of $W^\dag$. This is the subject of Lemma~\ref{lem:fraciso-degmeastrans}, and a substantial part of this section is devoted to building tools for its proof. Lemma~\ref{lem:fraciso-degmeastrans} allows us to show that $W^{\heartsuit}$ and $U^{\heartsuit}$ are fractionally isomorphic in the setting of Lemma~\ref{lem:FracIsoIFFProjFracIso}. A quick application of Theorem~\ref{thm:factorFI} then allows us to conclude that $U$ and $W$ are piecewise projectively fractionally isomorphic.

For our proof, we need to borrow two results from~\cite{Grebik2022}, namely Proposition~\ref{prop:didm-transform} and Lemma~\ref{lem:char-didm-transform}. These are stated for $L^\infty$-akernels, whereas the corresponding statements in~\cite{Grebik2022} are stated for akernels bounded by~1. The reduction is in both cases straightforward: Given an $L^\infty$-akernel $W$ in our setting, we consider an akernel $W':=\frac{W}{c}$, where $c:=\|W\|_\infty$. We now apply the corresponding result from~\cite{Grebik2022} to $W'$ and subsequently interpret it back for the akernel $W'$. Observe that there is an obvious transformation of the measure $\nu_W$ which is supported on $\P^c$ to $\nu_{W'}$ which is supported on $\P^1$. 

Since the objects $i_{W^\dag}(x)$ are elements of $\P^c$, it turns out to be useful to have for each $L^\infty$-akernel $W$ on $X$ a corresponding $L^\infty$-akernel $K[\nu_W]$ on $\P^c$. The following proposition combines generalizations of parts of Proposition~6.8 and Claim~6.9 in~\cite{Grebik2022}.

\begin{prop} \label{prop:didm-transform}
    Suppose that $\nu$ is a probability measure on $\P^c$ such that $\nu=\nu_W$ for some $L^\infty$-akernel $W$. Then $\mu_\alpha$ is absolutely continuous with respect to $\nu$ for $\nu$-almost every $\alpha\in\P^c$, with the Radon--Nikodym derivative satisfying $0 \le \frac{\diff\mu_{\alpha}}{\diff\nu} \le \|W\|_{\infty}$, and there is an akernel \label{notation:Knu}$K[\nu]\in L^{\infty}(\P^c\times\P^c,\nu\times\nu)$ such that $K[\nu](\alpha,-) = \frac{\diff\mu_{\alpha}}{\diff\nu}$ for $\nu$-almost every $\alpha\in\P^c$.
\end{prop}

The following lemma is a straightforward generalization of~\cite[Theorem 6.10]{Grebik2022} which formalizes the link between $W$ and $K[\nu_W]$. Given an $L^\infty$-akernel $W$ and a relatively complete sub-sigma-algebra $\C\in\Theta_\mu$, write $W_{\C}$ for $\Expectation[W\vert\B\times\C]$.
\begin{lem} \label{lem:char-didm-transform}
    For every $L^\infty$-akernel $W$ we have
    \begin{equation} \label{eq:itedegmeas-kernel}
        W_{\C(W)}(x,y) = K[\nu_W](i_W(x),i_W(y))
    \end{equation}
    for $(\mu\times\mu)$-almost every $(x,y)\in X\times X$.
\end{lem}

The function $\deg_{W^{\heartsuit}}$ plays a key role since it relates $W^\dag$ and $W^{\heartsuit}$. More specifically, we have 
\begin{equation}\label{eq:identitka}
W^\dag(x,y)\deg_{W^{\heartsuit}}(y) = W^{\heartsuit}(x,y)
\end{equation}
for each $x,y\in X$. For an $L^\infty$-kernel $W$ on $(X,\B,\mu)$ let \label{notation:chiW}$\chi_W$ be the measure on $(X,\B)$ given by 
\begin{equation}\label{eq:defChiW}
\diff\chi_W = \deg_{W^{\heartsuit}}\diff\mu
\end{equation}
and \label{notation:piW}$\pi_W = (i_{W^\dag})_*\chi_W$ be the pushforward measure on $\P^c$ ($c=\|W\|_\infty$) of $\chi_W$ via $i_{W^\dag}$. Note that $\chi_W$ and subsequently also $\pi_W$ are probability measures. The following lemma, which is the key intermediate step in our proof of Lemma~\ref{lem:fraciso-degmeastrans}, says that the measures $\pi_W$ and $\pi_U$ are equal whenever the Markov renormalizations $W^\dag$ and $U^\dag$ are fractionally isomorphic. This is the one place where we apply the theory of Markov chains.

\begin{lem} \label{lem:degmeasmatch}
    Suppose that $W$ and $U$ are $L^\infty$-kernels with positive minimum degree such that $W^{\dagger}$ and $U^{\dagger}$ are fractionally isomorphic. Write $c$ for $\|W^\dag\|_\infty$. Then the measures $\pi_W$ and $\pi_U$ on $\P^c$ are equal.
\end{lem}

\begin{proof}
    Write $\nu$ for $\nu_{W^\dagger}=\nu_{U^\dagger}$. Proposition~\ref{prop:didm-transform} returns an akernel $K[\nu]$.
    
    Before we dive into the proof, let us briefly motivate the steps in our proof. The basic idea is that we would like to run a Markov chain on $\P^c$ with the transpose of $K[\nu]$ as its transition probability kernel, apply Theorem~\ref{thm:MCconverge} to show that the Markov chain converges, and demonstrate that the unique limit probability measure is equal to both $\pi_W$ and $\pi_U$. The first significant hurdle we encounter is reducibility, that is, the existence of multiple connected components. We deal with this by splitting the space $\P^c$ into the connected components of $K[\nu]$, running a Markov chain in each of them, and then suitably amalgamating the Markov chains. The second obstacle is unsynchronized periodicity amongst countably infinitely many connected components. To handle this, we first analyze the behavior for any finite number $M$ of connected components and then utilize the fact that the contribution of the remaining infinite tail necessarily vanishes as $M$ tends to infinity.

    We begin by laying the groundwork for connected components in $K[\nu]$. Let $\F$ be the collection of Borel sets $B\subset\P^c$ such that
    \begin{equation} \label{eq:inseparable}
        \int_B\int_{\P^c\sm B} K[\nu](\alpha,\beta) \diff\nu(\alpha)\nu(\beta) = 0\;.
    \end{equation}
    This formula strongly resembles the formula in Definition~\ref{def:connected}\ref{en:connectedSep} defining separation between connected components, and indeed the intended meaning is the same. But we must be more careful as $K[\nu]$ is an akernel; indeed, the notion of connected components is more problematic in the asymmetric case, just as it is in oriented graphs. Here the connection between $K[\nu]$ and a symmetric kernel $W$ via the Markov renormalization $W^\dag$ as given in~\eqref{eq:itedegmeas-kernel} lets us draw upon structural information coming from connected components in $W$. In particular, we utilize this connection and the fact that $\nu=\nu_{W^\dag}$ to show in~\ref{item:complement-empty} below that~\eqref{eq:inseparable} is equivalent to the transposed version~\eqref{eq:transposed-bipartite-empty} given below. This symmetry will allow us to introduce the notion of `connected components of~$K[\nu]$'.
    
   	\begin{claim:lem:degmeasmatch:A}
   	The following hold for all $B\in\F$.
    \begin{enumerate}[label=(A\arabic{*})]
        \item \label{item:preimage-unioncomp} The set $i_{W^\dag}^{-1}(B)$ is a (possibly empty) union of connected components of $W$, modulo a $\mu$-null set.
        \item \label{item:complement-empty} We have $\P^c\sm B\in\F$.
    \end{enumerate}
    \end{claim:lem:degmeasmatch:A}
    \begin{proof}
    Write $S$ for $i_{W^\dag}^{-1}(B)$. Since $i_{W^\dag}$ is $\C(W^\dag)$-measurable by Lemma~\ref{lem:generate_CW}, by~\eqref{eq:itedegmeas-kernel} and~\eqref{eq:inseparable} we have
    \begin{align} \label{eq:bipartite-empty}
        \int_S\int_{X\sm S} W^\dag(x,y) \diff\mu(x)\mu(y) = \int_S\int_{X\sm S} W^\dag_{\C(W^\dag)}(x,y) \diff\mu(x)\mu(y)  =\int_B\int_{\P^c\sm B} K[\nu](\alpha,\beta) \diff\nu(\alpha)\nu(\beta) = 0\;.
    \end{align}
    Together with the definition of $W^\dag$, this means that $W$ is zero $(\mu\times\mu)$-almost everywhere on $S \times (X\sm S)$. Hence, \ref{item:preimage-unioncomp} follows. Now $W$ is symmetric, so both $W$ and $W^\dag$ are zero $(\mu\times\mu)$-almost everywhere on $(X\sm S) \times S$. By applying~\eqref{eq:itedegmeas-kernel} we obtain
    \begin{equation}\label{eq:transposed-bipartite-empty}
        \int_{\P^c\sm B}\int_B K[\nu](\alpha,\beta) \diff\nu(\alpha)\nu(\beta) = \int_{X\sm S}\int_S W^\dag_{\C(W^\dag)}(x,y) \diff\mu(x)\mu(y) = \int_{X\sm S}\int_S W^\dag(x,y) \diff\mu(x)\mu(y) = 0\;,
    \end{equation}
    so~\ref{item:complement-empty} follows.
    \end{proof}
    
    We say that two connected components of $W$ are \emph{inseparable} if for every $B\in\F$ we have that $i_{W^\dag}^{-1}(B)$, modulo a $\mu$-null set, contains both connected components or neither of them. Inseparability induces an equivalence relation on the collection of connected components of $W$; let $J$ be the index set for the countably many equivalence classes. Without loss of generality, we have $J=\N$ or $J=[N]$ for some $N\in\N$. 

    Let $k\in J$. For each $j\in J\sm\{k\}$ pick $A_{k,j}\in\F$ such that $i_{W^\dag}^{-1}(A_{k,j})$, modulo a $\mu$-null set, contains all of the connected components in the $k$th equivalence class and none of those in the $j$th equivalence class; this is possible by the definition of inseparability and the properties~\ref{item:preimage-unioncomp} and~\ref{item:complement-empty}. Let $S_k := \bigcap_{j\in J\sm\{k\}}A_{k,j}$ and \label{notation:Xk}$X_k := i_{W^\dag}^{-1}(S_k)$. We have $\nu(S_k) = \mu(X_k)$, and it is immediate from~\eqref{eq:inseparable} that $S_k\in\F$. Note that $S_k$ serves as the $k$th connected component of $K[\nu]$ (c.f. discussion after~\eqref{eq:inseparable}) and $X_k$ is the union of all the connected components in the $k$th equivalence class, modulo a $\mu$-null set.
    
    We now prepare to run a Markov chain on $S_k$. Write $\nu_k$ and $\mu_k$ for the probability measures $\frac{\nu_{\restriction S_k}}{\nu(S_k)}$ and $\frac{\mu_{\restriction X_k}}{\mu(X_k)}$ respectively. Now write \label{notation:Wk}$W_k$ for $W_{\restriction X_k}$. Here $W_k$ is an $L^\infty$-kernel on $(X_k\times X_k,\mu_k\times\mu_k)$ with positive minimum degree. This allows us to define akernels $W_k^\dagger:=(W_k)^\dagger$ and $W_k^\heartsuit:=(W_k)^\heartsuit$. Let us quickly derive some basic properties of these akernels. By Lemma~\ref{lem:MarkovBoundedDegree}, $W_k^\dagger$ is an $L^\infty$-akernel on $(X_k\times X_k,\mu_k\times\mu_k)$ with positive minimum degree. For $x,y\in X_k$, we have
    \begin{equation}\label{eq:relateXkdagger}
(W_k)^\dagger(x,y)=\frac{W_k(x,y)}{\deg_{W_k}(y)}=\frac{W(x,y)}{\mu(X_k)^{-1}\cdot\deg_{W}(y)}=\mu(X_k)\cdot W^\dagger(x,y)\;.
    \end{equation}

    We now define
    \label{notation:Knuk}$K[\nu]_k(\alpha,\beta) := (K[\nu])\llbracket S_k\rrbracket(\alpha,\beta)$. Since $W_k^\dagger$ is an $L^\infty$-akernel with positive minimum degree, we also get that $K[\nu]_k$ is an $L^\infty$-akernel on $(S_k\times S_k,\nu_k\times\nu_k)$ with positive minimum degree, say $\delta>0$, by Lemma~\ref{lem:char-didm-transform} and Lemma~\ref{lem:generate_CW}. By modifying $K[\nu]_k$ on a $(\nu_k\times\nu_k)$-null set, we may assume without loss of generality that $K[\nu]_k$ satisfies $\deg_{K[\nu]_k}(\alpha) \ge \delta$ for all $\alpha\in S_k$. Let $\Phi_k$ be a Markov chain on $S_k$ with transition probability kernel $Q_k(\alpha,A) = \int_A K[\nu]_k(\beta,\alpha)\diff\nu_k(\beta)$ for $A\in\B(S_k)$. We prepare to apply Theorem~\ref{thm:MCconverge} to $\Phi_k$.
    
    \begin{claim:lem:degmeasmatch:B}
    	The Markov chain $\Phi_k$ is $\nu_k$-irreducible.
    \end{claim:lem:degmeasmatch:B}
    \begin{proof}
For $A\in\B(S_k)$ and $\alpha\in S_k$ write $q_{\alpha,A}$ for the probability that $\Phi_k$ started at $\alpha$ satisfies $\Phi_k(n)\in A$ for some $n\in\N$. For the sake of contradiction, assume that there exists $A_0\in\B(S_k)$ with $\nu_k(A_0)>0$ such that $q_{\alpha,A_0} = 0$ for some $\alpha\in S_k$. Let $R_{A_0}$ be the set of all such $\alpha\in S_k$. Because $K[\nu]_k$ has positive minimum degree and bounded $L^\infty$-norm, we have $\nu_k(R_{A_0}), \nu_k(S_k\sm R_{A_0})>0$. Writing $Q_{k,\alpha}$ for the measure on $\B(S_k)$ given by $Q_{k,\alpha}(A) = Q_k(\alpha,A)$, for $\alpha\in R_{A_0}$ we have 
\begin{equation} \label{eq:OneStepReturnTimeDIDM}
        0=q_{\alpha,A_0} = Q_{k,\alpha}(A_0) + \int_{S_k\sm(A_0\cup R_{A_0})}q_{\beta,A_0}\diff Q_{k,\alpha}(\beta) \;.
\end{equation}    
Since $q_{\beta,A_0}>0$ for $\beta\in S_k\sm(A_0\cup R_{A_0})$, we get $Q_k(\alpha,S_k\sm R_{A_0}) = 0$ from~\eqref{eq:OneStepReturnTimeDIDM}. Hence, we have $R_{A_0}\in\F$. But now $i_{W^\dag}^{-1}(R_{A_0})$ and $i_{W^\dag}^{-1}(S_k\sm R_{A_0})$ are unions of connected components of $W$ with $R_{A_0},S_k\sm R_{A_0}\in\F$ by~\ref{item:preimage-unioncomp} and~\ref{item:complement-empty}, which contradicts the fact that $S_k\in\F$ represents an equivalence class of inseparable connected components in $W$.
\end{proof}
    
Let $\pi_k$ be the pushforward measure of $\chi_{W_k}$ via $(i_{W^\dag})_{\restriction X_k}$. We claim that $\pi_k$ is an invariant probability measure on $\B(S_k)$. Indeed, for every $A\in\B(S_k)$ we have
\begin{align} 
\nonumber
\pi_k(A) = ((i_{W^\dag})_{\restriction X_k})_*\chi_{W_k}(A) & = \int_{i_{W^\dag}^{-1}(A)} \deg_{W_k^{\heartsuit}}(x) \diff\mu_k(x) \\
\nonumber
& = \int_{i_{W^\dag}^{-1}(A)} \int_{X_k} W_k^\heartsuit(x,y)\diff\mu_k(y) \diff\mu_k(x) \\
\label{eq:Mar}
\JUSTIFY{by~\eqref{eq:identitka}}    & = \int_{X_k}\int_{i_{W^\dag}^{-1}(A)}  W_k^\dag(x,y) \diff\mu_k(x) \deg_{W_k^{\heartsuit}}(y) \diff\mu_k(y)\;.
\end{align}
Furthermore, for every $y\in X_k$ we have
\begin{equation}
\begin{split}\label{eq:Maru}
    \int_{i_{W^\dag}^{-1}(A)}  W_k^\dag(x,y) \diff\mu_k(x)&\eqByRef{eq:relateXkdagger}\int_{i_{W^\dag}^{-1}(A)}  \nu(S_k)W^\dag(x,y) \diff\mu_k(x)\\
    \JUSTIFY{$i_{W^\dag}$ is $\C(W^\dag)$-measurable by Lemma~\ref{lem:generate_CW}}&=\int_{i_{W^\dag}^{-1}(A)}  \nu(S_k)W^\dag_{\C(W^\dag)}(x,y) \diff\mu_k(x)\;.
\end{split}
\end{equation} 
Now by substituting both~\eqref{eq:Maru} and~\eqref{eq:defChiW} in~\eqref{eq:Mar}, we obtain
\[
\pi_k(A)
= \int_{X_k}\int_{i_{W^\dag}^{-1}(A)}  \nu(S_k)W^\dag_{\C(W^\dag)}(x,y) \diff\mu_k(x) \diff\chi_{W_k}(y)\;.
\]
Hence, we conclude that
\begin{equation*}
    \pi_k(A) = \int_{S_k} \int_A  K[\nu]_k(\beta,\alpha) \diff\nu_k(\beta) \diff\pi_k(\alpha)  = \int_{S_k} Q_k(\alpha,A)\diff\pi_k(\alpha)
\end{equation*}
This shows that $\pi_k$ is invariant.

We define a full irreducibility probability measure $\psi_k$ for $\Phi_k$ by
\[\psi_k(A) := \int_{S_k} \sum_{j=1}^\infty 2^{-j-1}Q_k^j(\alpha,A) \diff\nu_k(\alpha),\quad A\in\B(S_k)\;.\]
(Recall that in order for $\psi_k$ to satisfy the definition of a full irreducibility probability measure from Section~\ref{ssec:RWbackground}, we need $\nu_k$ to be an irreducibility probability measure. This is indeed the case by Claim~\ref*{lem:degmeasmatch}.B.)

Hence, $\Phi_k$ is a positive recurrent Markov chain. Let $p_k$ be the period of $\Phi_k$. Since $W_k$ has positive minimum degree, the measures $\chi_{W_k}$ and $\mu_k$ are equivalent. Thus, the measures $\pi_k$ and $\nu_k$ are equivalent. For each $n\in\N_0$ we define a measure $\pi_{k,n}$ on $\B(S_k)$ given by $\pi_{k,n}(A) = \int_{S_k} \frac{1}{p_k} \sum_{r=0}^{p_k-1} Q_k^{np_k+r}(\alpha,A) \diff\nu_k(\alpha)$. Here $\pi_{k,n}$ represents the distribution of the Markov chain $\Phi_k$ started (at time $0$) from the uniform distribution on $S_k$ and averaged over the time interval $\{np_k,\dots,(n+1)p_k-1\}$, which is one full period of $\Phi_k$. In particular, this has the exact form needed to apply Theorem~\ref{thm:MCconverge}, which we do to obtain
\begin{equation} \label{eq:comp-conv-norm-didm}
    \lim_{n\to\infty} \| \pi_{k,n} - \pi_k\|_{TV} = 0\;.
\end{equation}

We now extend~\eqref{eq:comp-conv-norm-didm} from one component of $K[\nu]$ to multiple components. For this, we would like to define a quantity that resembles $\pi_{k,n}$ in that it represents a distribution averaged over whole multiples of periods. For multiple components it is natural to take the least common multiple of the periods of the components, so in particular we start with finitely many components. Let $M\in J$. Write $S_{(M)}$ and $X_{(M)}$ for $\bigcup_{k\in[M]}S_k$ and $\bigcup_{k\in[M]}X_k$ respectively. Let $\nu_{(M)}$ and $\mu_{(M)}$ denote the probability measures $\frac{\nu_{\restriction S_{(M)}}}{\nu(S_{(M)})}$ and $\frac{\mu_{\restriction X_{(M)}}}{\mu(X_{(M)})}$ respectively. Let $W_{(M)}(x,y) := W_{\restriction X_{(M)}}(x,y)$ and $K[\nu]_{(M)}(\alpha,\beta) := \nu(S_{(M)})K[\nu]_{\restriction S_{(M)}}(\alpha,\beta)$ be akernels on $(X_{(M)}\times X_{(M)},\mu_{(M)}\times\mu_{(M)})$ and $(S_{(M)}\times S_{(M)},\nu_{(M)}\times\nu_{(M)})$ respectively. Let $p_{(M)}$ be the least common multiple of $\{p_i:i\in[M]\}$. Let $Q_{(M)}(\alpha,A) = \int_A K[\nu]_{(M)}(\beta,\alpha)\diff\nu_{(M)}(\beta)$ define a transition probability kernel. By unpacking the definitions of $Q_{(M)}$, $Q_k$, $\nu_{(M)}$ and $\nu_k$, recalling the definition of the $n$-step transition probability kernel and applying the defining property of $\P^c\sm S_k\in\F$ given in~\eqref{eq:inseparable}, we get
\begin{equation} \label{eq:multcomp-transitions}
    Q_{(M)}^n(\alpha,A) = Q_k^n(\alpha,A\cap S_k)
\end{equation}
for all $\alpha\in S_k$. In other words, the Markov chain $\Phi_{(M)}$ with transition probability kernel $Q_{(M)}$ acting on the disjoint union $S_{(M)}$ of connected components has exactly the same behavior as $\Phi_k$ on the connected component $S_k$, and this behavior is independent of that on all other connected components. Let $\pi_{(M)}$ be the pushforward measure of $\chi_{W_{(M)}}$ via $(i_{W^\dag})_{\restriction X_{(M)}}$. By unpacking the definitions of $\chi_{W_{(M)}}$ and $\chi_{W_k}$, we obtain
\begin{equation} \label{eq:MCfinite-distrdecomp}
    \pi_{(M)}(A) = \sum_{k\in[M]} \tfrac{\nu(S_k)}{\nu(S_{(M)})}\cdot\pi_k(A\cap S_k)\;.
\end{equation}
This is consistent with the expected limit behavior of the Markov chain $\Phi_{(M)}$; it says that $\pi_{(M)}$ is a convex combination of the functions $\pi_k$, with the factors $\tfrac{\nu(S_k)}{\nu(S_{(M)})}$ representing the relative mass in each connected component. For each $n\in\N_0$ let $\pi_{(M),n}$ be a measure on $\B(S_{(M)})$ given by $\pi_{(M),n}(A) = \int_{S_{(M)}} \frac{1}{p_{(M)}} \sum_{r=0}^{p_{(M)}-1} Q_{(M)}^{np_{(M)}+r}(\alpha,A) \diff\nu_{(M)}(\alpha)$. By applying~\eqref{eq:multcomp-transitions}, we get
\begin{equation} \label{eq:MCfinite-limitdecomp}
    \pi_{(M),n}(A) = \sum_{k\in[M]} \tfrac{\nu(S_k)}{\nu(S_{(M)})}\cdot\tfrac{p_k}{p_{(M)}} \sum_{\ell=1}^{p_{(M)}/p_k} \pi_{k,n+\ell-1}(A\cap S_k)\;.
\end{equation}
Again, this is consistent with the described behavior of the Markov chain $\Phi_{(M)}$; the factors $\tfrac{p_k}{p_{(M)}}$ and inner sum adjust for the different period lengths.

Now by applying~\eqref{eq:MCfinite-distrdecomp}, \eqref{eq:MCfinite-limitdecomp} and the definition of the total variation norm, we get
\[\| \pi_{(M),n} - \pi_{(M)}\|_{TV} \le \sum_{k\in[M]} \tfrac{\nu(S_k)}{\nu(S_{(M)})}\cdot\tfrac{p_k}{p_{(M)}} \sum_{\ell=1}^{p_{(M)}/p_k} \| \pi_{k,n+\ell-1} - \pi_k\|_{TV}\;.\]
Hence, by~\eqref{eq:comp-conv-norm-didm} we have 
\begin{equation} \label{eq:partsum-comp-conv-norm-didm}
    \lim_{n\to\infty} \| \pi_{(M),n} - \pi_{(M)}\|_{TV} = 0\;.
\end{equation}

The argument above deals with the case when $J$ is finite (by picking $M\in\N$ so that $J=[M]$), so it remains to consider the case $J=\N$. In this case, we would also like to run a Markov chain $\Phi$ on the entire space. Without loss of generality, we may assume that $k\mapsto\mu(X_k)$ is non-increasing. Let $Q(\alpha,A) = \int_A K[\nu](\beta,\alpha)\diff\nu(\beta)$ define a transition probability kernel. As before, we obtain
\begin{equation} \label{eq:countablecomp-transitions}
    Q^n(\alpha,A) = Q_k^n(\alpha,A\cap S_k)
\end{equation}
for all $\alpha\in S_k$ by unpacking the definitions of $Q$, $Q_k$, $\nu$ and $\nu_k$, recalling the definition of the $n$-step transition probability kernel and applying the defining property of $\P^c\sm S_k\in\F$ given in~\eqref{eq:inseparable}. Furthermore, by unpacking the definitions of $\chi_W$ and $\chi_{W_k}$ in analogy with~\eqref{eq:MCfinite-distrdecomp} and recalling that $\nu(\P^c) = 1$, we get
\begin{equation} \label{eq:MCfull-distrdecomp}
    \pi_W(A) = \sum_{k\in J} \nu(S_k)\cdot\pi_k(A\cap S_k)\;.
\end{equation}
In particular, note that the behavior of $\Phi$ on any finite subset $[M]\subseteq J$ of connected components is precisely what is described by our previous analysis for $\Phi_{(M)}$ (up to rescaling for relative mass).

The key obstacle to treating all components simultaneously is that we cannot generally hope to have a least common multiple of all (possible infinitely many distinct) period lengths.
This makes it challenging to define the analogue of $\pi_{(M),n}$ towards a version of~\eqref{eq:partsum-comp-conv-norm-didm}; here the analogue of $\pi_{(M)}$ is $\pi_W$ as defined in~\eqref{eq:MCfull-distrdecomp}. 
To handle this difficulty, we shall do the following. First, we define a measure $\rho_{M,n}$ on $\B(\P^c)$ that differs from $\pi_{(M),n}$ only in that it represents the behavior of the Markov chain $\Phi$ on the whole space instead of just $S_{(M)}$ (with suitable rescaling for relative mass). For each $n\in\N_0$ and $M\in J$ let $\rho_{M,n}$ be a measure on $\B(\P^c)$ given by $\rho_{M,n}(A) = \int_{\P^c} \frac{1}{p_{(M)}} \sum_{r=0}^{p_{(M)}-1} Q^{np_{(M)}+r}(\alpha,A) \diff\nu(\alpha)$. Second, we pick sequences $(M_m)_{m\in\N_0}$ and $(N_m)_{m\in\N_0}$ of non-negative integers so that $\nu(\P^c\sm S_{(M_m)}) \le 2^{-m-2}$ and for all $n\ge N_m$ we have $\| \pi_{(M_m),n} - \pi_{(M_m)}\|_{TV} \le 2^{-m-1}$; these are well defined because $\nu(\P^c) = 1$ is finite, the sequence $(\nu(S_k))_{k\in\N}$ is monotone decreasing with limit $0$ and we have~\eqref{eq:partsum-comp-conv-norm-didm}. For each $m\in\N_0$ let $\widetilde{\pi}_m = \rho_{M_m,N_m}$.

By applying~\eqref{eq:countablecomp-transitions}, recalling~\eqref{eq:MCfinite-limitdecomp} and splitting the expression for $\widetilde{\pi}_m$ into a main part corresponding to the connected components indexed by $[M_m]$ and another part representing the residual tail, for each $m\in\N_0$ we obtain
\begin{equation} \label{eq:MCfull-limitdecomp}
    \widetilde{\pi}_m(A) = \nu(S_{(M_m)})\cdot\pi_{(M_m),N_m}(A\cap S_{(M_m)}) + \sum_{k\in J\sm[M_m]} \nu(S_k) \cdot \tfrac{1}{p_{(M_m)}} \sum_{r=0}^{p_{(M_m)}-1} a_{m,k,r}(A\cap S_k) \;,
\end{equation}
where each $a_{m,k,r}$ is a probability measure on $S_k$ defined by
\begin{equation} \label{eq:MCfull-unsyncedtail}
    a_{m,k,r}(B) = \int_{S_k} Q_k^{N_mp_{(M_m)}+r}(\alpha,B) \diff\nu_k(\alpha)\quad\mbox{for each Borel set $B\subseteq S_k$}\;.
\end{equation}
Observe that we can also split the sum in~\eqref{eq:MCfull-distrdecomp} into a main part corresponding to the connected components indexed by $[M_m]$ and another part representing the residual tail as
\begin{equation} \label{eq:MCfull-distrdecomp-split}
    \pi_W(A) = \nu(S_{(M_m)})\cdot\pi_{(M_m)}(A\cap S_{(M_m)}) + \sum_{k\in J\sm[M_m]} \nu(S_k)\cdot\pi_k(A\cap S_k)\;.
\end{equation}

We shall compute a bound on $\|\widetilde{\pi}_m - \pi_W\|_{TV}$ by comparing the terms of the sums in~\eqref{eq:MCfull-limitdecomp} and~\eqref{eq:MCfull-distrdecomp-split}. Applying the definition of total variation norm in~\eqref{eq:TVNorm} and countable additivity of measures, we get
\begin{equation*}
    \|\widetilde{\pi}_m - \pi_W\|_{TV} \le \nu(S_{(M_m)})\cdot\left\|\pi_{(M_m),N_m} - \pi_{(M_m)}\right\|_{TV} + \sum_{k\in J\sm[M_m]} \nu(S_k) \cdot \left\|\tfrac{1}{p_{(M_m)}} \cdot \left(\sum_{r=0}^{p_{(M_m)}-1} a_{m,k,r}\right) - \pi_k\right\|_{TV} \;.
\end{equation*}
Since the measures $\pi_k$ and $a_{m,k,r}$ (given by~\eqref{eq:MCfull-unsyncedtail}) are all probability measures, we obtain
\[\|\widetilde{\pi}_m - \pi_W\|_{TV} \le \nu(S_{(M_m)})\cdot\left\|\pi_{(M_m),N_m} - \pi_{(M_m)}\right\|_{TV} + 2\nu\left(\P^c \sm S_{(M_m)}\right) \le 2^{-m}\;.\]
Hence, we have $\lim_{m\to\infty} \|\widetilde{\pi}_m - \pi_W\|_{TV} = 0$. 
    
An analogous argument for $U$ yields $\lim_{m\to\infty} \| \widetilde{\pi}_m - \pi_U\|_{TV} = 0$, so $\pi_W$ and $\pi_U$ are equal.
\end{proof}

The following lemma shows that the transformation from $i_{W^\dag}$ to $i_{W^{\heartsuit}}$ depends only on the fractional isomorphism class of the Markov renormalization $W^{\dag}$.

\begin{lem} \label{lem:fraciso-degmeastrans}
    Suppose that $W$ and $U$ are $L^\infty$-kernels with positive minimum degree such that $W^{\dag}$ and $U^{\dag}$ are fractionally isomorphic. Then there is a function $f\colon\P^c\to\P^c$ such that $i_{W^{\heartsuit}} = f\circ i_{W^\dag}$ and $i_{U^{\heartsuit}} = f\circ i_{U^\dag}$.
\end{lem}

\begin{proof}
    Write $\nu$ for $\nu_{W^\dag}=\nu_{U^\dag}$. Proposition~\ref{prop:didm-transform} returns an akernel $K[\nu]$. Write $\pi$ for $\pi_W = \pi_U$, which is well-defined by Lemma~\ref{lem:degmeasmatch}. We shall inductively construct a sequence $(f_k)_{k\in\N_0}$ of functions with the following properties. First, for all $k\in\N_0$ we have $f_k\colon\P^c\to\P^c_k$. Second, for $k\in\N$ and $j\in[k-1]_0$ we have $f_k(\alpha)(j) = f_{k-1}(\alpha)(j)$. Third, for all $k\in\N_0$ we have $i_{W^{\heartsuit},k} = f_k\circ i_{W^\dag}$ and $i_{U^{\heartsuit},k} = f_k\circ i_{U^\dag}$. Observe that the definitions of the functions $i_{W^{\heartsuit}}$ and $i_{U^{\heartsuit}}$ imply that the function $f\colon\P^c\to\P^c$ given by $f(\alpha)(k) = f_k(\alpha)(k)$ for $k\in\N_0$ has the desired properties.
    
    For $k=0$ let $f_0\colon\P^c\to\P^c_0$ be given by $f_0(\alpha) = \OneElementSpace$. This function has the desired properties because $i_{W^{\heartsuit},0}$ and $i_{U^{\heartsuit},0}$ are both constant functions mapping to $\OneElementSpace$. For $k\in\N$ we set $f_k(\alpha)(j) = f_{k-1}(\alpha)(j)$ for $j\in[k-1]_0$ and define $f_k(\alpha)(k)$ as follows. For measurable $A\subset\P^c_{k-1}$ set
    \[f_k(\alpha)(k)(A) = \int_{f_{k-1}^{-1}(A)} K[\nu](\alpha,\beta)\diff\pi(\beta)\;.\]
    We now verify that $f_k$ has the requisite properties. The first two are immediate, so we focus on the third one. Take $x\in X$. The domain of the functions $i_{W^{\heartsuit},k}$, $i_{U^{\heartsuit},k}$ and $f_k\circ i_{U^\dag}$ has components indexed by $[k]_0$. We first verify the claimed identity for components indexed by $j\in[k-1]_0$. This claim easily follows by induction. Indeed, we have
    \[(f_k\circ i_{W^\dag})(x)(j) = (f_{k-1}\circ i_{W^\dag})(x)(j) = i_{W^{\heartsuit},k-1}(x)(j) = i_{W^{\heartsuit},k}(x)(j)\;,\]
    as was needed. We now turn to the case $j=k$.
    For an arbitrary measurable $A\subset\P^c_{k-1}$ we have
    \begin{align*}
        f_k(i_{W^\dag}(x))(k)(A) &= \int_{f_{k-1}^{-1}(A)} K[\nu](i_{W^\dag}(x),\beta)\diff\pi(\beta)  = \int_{i_{W^\dag}^{-1}(f_{k-1}^{-1}(A))} W^\dag_{\C(W^\dag)}(x,y)\diff\chi_W(y) \\
        \JUSTIFY{\text{Lemma~\ref{lem:generate_CW},  $i_{W^{\heartsuit},k-1} = f_{k-1}\circ i_{W^\dag}$ and~\eqref{eq:defChiW}}}& = \int_{i_{W^{\heartsuit},k-1}^{-1}(A)} W^\dag(x,y)\deg_{W^{\heartsuit}}(y)\diff\mu(y)  \\
        \JUSTIFY{by~\eqref{eq:identitka}}&= \int_{i_{W^{\heartsuit},k-1}^{-1}(A)} W^\heartsuit(x,y)\diff\mu(y)=i_{W^{\heartsuit},k}(x)(k)(A)\;,
    \end{align*}
    as was needed.
    Hence, we have $i_{W^{\heartsuit},k} = f_k\circ i_{W^\dag}$. An analogous line of reasoning yields $i_{U^{\heartsuit},k} = f_k\circ i_{U^\dag}$, thereby completing the proof.
\end{proof}

We conclude the section with the proof of Lemma~\ref{lem:FracIsoIFFProjFracIso}.

\begin{proof}[Proof of Lemma~\ref{lem:FracIsoIFFProjFracIso}]
    Let $W$ and $U$ be $L^\infty$-kernels with positive minimum degree such that $W^{\dagger}$ and $U^{\dagger}$ are fractionally isomorphic. Lemma~\ref{lem:fraciso-degmeastrans} returns a function $f\colon\P^c\to\P^c$ such that $i_{W^{\heartsuit}} = f\circ i_{W^\dag}$ and $i_{U^{\heartsuit}} = f\circ i_{U^\dag}$. For an arbitrary measurable $A\subset\P^c$ we have
    \[\nu_{W^{\heartsuit}}(A) = \mu(i_{W^{\heartsuit}}^{-1}(A)) = \mu(i_{W^\dag}^{-1}(f^{-1}(A))) = \nu_{W^\dag}(f^{-1}(A))\]
    and analogously we have $\nu_{U^{\heartsuit}}(A) = \nu_{U^\dag}(f^{-1}(A))$. Since $W^{\dag}$ and $U^{\dag}$ are fractionally isomorphic, we have $\nu_{W^\dagger} = \nu_{U^\dagger}$ and so $\nu_{W^{\heartsuit}}(A) = \nu_{U^{\heartsuit}}(A)$. Hence, $W^{\heartsuit}$ and $U^{\heartsuit}$ are fractionally isomorphic. 

    Now fix a decomposition $X = \bigsqcup_{i\in I}\Lambda_i$ into the countably many connected components of $W$. Fix an analogous decomposition $X = \bigsqcup_{i\in I'}\Omega_i$ into the countably many connected components of $U$. In particular, by Theorem~\ref{thm:factorFI}, we have for every connected bounded degree kernel $\Gamma^*$ that
    \begin{equation*}
        \sum_{i\in I':\text{$U^{\heartsuit}\llbracket \Omega_i\rrbracket$ is frac.\ iso.\ to $\Gamma^*$}}\mu(\Omega_i) = \sum_{i\in I:\text{$W^{\heartsuit}\llbracket \Lambda_i\rrbracket$ is frac.\ iso.\ to $\Gamma^*$}}\mu(\Lambda_i) \;.    
    \end{equation*}
We can rewrite this by replacing restrictions of type $\llbracket\cdot\rrbracket$ with restrictions of type $\restriction$ to obtain
\begin{equation}\label{eq:RZVD}
\sum_{i\in I':\text{$\mu(\Omega_i)(U^{\heartsuit})_{\restriction \Omega_i}$ is frac.\ iso.\ to $\Gamma^*$}}\mu(\Omega_i)
=
\sum_{i\in I:\text{$\mu(\Lambda_i)(W^{\heartsuit})_{\restriction \Lambda_i}$ is frac.\ iso.\ to $\Gamma^*$}}\mu(\Lambda_i)
\;.    
\end{equation}
    
\begin{claim:lem:FracIsoIFFProjFracIso:A}
    The following are equivalent for each $i\in I$ and each connected bounded-degree kernel $\Gamma$.
    \begin{enumerate}[label=(\roman*)]
    \item\label{it:Bud1} $W_{\restriction \Lambda_i}$ is projectively fractionally isomorphic to $\Gamma$,
    \item\label{it:Bud2} $(W^\heartsuit)_{\restriction \Lambda_i}$ is projectively fractionally isomorphic to $\Gamma$, and
    \item\label{it:Bud3} $\mu(\Lambda_i)\|\Gamma \|_1\cdot (W^\heartsuit)_{\restriction \Lambda_i}$ is fractionally isomorphic to $\Gamma$.
    \end{enumerate}
\end{claim:lem:FracIsoIFFProjFracIso:A}

\begin{proof}
    The equivalence between~\ref{it:Bud1} and~\ref{it:Bud2} is obvious, since $W_{\restriction \Lambda_i}$ and $(W^\heartsuit)_{\restriction \Lambda_i}$ differ just by a scalar multiple. The direction~\ref{it:Bud3}$\Rightarrow$\ref{it:Bud2} is obvious. For the direction \ref{it:Bud2}$\Rightarrow$\ref{it:Bud3}, suppose that $(W^\heartsuit)_{\restriction \Lambda_i}$ is projectively fractionally isomorphic to $\Gamma$. That is, $c\cdot (W^\heartsuit)_{\restriction \Lambda_i}$ is fractionally isomorphic to $\Gamma$ for some $c>0$. It is our task to show that $c=\mu(\Lambda_i)\|\Gamma \|_1$. Since two fractionally isomorphic kernels must have the same $L^1$-norm, we have $c\cdot \|(W^\heartsuit)_{\restriction \Lambda_i}\|_1=\|\Gamma\|_1$. Now by applying~\eqref{not:heart} and noting that the restriction operation $\restriction \Lambda_i$ reweights the square $\Lambda_i\times \Lambda_i$ by a factor $\mu(\Lambda_i)^{-2}$, we obtain $\|(W^\heartsuit)_{\restriction \Lambda_i}\|_1=\frac1{\mu(\Lambda_i)}$, completing the proof.
\end{proof}

To finish the proof, we need to show that for every connected bounded-degree kernel $\Gamma$, we have 
\[
\sum_{i:\text{$U_{\restriction \Omega_i}$ is proj.\ frac.\ iso.\ to $\Gamma$}}\mu(\Omega_i)
=
\sum_{i:\text{$W_{\restriction \Lambda_i}$ is proj.\ frac.\ iso.\ to $\Gamma$}}\mu(\Lambda_i)
\;.    
\]
By Claim~\ref*{lem:FracIsoIFFProjFracIso}.A (and the obvious counterpart of this claim for the kernel $U$), this is equivalent to
\[
\sum_{i:\text{$\mu(\Omega_i)\|\Gamma \|_1\cdot (U^\heartsuit)_{\restriction \Omega_i}$ is frac.\ iso.\ to $\Gamma$}}\mu(\Lambda_i)
=
\sum_{i:\text{$\mu(\Lambda_i)\|\Gamma \|_1\cdot (W^\heartsuit)_{\restriction \Lambda_i}$ is frac.\ iso.\ to $\Gamma$}}\mu(\Lambda_i)
\;.    
\]
This last formula indeed holds by~\eqref{eq:RZVD} applied to $\Gamma^*:=\frac{\Gamma}{\|\Gamma\|_1}$.
\end{proof}

\section{Concluding remarks}\label{sec:concludingremarks}
\subsection{On boundedness}\label{ssec:boudnednessassumption}
Our first main result, Theorem~\ref{thm:mainBJR}, provides a characterization of akernels that yield the same distribution of the corresponding Bollob\'as--Janson--Riordan branching processes, under the assumption of bounded degrees. We do not know the extent to which this assumption may be relaxed. More specifically, the definition of the random rooted tree $\fX_W$ is sensible for every kernel $W$ for which $\deg_W(x)<\infty$ for almost every $x\in X$. The theory developed in~\cite{bollobas2007PhaseTransition,MR2659281} applies to all $L^1$-kernels.\footnote{It is a curious question which of the results from~\cite{bollobas2007PhaseTransition,MR2659281}  extend to kernels whose degrees are almost all finite.} However, our Theorem~\ref{thm:mainBJR} applies only to bounded degree akernels, a class somewhat broader than akernels of bounded $L^\infty$-norm (c.f.\ Footnote~\ref{foot:Example}). As we explained in Remark~\ref{rem:PcBounded}, the reason we were not able to extend Theorem~\ref{thm:mainBJR} to $L^1$-kernels is our reliance on the Stone--Weierstraß Theorem. We leave extending Theorem~\ref{thm:mainBJR} to $L^1$-kernels or, alternatively, finding a counterexample in the same setting as an open problem. Note that the $(\Leftarrow)$-direction as formulated in Proposition~\ref{prop:CorrespondenceTwoBranchings} extends easily.

Our second main result, Theorem~\ref{thm:mainUSTGeneral}, applies only to $L^\infty$-kernels of positive minimum degree, whereas the maximum potential generality seems to be nondegenerate $L^1$-kernels. There are three places where we rely on these degree assumptions. The first one is that we need the corresponding Markov renormalization to be of bounded degree (per Lemma~\ref{lem:MarkovBoundedDegree}) so that we are in the setting of Theorem~\ref{thm:mainBJR}. We discussed a possible extension of Theorem~\ref{thm:mainBJR} above. The second one is the extinction Lemma~\ref{lem:extinct}. It could be that Lemma~\ref{lem:extinct} could be extended to the needed setting quite easily. Again, the proof of the $(\Leftarrow)$ direction as formulated in Section~\ref{ssec:TheEasyDirectionUST} extends easily. The third place concerns the invariant sub-sigma-algebras introduced in~\cite{Grebik2022} and recalled in Section~\ref{ssec:MorePrelFI}. This theory was worked out in~\cite{Grebik2022} by Grebík and Rocha in the context of $L^2$-kernels. In a personal communication, Grebík suggested that the correct generality may be that of $L^1$-kernels.

\subsection{Random rooted trees versus branching processes}\label{ssec:ConclusionRandomRootedTreesVSBranchingProcess}
We could ask a version of our main questions in the setting of branching processes rather than random rooted tree, that is, we could ignore the tree structure and ask only about the total population in each generation. More formally, and in the setting of Theorem~\ref{thm:mainBJR}, if $(T,r)$ is a random rooted tree of $\fX_A$ for a given akernel $A$, then $\mathbf{X}_A$ is a stochastic process indexed by $\N_0$, and defined as
\[
(\mathbf{X}_A)_\ell=|\{v\in V(T):\mathrm{dist}_T(v,r)=\ell \}|
\quad\mbox{for $\ell\in\N_0$.}
\]
Obviously, the equality of distributions of $\fX_{A_1}$ and $\fX_{A_2}$ implies the equality of distributions of $\mathbf{X}_{A_1}$ and $\mathbf{X}_{A_2}$. It follows from Theorem~\ref{thm:mainBJR} that if $A_1$ and $A_2$ are fractionally isomorphic bounded-degree akernels, then $\mathbf{X}_{A_1}$ and $\mathbf{X}_{A_2}$ have the same distributions. We leave it as an open problem whether the converse holds as well. Likewise, we ask whether the branching process version of Theorem~\ref{thm:mainUSTGeneral} holds.

\subsection{Extremal questions}\label{ssec:extremal}
Many interesting extremal questions can be explored concerning $\fX_W$ and $\fU_W$. These problems are of independent interest, and become even more compelling in light of the connection between $\fX_W$ and the giant component in sparse inhomogeneous random graphs (Fact~\ref{fact:localstructureInhomo}), as well as the connection between $\fX_W$ and the uniform spanning tree (Theorem~\ref{thm:HNT}).

We are aware of only two existing nontrivial extremal results, one for $\fX_W$ and one for $\fU_W$. 

The first extremal result appears in~\cite{bollobas2007PhaseTransition} as Corollary~3.18. Before stating the result, let us give some explanatory background. It is a well-known fact in the theory of the giant component, known as the `barely supercritical regime', that the giant component in $\G(n,\frac{(1+\eps)}{n})$ for small $\eps>0$ is of order $2(\eps+o(\eps))n$. In~\cite{bollobas2007PhaseTransition}, it is shown that this gives the fastest possible growth rate of the giant component among inhomogeneous random graph models. More specifically, let $W\neq 0$ be an $L^2$-kernel, and let $s(W)>0$ be the supremum over $t$ with $\Pr[|\fX_{t W}|=\infty]=0$. That is, $s(W)\cdot W$ is the kernel at which the giant component starts to emerge. Then, the extremal result of~\cite{bollobas2007PhaseTransition} says that 
\[
\lim_{\eps\searrow 0}\frac{\Pr[|\fX_{(1+\eps)\cdot s(W)\cdot W}|=\infty]}{\eps}\le \lim_{\eps\searrow 0}\frac{\Pr[|\fX_{(1+\eps)}|=\infty]}{\eps}=2\;.
\]

The next extremal result about $\fU_W$ appeared in~\cite{MR3876899}. In particular, it is proven there, that for every graphon $W$ with $\degmin(W)>0$ we have
\begin{align*}
    \Probability[\deg_{\fU_W}(root)=1]&\ge \exp(-1)\;, \\
    \Probability[\deg_{\fU_W}(root)=2]&\le \exp(-1)\;\mbox{, and}\\
    \Probability[\deg_{\fU_W}(root)=k]&\le \exp(2-k)\cdot \frac{(k-2)^{k-2}}{(k-1)!}\;\mbox{for $k\ge 3$.}
\end{align*}
As we remarked earlier, these can be interpreted as a lower-bound on the number of leaves and upper-bound on the number of vertices of degree $k\ge 2$ in the uniform spanning tree of dense graphs.

Perhaps the most interesting extremal question is about maximizing the order of the giant component in an inhomogeneous random graph whose kernel has a prescribed $L^1$-norm.
\begin{quest}\label{quest:maxgiant}
Suppose that $a>0$. Among all $L^\infty$-kernels with $\|W\|_1=a$, find the supremum of $\Probability[|\fX_W|=\infty]$.
\end{quest}

The results and problems mentioned in this section are not linked to our main results directly. However, in solving these problems, Theorems~\ref{thm:mainBJR} and~\ref{thm:mainUST} would provide a tool to assert that any extremal kernel found for a specific problem is unique up to fractional isomorphism or up to projective fractional isomorphism.

We thank Matas Šileikis for allowing us to include this part; some work on these problems has been done in collaboration with him.

%% file: summaryofnotation.tex
\begin{longtable}{ll}
$(X,\B,\mu)$&ground space for all kernels and akernels, see p.\pageref{notation:XBmu}\\
$\deg_W(x)$& degree of an akernel at $x\in X$ (measured on the first coordinate), see p.\pageref{notation:degree}\\
$\degmin(W)$, $\degmax(W)$&minimum and maximum degree of $W$, see p.\pageref{notation:minmaxdegre}\\
$\fX_W$& Bollobás--Janson--Riordan branching process from akernel $W$, see p.\pageref{notation:fXW}\\
$\fU_W$& uniform spanning tree branching process from kernel $W$, see p.\pageref{def:HNTbranching}\\
$\fU_W^-$& descendant part of $\fU_W$, see p.\pageref{notation:HNTminus}\\
$\P_n^c$, $\P^c$&spaces of iterated degree measures, see p.\pageref{eq:Pfinite}, eq.\eqref{eq:Pfinite} and p.\pageref{eq:Pinfty}, eq.\eqref{eq:Pinfty}\\
$p_{\ell,k}$&projection from $\P_k$ to $\P_\ell$, see p.\pageref{notation:projection}\\
$\B(\P^c)$, $\B(\P_n^c)$&Borel sets on $\P^c$ and $\P^c_n$, see p.\pageref{notation:BorelP}\\
$D(\alpha)$& `total mass' of $\alpha\in\P^c$, see p.\pageref{notation:degreeMeasure}\\
$\mu_\alpha$& lifting of an element $\alpha\in\P^c$ to an element in $\M_{\le c}(\P^c)$, see p.\pageref{notation:lifting}\\
$i_{W,n}(x)$, $i_W(x)$ & degree information at vertex $x$, see p.\pageref{notation:i}\\
$\nu_W$&iterated degree measure, see p.\pageref{notation:iterateddegreemeasure}\\
$\fB(\alpha)$&branching process from iterated degree measure $\alpha\in \P^c$, see p.\pageref{notation:fBallpha}\\
$\fB_W$&branching process $\fB(\alpha)$, where $\alpha$ is sampled from akernel $W$ using $i_W$, see p.\pageref{notation:fBW}\\
$T^\uparrow$&introducing a new root preceding the old root, see p.\pageref{notation:uparrow}\\
$T_1\oplus T_2$&pasting two rooted trees along the root, see p.\pageref{notation:oplus}\\
$\T_n$&trees of height at most $n$, see p.\pageref{notation:Tn}\\
$\widehat{\nu}$ & exponential tilting of measure $\nu$, see p.\pageref{notation:exponentialtilting}\\
$W^\dagger$& Markov renormalization of a kernel, see p.\pageref{notation:Markovrenormalization}\\
$e_F$&coefficient associated to a tree, see p.\pageref{eq:constE}, eq.\eqref{eq:constE}\\
$f_{T,n}$, $\overline{f_{T,n}}$ & tree functions, see p.\pageref{def_treefunc}, Definition~\ref{def_treefunc}\\
$\F_n$, $\overline{\F_n}$& collections of functions $f_{T,n}$ and $\overline{f_{T,n}}$, see p.\pageref{notation:cF}\\
$\C^W_n$, $\C^W$& pullback sigma-algebras $\P^c_{n}$ via $i_{W,n}$ or of $\P^c$ via $i_W$, see p.\pageref{not:CW}\\
$U_{\restriction Y}$&restriction of a kernel, see p.\pageref{notation:restrictedgraphon}\\
$U\llbracket Y\rrbracket$&restricted and rescaled kernel, see p.\pageref{notation:restrictedrescaledgraphon}\\
$\P^{a,b}$&measures in $\P^b$ `separated from anything in $\P^a$', see p.\pageref{notation:Pab}\\
$\alpha^\ddagger$&degree tilting, see p.\pageref{notation:degreetilting}\\
$\lambda_{W,x}$&neighborhood measure of $x$ in $W$, see p.\pageref{notation:lambda}\\
$\mathfrak{A}_W$&$\P$-valued stochastic process arising from a random walk on $W$, see p.\pageref{notation:frakA}\\
$W^\heartsuit$&a renormalization inside components of $W$, see p.\pageref{not:heart}\\
$K[\nu]$ & the simplest akernel whose iterated degree measure is $\nu$, see p.\pageref{notation:Knu}\\
$\chi_W$&probability measure on $X$, $\diff\chi_W = \deg_{W^{\heartsuit}}\diff\mu$, see p.\pageref{notation:chiW}\\
$\pi_W$&probability measure on $\P^c$, $\pi_W = (i_{W^\dag})_*\chi_W$, see p.\pageref{notation:piW}\\
$X_k$&collection of `inseparable components' of $W$, see p.\pageref{notation:Xk}\\
$W_k,W_k^\dagger,W_k^\heartsuit$ & $W_{\restriction X_k}$ and versions thereof, see p.\pageref{notation:Wk}\\
$K[\nu]_k$& $=(K[\nu])\llbracket S_k\rrbracket$, see p.\pageref{notation:Knuk}
\end{longtable}

%% file: main.bbl
\begin{thebibliography}{10}

\bibitem{MR1085326}
D.~Aldous.
\newblock The continuum random tree. {I}.
\newblock {\em Ann. Probab.}, 19(1):1--28, 1991.

\bibitem{MR1885388}
N.~Alon and J.~H. Spencer.
\newblock {\em The probabilistic method}.
\newblock Wiley-Interscience Series in Discrete Mathematics and Optimization.
  Wiley-Interscience [John Wiley \& Sons], New York, second edition, 2000.
\newblock With an appendix on the life and work of Paul Erd\H{o}s.

\bibitem{GlobalUSTGraphon}
E.~Archer and M.~Shalev.
\newblock The {GHP} scaling limit of uniform spanning trees of dense graphs.
\newblock {\em Random Structures \& Algorithms}, 65(1):149--190, 2024.

\bibitem{MR1700749}
P.~Billingsley.
\newblock {\em Convergence of probability measures}.
\newblock Wiley Series in Probability and Statistics: Probability and
  Statistics. John Wiley \& Sons, Inc., New York, second edition, 1999.
\newblock A Wiley-Interscience Publication.

\bibitem{MR2599196}
B.~Bollob\'{a}s, C.~Borgs, J.~Chayes, and O.~Riordan.
\newblock Percolation on dense graph sequences.
\newblock {\em Ann. Probab.}, 38(1):150--183, 2010.

\bibitem{MR2659281}
B.~Bollob\'{a}s, S.~Janson, and O.~Riordan.
\newblock The cut metric, random graphs, and branching processes.
\newblock {\em J. Stat. Phys.}, 140(2):289--335, 2010.

\bibitem{bollobas2007PhaseTransition}
B.~Bollobás, S.~Janson, and O.~Riordan.
\newblock The phase transition in inhomogeneous random graphs.
\newblock {\em Random Structures \& Algorithms}, 31(1):3--122, 2007.

\bibitem{GraphLimitsBSLSV}
C.~Borgs, J.~T. Chayes, L.~Lov\'{a}sz, V.~T. S\'{o}s, and K.~Vesztergombi.
\newblock Convergent sequences of dense graphs. {I}. {S}ubgraph frequencies,
  metric properties and testing.
\newblock {\em Adv. Math.}, 219(6):1801--1851, 2008.

\bibitem{MR3829971}
H.~Dell, M.~Grohe, and G.~Rattan.
\newblock Lov\'{a}sz meets {W}eisfeiler and {L}eman.
\newblock In {\em 45th {I}nternational {C}olloquium on {A}utomata, {L}anguages,
  and {P}rogramming}, volume 107 of {\em LIPIcs. Leibniz Int. Proc. Inform.},
  pages Art. No. 40, 14. Schloss Dagstuhl. Leibniz-Zent. Inform., Wadern, 2018.

\bibitem{MR2668548}
Z.~Dvo\v{r}\'{a}k.
\newblock On recognizing graphs by numbers of homomorphisms.
\newblock {\em J. Graph Theory}, 64(4):330--342, 2010.

\bibitem{MR0125031}
P.~Erd\H{o}s and A.~R\'{e}nyi.
\newblock On the evolution of random graphs.
\newblock {\em Magyar Tud. Akad. Mat. Kutat\'{o} Int. K\"{o}zl.}, 5:17--61,
  1960.

\bibitem{FRIEZE198547}
A.~M. Frieze.
\newblock On the value of a random minimum spanning tree problem.
\newblock {\em Discrete Applied Mathematics}, 10(1):47--56, 1985.

\bibitem{Grebik2022}
J.~Greb{\'i}k and I.~Rocha.
\newblock Fractional isomorphism of graphons.
\newblock {\em Combinatorica}, 42(3):365--404, Jun 2022.

\bibitem{Gri:RandomTree}
G.~R. Grimmett.
\newblock Random labelled trees and their branching networks.
\newblock {\em J. Austral. Math. Soc. Ser. A}, 30(2):229--237, 1980/81.

\bibitem{MR3876899}
J.~Hladk\'{y}, A.~Nachmias, and T.~Tran.
\newblock The local limit of the uniform spanning tree on dense graphs.
\newblock {\em J. Stat. Phys.}, 173(3-4):502--545, 2018.

\bibitem{InhomoMST}
J.~Hladký and G.~Viswanathan.
\newblock Random minimum spanning tree and dense graph limits.
\newblock arXiv:2310.11705.

\bibitem{ConnectednessGraphons}
S.~Janson.
\newblock Connectedness in graph limits.
\newblock arXiv:0802.3795.

\bibitem{MR2908619}
S.~Janson.
\newblock Simply generated trees, conditioned {G}alton-{W}atson trees, random
  allocations and condensation.
\newblock {\em Probab. Surv.}, 9:103--252, 2012.

\bibitem{MR2816939}
S.~Janson and O.~Riordan.
\newblock Duality in inhomogeneous random graphs, and the cut metric.
\newblock {\em Random Structures Algorithms}, 39(3):399--411, 2011.

\bibitem{MR1068492}
R.~M. Karp.
\newblock The transitive closure of a random digraph.
\newblock {\em Random Structures Algorithms}, 1(1):73--93, 1990.

\bibitem{Kolchin1977}
V.~F. Kolchin.
\newblock Branching processes, random trees, and a generalized scheme of
  arrangements of particles.
\newblock {\em Mathematical notes of the Academy of Sciences of the USSR},
  21(5):386--394, May 1977.

\bibitem{MR3012035}
L.~Lov\'asz.
\newblock {\em Large networks and graph limits}, volume~60 of {\em American
  Mathematical Society Colloquium Publications}.
\newblock American Mathematical Society, Providence, RI, 2012.

\bibitem{GraphLimitsLSz}
L.~Lov{\'a}sz and B.~Szegedy.
\newblock {Limits of dense graph sequences}.
\newblock {\em J. Combin. Theory Ser. B}, 96(6):933--957, 2006.

\bibitem{MeynTweedieBook}
S.~P. Meyn and R.~L. Tweedie.
\newblock {\em Markov chains and stochastic stability}.
\newblock Communications and Control Engineering Series. Springer-Verlag
  London, Ltd., London, 1993.

\bibitem{MR4334549}
J.~Petit, R.~Lambiotte, and T.~Carletti.
\newblock Random walks on dense graphs and graphons.
\newblock {\em SIAM J. Appl. Math.}, 81(6):2323--2345, 2021.

\bibitem{MR1297385}
M.~V. Ramana, E.~R. Scheinerman, and D.~Ullman.
\newblock Fractional isomorphism of graphs.
\newblock {\em Discrete Math.}, 132(1-3):247--265, 1994.

\bibitem{rudin1976principles}
W.~Rudin et~al.
\newblock {\em Principles of mathematical analysis}, volume~3.
\newblock McGraw-hill New York, 1976.

\bibitem{MR0843938}
G.~Tinhofer.
\newblock Graph isomorphism and theorems of {B}irkhoff type.
\newblock {\em Computing}, 36(4):285--300, 1986.

\end{thebibliography}
